\documentstyle[12pt]{report} 

\def\lll{\label}

\def \Box {\!\sqcup \mkern -12mu \sqcap\mkern .6mu}

\def \dfll {\leaders \hbox to 1em {\hss.\hss}\hfill}
\def\beq {\begin{equation}}
\def\eeq {\end{equation}}
\newtheorem{theorem}{Theorem}
\newtheorem{cor}[theorem]{Corollary}

\newtheorem{defi}{Definition}
\newtheorem{qes}{Question}
\newtheorem{lemma}{Lemma}

\def\blm{\begin{lemma}}
\def\elm{\end{lemma}}
\def\bdf{\begin{defi}}
\def\edf{\end{defi}}
\def\btm{\begin{theorem}}
\def\etm{\end{theorem}}
\def\bpp{\begin{propos}}
\def\epp{\end{propos}}
\def\bQ {\begin{qes}}
\def\eQ {\end{qes}}
\def\btm{\begin{theorem}}
\def\etm{\end{theorem}}
\def\ben{\begin{enumerate}}
\def\een{\end{enumerate}}
\def\bar{\begin{array}}
\def\ear{\end{array}}

\def \cn {^{conn}}
\def \nhopf {\mbox{$^0\bigcirc\mkern-13mu\bigcirc\,$}}
\def \ep {\hfill$\Box$\ms}
\def \Rm {\mbox{${\tt R}\!-\!mod$}}

\def\hookfill{$\mathsurround=0pt \mathord \lhook \mkern-4mu 
\cleaders\hbox{$\mkern-2mu
\mathord- \mkern-2mu$}\hfill \mkern-6mu\mathord\to$}
\newcommand {\into} [1]{\hbox to #1pt{\hookfill}}
\newcommand {\INTO} [2]{\stackrel{\mbox{#1}}{\into {#2}}}

\def\twoheadrightarrow{\to\mkern-20mu\to}
\def\thrafill{$\mathsurround=0pt \mathord- \mkern-6mu 
\cleaders\hbox{$\mkern-2mu
\mathord- \mkern-2mu$}\hfill \mkern-6mu\mathord\twoheadrightarrow$}
\newcommand {\onto} [1]{\hbox to #1pt{\thrafill}}
\newcommand {\ONTO} [2]{\stackrel{\mbox{#1}}{\onto {#2}}}
\newcommand{\TO}[2]{\stackrel {\mbox{#1}}{\hbox to #2pt{\rightarrowfill}}}

\newcommand {\prg}[2]{\vspace*{.6cm}\noindent
{\bf #1)\ #2\ :\hfill\ }\lll{pg-#1}\newline\smallskip\vspace*{-.3cm}}
\newcommand {\cht}[2]{\setcounter{chapter}{#1}
\subsection*{#1)\  #2}\lll{pg-#1}}
\newcommand {\cit}[2]{ {\em\qquad #1)\ #2}\dfll\pageref{pg-#1}}

\def \be {\beta^{in\!t}}
\def \pint {\pi_1^{in\!t}}
\def \pc {^{\bullet}}
\def \ppc {^{\bullet\bullet}}

\def\ub{\underline}
\def \R {\mbox{\tt R}}

\def \H {\mbox{$H_1^{in\!t}$}}
\def \F {I\mkern -6.2mu  F}
\def \P {I\mkern -6.2mu  P}
\def \V {\mbox{\boldmath$\cal V$\unboldmath}}
\def \S {\mbox{${\cal S}_{\!l}$}}

\def\id{{1\mkern-5mu {\rm I}}}
\def\isto {\widetilde{\longrightarrow}}

\def \lz  {\langle}
\def \rz  {\rangle}

\def\B{{\cal B}}

\def\ms{\medskip}

\def\Ii {\mbox{\raise .4 ex\hbox{$\int$}$\!\! I$}}
\def \bm {\boldmath}
\def \ubm {\unboldmath}
\def \ct {\mbox{\bm${\cal T}\! gl($\ubm}}
\def \T {\mbox{\bm${\cal T}$\ubm}}

\def \ctinf {\mbox{\bm${\cal T}\! gl^{\infty}($\ubm}}
\def \G{\mbox{\boldmath$\Gamma$\unboldmath}}

\newcommand {\Cb} [1] {\mbox{$Cob_3(#1)$}} 
\newcommand {\Cc} [1] {\mbox{$Cob_3^{conn}(#1)$}}

\def \btimes {{\,\raise -.28 ex\hbox{$\Box$}\mkern -13.8mu \times\,}}

\begin{document} 

\vspace*{0cm}

\begin{center}

\section*{On the Connectivity  of Cobordisms and  \\
  Half-Projective TQFT's \\ }
\bigskip

\medskip

{\large THOMAS KERLER}\\
 \vspace*{0cm}

\bigskip

February, 1996

\end{center}
\vspace*{1.4cm}

{\small \noindent{\bf Abstract} : We consider a generalization of 
the axioms of a TQFT, so called half-projective TQFT's, where we
inserted an anomaly, ${\tt x}^{\mu_0}\,$, in the composition law.
Here $\mu_0$ is a coboundary (in a group cohomological sense)
on the cobordism categories 
with  non-negative, integer values. The element {\tt x} of the ring 
over which the TQFT is defined does not have to be 
invertible. In particular, it may be zero.

This modification makes it possible to extend quantum-invariants,
which vanish on $S^1\times S^2$, to  non-trivial  TQFT's. Note, that a TQFT in
the ordinary sense of Atiyah with this property has to be trivial
all together. We organize our discussions such that the notion of a
half-projective TQFT is extracted as the only possible generalization under 
a few very natural assumptions.

Based on separate work with Lyubashenko on connected TQFT's, we construct
a large class of half-projective TQFT's with ${\tt x}=0$. Their invariants
all vanish on $S^1\times S^2$, and they coincide with the Hennings invariant
for non-semisimple Hopf algebras and, more generally, with the 
Lyubashenko invariant for non-semisimple categories. 

We also develop a few topological tools that allow us to 
determine the cocycle $\mu_0$ and
find  numbers, $\varrho(M)\,$, such that the  linear map associated to 
a cobordism, $M$, is of the form ${\tt x}^{\varrho(M)}f_M\,$. 
They are concerned with connectivity properties
of cobordisms, as for example maximal non-separating surfaces.
We introduce in particular the notions of 
``interior'' homotopy and homology groups,
and of coordinate graphs, which are functions on cobordisms
with values in the morphisms of a graph category.

For  applications we will prove that half-projective TQFT's with
${\tt x}=0$ vanish on cobordisms with infinite interior  homology,
and we argue that the order of
divergence of the TQFT on a cobordism, $M\,$, in the ``classical limit''
can be estimated by the rank of its maximal free interior group, which 
coincides with $\varrho(M)\,$.}
\bigskip

\subsection*{Contents}\

{\bf 1.) Introduction}\hfill\pageref{pg-1}

\cit{1.1} {Survey of Contents} 

\ms

{\bf 2.) Cobordism Categories, and Half-Projective TQFT's}\hfill\pageref{pg-2}

\cit{2.1}{ Categories of Cobordisms, the Structure of \Cb * }

\cit{2.2}{ Elementary Compositions, and the Cocycle $\mu_0$ }

\cit{2.3}{ Half-Projective TQFT's, and Generalizations}

\cit{2.4}{ TQFT's for Cobordisms with Corners}

\ms

{\bf 3.) Non-separating Surfaces, Interior Fundamental Groups, \newline
\hphantom{xxxxxxx} and Coordinate-Graphs}\hfill\pageref{pg-3}

\cit{3.1} { r-Diagrams of Non-Separating Surfaces}

\cit{3.2}{ Interior Fundamental Groups, and an A-Priori Estimate on $\varrho(M)$}

\cit{3.3}{ The Graph-Category \G\  }

\cit{3.4}{ Coordinate Graphs of Cobordisms }

\cit{3.5}{ Existence of Coordinate Graphs from Interior Groups}

\ms

{\bf 4.) Construction of Half-Projective TQFT's}\hfill\pageref{pg-4}

\cit{4.1}{ Surface-Connecting Cobordisms}

\cit{4.2}{ Basic Constraints on Generalized TQFT's}

\cit{4.3}{ The Example of Extended TQFT's}

\cit{4.4}{ Integrals, Semisimplicity, and ${\tt x}=\V(S^1\times S^2)\,$}

\cit{4.5}{ Main Result, and Hints to Further Generalizations and Applications}

\ms

{\bf Appendix}\hfill\pageref{pg-A}

\cit{A.1}{ Proofs of Section 3.4}

\cit{A.2}{ The Spaces $\H(M,{\bf G})$, the Numbers $\be_j(M)\,$, and Further
Anomalies}

\cit{A.3}{ Summary of Tangle Presentations}

\ms

{\bf References}\hfill\pageref{pg-R}

\bigskip

\medskip

\cht 1 {Introduction}

Although physical motivations were at its origin, the notion of topological
quantum field theories (TQFT's) has become a part of algebraic topology,
since it was axiomatically defined by Atiyah [A]. In the same way as 
 for example homology
or homotopy, it is given as a functor from a topological category into an
algebraic category. More precisely, it is a functor of the following 
symmetric tensor categories:
\beq\lll{eq-I0}
\V\,:\;Cob_{d+1}\,\longrightarrow\,{\tt R}-mod\;,
\eeq
where $Cob_{d+1}$ has as objects d-manifolds and as morphisms $d+1$-dimensional
cobordisms between them, and ${\tt R}-mod$ is the usual category of (free)
{\tt R}-modules and {\tt R}-linear maps. Their tensor products are the 
disjoint union $\sqcup$ and $\otimes_{\tt R}$, respectively.
(Quite often we only require {\tt R} to be a ring).

We may think of \V\ as a representation of the algebra of cell-attachments
to the boundaries of  $d+1$-manifolds. But unlike, e.g., homology it detects
algebraically much more involved relations between the cells than just
their intersection numbers. In this paper we shall be exclusively concerned
with the case $d=2$, where this algebra corresponds to quantum groups. 
We will make extensive reference to known TQFT's in 2+1 dimensions, and
use results of three-dimensional, geometric  topology.
\medskip

The original purpose of this paper is to resolve a paradox that occurs
in several different examples of ``quantum-invariants'' of 3-manifolds
and concrete quantum field theories. It is about  a degeneracy that 
at first sight  appears to prevent us from
constructing a TQFT in the rigorous, axiomatic sense. The problem is resolved
by a seemingly minor modification of the axioms of a TQFT, yielding what
we shall call {\em half-projective} TQFT's. The generalization naturally
leads us to several questions about the connectivity of cobordisms,
for which we will develop several tools that should have applications 
also to other topological problems.
\medskip

Specifically, the phenomena that we are interested in is that 
sometimes the invariant of a ``quantum-theory'' vanishes on the product of
sphere and circle, i.e.,
\beq\lll{eq-I1}
{\tt x}\,:=\,\V(S^1\times S^2)\;=\;0\qquad.
\eeq
It is an elementary implication of the axioms of a TQFT, observed in [Wi]
but also [A],
that for a surface, $\Sigma$, the invariant of the circle product  is
the dimension of the associated vector space, i.e., 
\beq\lll{eq-I2}
\V(S^1\times\Sigma)\;=\;dim\Bigl(\V(\Sigma)\Bigr)\qquad.
\eeq
Hence (\ref{eq-I1}) entails also triviality of the vector space,
$\V(S^2)=0\,$. A dramatic consequences of this for a TQFT is that 
$\V\equiv 0$ on all surfaces and cobordisms. The reason is easily seen,
if we remove from a general cobordism $M:\Sigma_s\to\Sigma_t$ a ball
so that $M^*:=M-D^3:\,\Sigma_s\to\,\Sigma_t\sqcup S^2\,$. Expressing the 
regluing of $D^3$ as a composition of cobordisms, we obtain 
\beq\lll{eq-I3}
M\,\,:\;\;\Sigma_s\;\TO{$M^*$}{30}\;\Sigma_t\sqcup S^2\;\TO{$D^3$}{30}\;
\Sigma_t\qquad.
\eeq
An application of \V\ to this yields the assertion, since the middle
vector space is mapped to zero. Clearly, this means that a (non-trivial) 
invariant can be extended to a TQFT only if it does not vanish on 
$S^1\times S^2\,$. 
\medskip

However, the examples, in which the degeneracy of (\ref{eq-I1}) is 
encountered, do appear to have a lot of the structural properties of
a TQFT, and are very closely related to situations, where non-trivial
TQFT's actually
exist.

An example of a more algebraic nature is the Hennings invariant of 
3-manifolds, see [H], for a finite dimensional, quasi-triangular 
 Hopf algebra, $\cal A\,$.
Its construction is analogous to that of the invariant of
Reshetikhin-Turaev in [RT], also [T], except that special elements of 
$\cal A$ are used directly instead of the representation theory of $\cal A\,$.
 The two invariants can be put on the same footing [Ke3] via the 
Lyubashenko invariant [L2], which is defined for abelian braided tensor
categories. The invariant of [RT] can be extended to a TQFT, which is usually
identified with the Chern Simons quantum field theory. Nevertheless, it
is easy to see that the (non-trivial) 
Hennings invariant vanishes on $S^1\times S^2$ if 
(and only if) $\cal A$ is not semisimple.

Subsequent studies in [L1], [L2], [KL], and [Ke3] showed that the
latter  invariants
can still be associated to representations of mapping class groups and,
more generally, TQFT's for cobordisms of {\em connected} surfaces. 
This reveals that the vanishing paradox in (\ref{eq-I1}) has to
have its origin in basic connectivity properties  of cobordisms.

The situation is analogous for the Kuperberg invariant [Ku] for a
Hopf algebra $\cal B$, which is believed to be closely related to the 
Hennings invariant. Again, if $\cal B$ is not semisimple, it vanishes
on $S^1\times S^2\,$, although it is constructed in a similar way as the
Turaev-Viro invariant [TV], for which we always have TQFT's.
\medskip

A second example in the concrete physical framework of
 quantum field theory is
given by conformal field theories and their
corresponding Chern Simons theories,
whose gauge groups are supergroups. The case of $U(1,1)$-WZW-models has
been worked out in [RS]. In this model we have the same vanishing for
 the ratios of partition functions 
$\V(S^1\times S^2):=\,{\cal Z}(S^1\times S^2)/
{\cal Z}(S^3)\,=0\,$. In the attempt to construct the operators of a TQFT
(in some regularization with parameter $\varepsilon\to 0$) one is faced
in [RS] with serious singularities of matrix elements in $\varepsilon$
that can only be removed at the price of having degenerate inner products 
for the physical state spaces $\V(\Sigma)\,$.
\medskip

We show in this article that there is essentially only one way to modify
the axioms of a TQFT, such that we preserve the tensor product rule for \V,
and  ensure that \V\ respects gluings of cobordisms over {\em connected } 
surfaces. The answer is 
provided by  the notion of a half-projective TQFT. By this we mean a
map between the category of cobordisms to the category of {\tt R}-modules,
which is a functor as in (\ref{eq-I0}), except that the composition law 
is of the form 
\beq\lll{eq-I4}
\V\bigl(M\circ N\bigr)\,\;=\;\,{\tt x}^{\mu_0(M,N)}\,\V(M)\V(N)\qquad.
\eeq
Here $\mu_0$ is a ``coboundary'' on $Cob_3$ in the sense of 
group cohomology, when we view categories as generalizations
of groupoids.
$\mu_0(M,N)\,$ can be   computed from basic
connectivity data of $M$ and $N\,$, and it has values only in the 
non-negative integers, $\mu_0\in{\bf Z}^{0,+}\,$.

If the number 
${\tt x}\in {\tt R}\,$ is invertible (e.g., ${\tt x}\neq 0$ and 
$k={\tt R}\,$ is a field) then the anomaly can of course be removed
by rescaling \V. Also, if $M$ and $N$ are composed over only one 
connected component, we find that $\mu_0(M,N)=0\,$ so that \V\ behaves
like an honest functor. This is consistent with the connected TQFT-functors
in [KL]. The identity (\ref{eq-I2}), however, is now modified. Repeating 
the original derivation with (\ref{eq-I4}) we find
\beq\lll{eq-I5}
\V(S^1\times\Sigma)\;=\;{\tt x}\,dim\Bigl(\V(\Sigma)\Bigr)\qquad.
\eeq
Hence the above examples are not in contradiction with extensions
to half-projective TQFT's, if we set ${\tt x}=0\,$.
\medskip

One of the main results of this paper is the construction of a large class
of non-trivial, half-projective TQFT's with ${\tt x}=0\,$. Our starting
point here are the connected TQFT's from [KL], but as in the
$U(1,1)$-WZW-models we have to deal with degenerate pairings of the
spaces $\V(\Sigma)\,$.
The  TQFT's we find extend, in particular, the Hennings invariant 
for an arbitrary
non-semisimple, finite-dimensional, modular Hopf algebra.
\medskip

In the quantum-algebraic framework, the element {\tt x} (and, especially,
the fact whether it is trivial or not) is intimately related to semisimplicity
 and cointegrals of the respective categories or  Hopf algebras.
In a concrete quantum field theory {\tt x} may be seen as a parameter
for the renormalization of the  product of field operators.
\medskip

It is interesting to observe that the two non-semisimple
examples from above share a few more common features beyond (\ref{eq-I1}).
In both cases we find that the representations of the mapping class group 
$SL(2,{\bf Z})\,$  of the torus contains algebraic summands and tensor factors
(the semisimple ones only produce finite representations), and that the
invariants of lens spaces and Seifert-manifolds
 are proportional to the order of the 
first homology group, see [RS] and [Ke3], and references therein.
We will investigate these properties in the general, axiomatic setting in
separate work.
\bigskip

\smallskip

\prg{1.1}{Survey of Contents}
\medskip

In Chapter 2 we provide the definitions of the cobordism category $\Cb *$
(Section 2.1), the special cocycle $\mu_0$ (Section 2.2), the notion of
half-projective TQFT's (Section 2.3), and various versions of extended
(half-projective) TQFT's (Section 2.4). In Section 2.2 we also
give an explicit formula (see Lemma~\ref{lm-mu0=dim}) and an algorithm
as in (\ref{eq-mu0-recurs}) for the computation  of $\mu_0\,$.
Section 2.3  includes a discussion of the basic implications of 
the axioms of a half-projective TQFT. In particular, it is shown in
Corollary~\ref{cor-altern} that, generically, the only properly 
half-projective, indecomposable  TQFT's are those with ${\tt x}=0\,$,
even for general rings, {\tt R}.
\medskip

The purpose of Chapter 3 is to develop the topological means that allow 
us to treat the connectivity properties of cobordisms relevant to the
formalism of half-projective TQFT's. In this we are mainly motivated by
the result in Lemma~\ref{lm-ord-rho}, which asserts that 
 $\V(M)\,=\,{\tt x}^{\varrho(M)}f_M\,$ for some ``regularized'' {\tt R}-linear
map $f_M\,$, where $\varrho(M)\,$ is given by the maximal number of 
non-separating surfaces in $M\,$. Note, that in general 
${\tt x}^{\varrho(M)}\,$ may generate  non-trivial ideals in the 
respective space of linear maps, seen as an {\tt R}-module.
In order to be able to compute the number $\varrho(M)\,$, 
we show in Theorem~\ref{thm-spc} of 
Section 3.5 that it 
is identical with the maximal rank of a
{\em free } interior group. 
 The notion of {\em interior fundamental groups},
where we divide by the subgroup coming from the surfaces,
 is introduced in Section 3.2 , and a basic gluing-property under compositions
of cobordisms is described in Lemma~\ref{lm-comp-freefac}.

In Section 3.4 we define {\em coordinate graphs} of manifolds with boundary, 
which is a rather useful tool to the end of encoding the connectivity 
properties  of a cobordism in a combinatorial way. Coordinate graphs
are given by 
(Morse) functions on  cobordisms with values in graphs that belong to
the graph category from Section 3.3. A result of particular interest is
Lemma~\ref{lm-decomp-ogr}, which ensures that to a decomposition of 
coordinate graphs we can always find a corresponding {\em connected}
decomposition of the cobordisms.

An interesting application of the results in Chapter 3 
is Corollary~\ref{cor-O},
which asserts that if \V\ is a general,
 half-projective TQFT with ${\tt x}=0\,$,
then $\V(M)=0\,$ for any cobordism, for which $\be_1(M)\neq 0\,$, i.e.,
with infinite, {\em interior} homology. For the special case of the 
$U_q(s\ell_2)$-Hennings invariant ($q$ a root of unity), evaluated
on closed manifolds, $M:\emptyset\to\emptyset\,$,  this vanishing phenomenon
was also observed by a direct calculation in [O].
\medskip

In Chapter 4 we show how non-trivial, half-projective TQFT's can be 
constructed from connected ones as, e.g., those
in [KL]. We start in Section 4.1
with the discussion of an algebra of special cobordisms between a surface
and the connected sum of its  components. Using the existence of
decompositions as in Lemma~\ref{lm-conmor-dec}, these cobordisms
 allow us to express
any cobordism by one that cobords only connected surfaces.

In Section 4.2 we start with the list of 
 Axioms V1-V5 for a generalized TQFT, \V, which
essentially state that \V\ respects tensor products aa well as
 compositions over
connected surfaces. We show that \V\ necessarily has to be a half-projective
TQFT. Moreover, we exhibit a list of eight properties that have to hold for a
connected TQFT, if it should extend to a disconnected one  satisfying V1-V5. 
The latter is true if it is, e.g.,  a specialization  from a 
half-projective one.
In Theorem~\ref{thm-conn-half} we show that these properties are in fact 
also sufficient, in order to guarantee the existence of such a 
 half-projective TQFT.

The purpose of Section 4.3 is to show that all but one of these properties are
automatically fulfilled, if the connected TQFT descends from an 
{\em extended} TQFT. In particular, we will use that \V\ is a factorization of
a functor $\V_1:\Cb 1\to{\cal C}\,$, where the objects of \Cb 1\ are surfaces
with one hole, and $\cal C$ is an abelian, braided tensor category.
For the connected TQFT's the vector spaces $\V(\Sigma)$ are thus identified
with invariances of special objects 
 in $\cal C$. When we pass to the disconnected case we have
to divide by the null spaces of the pairing with the respective coinvariance,
as it is stated in the summary in Lemma~\ref{lm-sec-43}. In the derivation
we actually first construct in (\ref{eq-quot-mor}) the linear morphism 
spaces, in which the  $\tau$-move of closed surfaces is realized, and then
identify these in Lemma~\ref{lm-H-on-vs} 
as linear maps between the quotiented vector spaces.

In Section 4.4 we tie the last remaining property, 
regarding the projective factor, to the existence of a natural transformation
of the identity functor of $\cal C\,$, whose image for each object 
is a multiple of the unit, see (\ref{eq-moni-epi-lambda}). In 
Lemma~\ref{lm-x-Css} we show that non-triviality of the value of 
such a transformation on the unit object itself (which will be
the same as {\tt x}) is a necessary and sufficient condition for 
the semisimplicity of $\cal C\,$. In the remainder of this section we 
establish the existence of such a transformation by identification
with  the {\em cointegral} of the coend $F=\int X^{\vee}\btimes X\,$,
assuming that  this is contained  in $\cal C\,$.

Finally, in Section 4.5 we combine the results of the previous sections and
of [KL] in Theorem~\ref{thm-main}, in which we establish the existence of
a large class of truly half-projective TQFT's. We also use the last section 
to speculate on the possibility of constructing generalized TQFT's, where we 
consider besides the tensor products also derived functors, like {\bf Tor},
 whose 
contributions may allow us to salvage some of the TQFT data that is lost
in the  division by the null spaces. 
As a further possible application of the formalism of half-projective TQFT's
we give a brief discussion of  classical limits, for which
${\tt x}\longrightarrow\infty\,$. We check for circle products 
the quality of the  estimate $\|\V(M)\|\,\geq const.\,{\tt x}^{\varrho(M)}\,$,
which is suggested by the Lemma~\ref{lm-ord-rho} and the normalizations used
in the canonical construction of the invariants. The estimates turn out to
holds in all of the considered cases, and they are
roughly half of the true order of divergence.
\medskip

The proofs for the basic, technical lemmas on coordinate graphs are
delivered in Appendix A.1. In Appendix A.2 we compute formulas for the
coboundaries $\mu_1$ and $\mu_{\partial}:=\mu_1-\mu_0\,$, which 
are generalizations in homology of $\mu_0\,$. We find 
$\mu_{\partial}\in{\bf Z}^{0,+}\,$. The corresponding anomaly in 
homotopy $\mu_{\pi}\,$ (see Section 3.2) counts the number of 
additional, non-separating surfaces in a product of cobordisms
that do not stem from the composites. Thus we pick up an additional 
factor, ${\tt x}^{\mu_{\pi}}\,$,  besides the one
 from the usual anomaly from (\ref{eq-I4}).
The tangle presentations of cobordisms from [Ke2], which we refer to
in Sections 4.1 and 4.4, are summarized in Appendix A.3.
\bigskip

\smallskip

\noindent
{\bf Acknowledgements:} I am indebted to P. Deligne and, especially,
V. Lyubashenko for many useful discussions, and, in particular, for
contributing Lemma~\ref{lm-int=1+1}. I also thank L. Rozansky for making
me aware of the example of supergroups.  Further, I wish to thank 
A. Beilinson, L. Kauffman, D. Kazhdan, and T. Kimura for interest and
encouraging discussions, and J. Stasheff for reading the
manuskript and making suggestions for improvements.
The author is supported by NSF-grant  DMS-9304580.

\bigskip

\bigskip

\cht 2 {Cobordism Categories, and Half-Projective TQFT's}

In this chapter we shall define and discuss generalizations of the 
TQFT-axioms of Atiyah. To this end we first introduce in Section 2.1 
the cobordism category \Cb *\ , whose objects are compact Riemann surfaces
with boundaries, and whose morphisms are homeomorphism classes of cobording
3-manifolds. Decompositions into connected components are expressed in
Section 2.2 in terms of the symmetric tensor structure of \Cb *. Here we also
introduce a coboundary, $\mu_0:=-\delta\be_0\,$, on \Cb *, where 
the coboundary operator  $\delta$ is again to be uncerstood in the 
generalized, group  cohomological sense, see (\ref{def-cocyc}).
The value of $\mu_0(M,N)\,$ is always a non-negative
integer, and it is $K-1$, if $M$ and $N$ are two connected cobordisms 
that are glued together along $K$ connected surfaces. 

This allows us in Section 2.3 to define the notion of half-projective  
TQFT's, $\V:\Cb 0\longrightarrow {\tt R}-mod\,$, on the category of
closed surfaces, $\Cb 0\,$, by inserting an anomaly
of the form ${\tt x}^{\mu_0}\,$ into the composition law for \V, where
{\tt x} does not have to be invertible. We show that
- except for specified, exceptional situations -  a half-projective
TQFT-functor, \V, is the sum of functors, $\V^j\,$, where each $\V^j\,$ maps
into the free modules of a summand, ${\tt R}_j\subset {\tt R}\,$,
and the component of {\tt x} in ${\tt R}_j\,$ is either zero or invertible.

Finally, in Section 2.4 we also discuss the various formalisms of extended TQFT's,
and how the notion of half-projective TQFT's can be extended to \Cb *.

\newpage

\prg{2.1}{ Categories of Cobordisms, the Structure of \Cb * }

The cobordism categories, which we wish to consider here, are slightly
more general than the ones defined, e.g., in [Ke1] or [KL]. 
The objects are as usual given by a set of inequivalent, compact,
oriented Riemann surfaces, $\Sigma\,$. Here we are only interested in
$\Sigma$ as a topological manifold. Moreover, we assume that $\Sigma$ is 
equipped with an ordering of its components.

We also fix parametrizations of
 the  boundary components,  $\partial\Sigma \cong \sqcup^n S^1\,$.
A morphism, $M\,:\,\Sigma_s\to\Sigma_t\,$, is  now defined 
 between any two such surfaces with $n_s$ and $n_t$ holes, respectively,
if the total number, $n_s+n_t$, of boundary components
 of $\partial M$ is even. 

We may organize the the set of holes into pairs, such that only holes of
 different surfaces are matched.
For any such choice  we glue in  cylinders 
connecting the boundary components of a pair  so that
we obtain a closed surface as follows:
\beq\lll{eq-sig-cl}
\Sigma^{cl}\,=\;-\Sigma_s\,\sqcup_{\sqcup^{n_s}S^1}\,
\Bigl\{{\bigsqcup}^{\frac{n_s+n_t} 2} S^1\times[0,1]\Bigr\}
\,\sqcup_{\sqcup^{n_t}S^1}\,\Sigma_t\qquad.
\eeq

A cobordism consists  now of an oriented, compact
3-manifold, $M\,$, and an orientation preserving 
homeomorphism, 
$\psi^0\,:\;\Sigma^{cl}\isto\partial M\,$.
Let us denote also the resulting inclusion  of  the source and target
surfaces into the cobordism:
\beq\lll{eq-psi-def}
\psi\,:\;-\Sigma_s\,\sqcup\,\Sigma_t\,\into {20} \,M\qquad.
\eeq
These maps will be sometimes called  {\em charts} or parametrizations.
\medskip
 
We shall denote the set of homeomorphie classes of such cobordisms with
fixed numbers, $n_s$ and $n_t$, of source and target surface components by
$\Cb {n_s,n_t}\,$.
\medskip

Below we depict a typical situation of how the boundary 
$\Sigma^{cl}\cong\partial M\,$ is built up, where $M:\,\Sigma_s\to\Sigma_t\,$
is a cobordism in $\Cb {3,5}\,$.
\smallskip

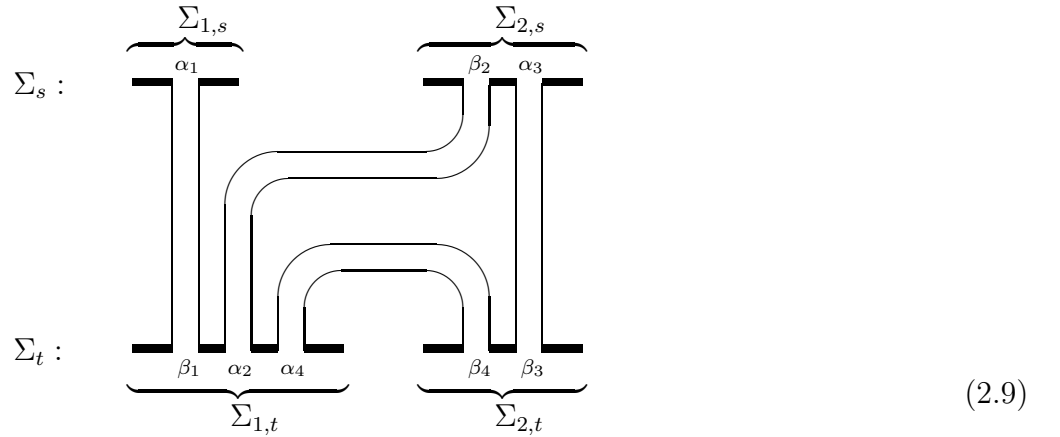
\beq\lll{diag-cob-ex}
\begin{picture}(170,140)

\put(-40,16){$\Sigma_t:$}

\put(-40,117){$\Sigma_s:$}

\put(5,19){\rule{15pt}{1mm}}
\put(30,19){\rule{10pt}{1mm}}
\put(50,19){\rule{10pt}{1mm}}
\put(70,19){\rule{15pt}{1mm}}

\put(115,19){\rule{15pt}{1mm}}
\put(140,19){\rule{10pt}{1mm}}
\put(160,19){\rule{15pt}{1mm}}

\put(20,19){\line(0,1){102}}
\put(30,19){\line(0,1){102}}

\put(150,19){\line(0,1){102}}
\put(160,19){\line(0,1){102}}

\put(60,20){\line(0,1){10}}
\put(70,20){\line(0,1){15}}

\put(130,20){\line(0,1){15}}
\put(140,20){\line(0,1){10}}

\put(95,35){\oval(50,30)[lt]}
\put(105,35){\oval(50,30)[rt]}
\put(95,50){\line(1,0){10}}
\put(100,30){\oval(80,60)[t]}

\put(50,20){\line(0,1){50}}
\put(75,70){\oval(50,30)[lt]}

\put(40,20){\line(0,1){45}}
\put(80,65){\oval(80,60)[lt]}

\put(75,85){\line(1,0){25}}
\put(80,95){\line(1,0){25}}

\put(105,110){\oval(50,30)[rb]}
\put(100,115){\oval(80,60)[rb]}

\put(130,120){\line(0,-1){10}}
\put(140,120){\line(0,-1){5}}

\put(3,8){$\underbrace{\hspace*{84pt}}$}
\put(42,-9){$\Sigma_{1,t}$}

\put(113,8){$\underbrace{\hspace*{64pt}}$}
\put(142,-9){$\Sigma_{2,t}$}

\put(5,120){\rule{15pt}{1mm}}
\put(30,120){\rule{15pt}{1mm}}

\put(132,11){$\scriptstyle\beta_4$}
\put(152,11){$\scriptstyle\beta_3$}

\put(22,11){$\scriptstyle\beta_1$}
\put(41,11){$\scriptstyle\alpha_2$}
\put(61,11){$\scriptstyle\alpha_4$}

\put(115,120){\rule{15pt}{1mm}}
\put(140,120){\rule{10pt}{1mm}}
\put(160,120){\rule{15pt}{1mm}}

\put(132,126){$\scriptstyle\beta_2$}
\put(151,126){$\scriptstyle\alpha_3$}

\put(21,126){$\scriptstyle\alpha_1$}

\put(3,132){$\overbrace{\hspace*{44pt}}$}
\put(22,142){$\Sigma_{1,s}$}

\put(113,132){$\overbrace{\hspace*{64pt}}$}
\put(142,142){$\Sigma_{2,s}$}

\end{picture}
\eeq
\medskip

In this example the source and target surface both have two components, i.e.,
$$
\Sigma_s^{(3)}\,=\;\Sigma_{1,s}^{(1)}\sqcup\Sigma_{2,s}^{(2)}
\qquad\;{\rm and}\;\qquad
\Sigma_t^{(5)}\,=\;\Sigma_{1,t}^{(3)}\sqcup\Sigma_{2,t}^{(2)}\quad,
$$
where the superscript at a surface  displays
the number of its holes. A surface component is
indicated in diagram (\ref{diag-cob-ex}) by a fat, 
horizontal line with an interruption for each hole.
We also introduced a  labelling of holes  so that
the hole with label $\alpha_j$ is connected to the
hole with label $\beta_j$ by the $j$-th
 cylindrical piece from (\ref{eq-sig-cl}), with $j=1,\ldots,4\,$.
Each of these pieces  are depicted in (\ref{diag-cob-ex}) by a
pair of thinner lines. For simplicity we have omitted in our 
example the possibility of crossings, although the cylinders may
be arbitrarily braided and knotted inside the cobordism.
\medskip

Note that in this definition we are allowed to have a cylinder 
connect two holes from two different target (or source) surfaces, as for
example the fourth cylinder between $\Sigma_{1,t}$ and $\Sigma_{2,t}$
in  (\ref{diag-cob-ex}). 
 For these we shall also specify a direction so that there
is a distinguished start- (or end-) hole. In the above example
we thus have to decide whether $\alpha_4$ or $\beta_4$ is the start
hole.

The union of all $\Cb {n,m}\,$ shall constitute $\Cb *\,$ as a set.
The composition in $\Cb *\,$ is  given by gluing the two 3-manifolds 
along the intermediate surface, requiring that end-holes
are glued to start-holes. 

In the process we may  encounter a situation, where
several cylindrical parts  combine  to form a  closed (connected) 
surface in $\partial M$, which can only be a  torus, $T^2$.
By construction 
this torus has a distinguished long meridian, with a given direction,
and a distinguished short meridian, which carries also a direction due to
the induced orientation. If we have similar data fixed on the respective
surface, $T^2$, that we chose as an object of $\Cb *$, then there is 
(up to isotopy) a unique homeomorphism between the two tori. Hence
we can add $T^2$ to either the start- of target-surface of the 
composite cobordism in a well defined manner. See also Section 2.4, where 
these tori are interpreted as ``horizontal objects'' of a 2-morphism.

On the morphism set we can also define a filling map 
\beq\lll{eq-fill}
\phi_0\,:\;\Cb{n_s,n_t}\;\to\;\Cb 0\qquad,
\eeq
which behaves nicely under compositions. It is given by gluing
a full tube $D^2\times [0,1]\,$ into the cylindrical parts, 
such that the holes of $-\Sigma_s\sqcup\Sigma_t\,$
are closed with discs $D^2\times\{0,1\}\,$. The cobordism $\phi_0(M)$ 
 is thus between the same surfaces without punctures.

In the definitions of [Ke2] and [KL] 
we only considered the case $n_s=n_t\,$, and the cylinders had to 
connect a source hole to a target hole. There $\phi$ is a 
true functor. Also, we considered  a 
central extension of $Cob_3(*)\,$
by the cobordism-group $\Omega_4$. This was naturally constructed by
 considering first also four-folds bounding the cobordism, and then
retaining their signature as an additional structure to the 3-cobordism,
see [Ke2] for details.  
We shall tacitly assume this extension  here, too.
\medskip

\prg{2.2}{Elementary Compositions, and the Cocycle $\mu_0$ }

In this section we shall introduce the coboundary $\mu_0$ on the 
cobordism category \Cb *, which enters the definition of a half-projective
TQFT. We start with the basics of the symmetric tensor structure of \Cb *.
We make the following  straightforward observations:
\medskip

A category of cobordisms between closed surfaces admits a natural
 tensor product given by the disjoint union. 
This tensor product is extends to the morphisms by disjoint unions in 
a functorial  way, and it is  obviously strictly associative.

Recall that for $\Cb 0$
we also assumed the connected components to be ordered. Thus if $\Sigma_j$ and
$\Sigma'_j$ are connected surfaces the tensor product of their ordered union
$\Bigl(\Sigma_1\sqcup\ldots\sqcup\Sigma_a\Bigl)\otimes
\Bigl(\Sigma_1'\sqcup\ldots\sqcup\Sigma_b'\Bigl)$ shall be the ordered union
$\Sigma_1\sqcup\ldots\sqcup\Sigma_a\sqcup\Sigma_1'\sqcup\ldots\sqcup\Sigma_b'\,$.

In particular, this means that $\Sigma\otimes\Sigma'$
and $\Sigma'\otimes\Sigma\,$ are different objects. However, we can find
a natural isomorphism, $\gamma:\,\Sigma\otimes\Sigma'\isto\Sigma'\otimes\Sigma\,$,
  in between them. In $\Cb 0$ this means that $\gamma$ is a cobordism, which has
a two-sided inverse. As a three-manifold with boundary it is given by
$\mbox{$\Bigl(\Sigma\times[0,1]\Bigr)$}\sqcup
\mbox{$\Bigl(\Sigma'\times[0,1]\Bigr)$}\,$ and
the boundary identifications 
$\Sigma\times\!\{0\}\sqcup\Sigma'\times\!\{0\}\isto\Sigma\sqcup\Sigma'$
and $\Sigma\times\!\{1\}\sqcup\Sigma'\times\!\{1\}\isto\Sigma'\sqcup\Sigma\,$.
are the obvious canonical maps. Also, it is clear that $\gamma\circ\gamma=\id\,$.

The system of morphisms $\gamma=\gamma(\Sigma,\Sigma')\,$ 
has in fact all the properties of a commutativity contraind of 
a symmetric tensor category. Specifically, it is natural with respect 
to bothe arguments, and the triangle equality is readily verified.

Let us summarize in the following lemma
the basic terminology that is natural from a 
categorial point of view and which provides a convenient language 
to organize reorderings of surface components  and decompositions
of cobordisms.

In the statements for  \Cb *\ we actually 
have to apply the adequate generalizations that take care
of the non-trivial vertical 1-arrows, like the extra tori discussed
in the previous section or renumberings of holes, see Section 2.4.

\blm\lll{lm-tens-cob}\  
\ben
\item \Cb *\  and \Cb 0\  are  strict, symmetric tensor categories. 
\item For any set, $\{\Sigma_j:\,j=1,\ldots,K\}\,$, of surfaces and any
permutation, $\pi\in S_K\,$, we have morphisms 
$\pi^*\,:\;\Sigma_1\otimes\ldots\otimes\Sigma_K
\,\to\,\Sigma_{\pi^{-1}(1)}\otimes\ldots\otimes 
\Sigma_{\pi^{-1}(K)}\,$, with $\pi\in S_K\,$, which are natural in every
argument, and which represent $S_K$.
\item For every morphism there is a unique, maximal number, 
$b=\beta_0\bigl(M\bigr)$, 
such that there is a  decomposition in the form 
$M\,=\,\pi_1^*\circ\Bigl(M_1\otimes\ldots\otimes M_b\Bigr)
\circ\pi_2^*\,$.
The cobordisms $M_j$ are then all connected.
\een
 \elm

The second remark in the lemma follows from the observation that cobordisms
of the form $\sigma^*=\id\otimes\gamma\otimes\id$, where $\gamma$ is the 
commutativity constraint on the union of two consecutive, connected components,
automatically fulfill the relations of the usual generators of the symmetric group.
Hence for a permutation $\pi\in S_K\,$ we can define the  
cobordism $\pi^*\,$ as the respective composition of the $\sigma^*$'s.

The last remark is simply expressing in categorial terms
the fact that, up to 
reordering of the boundary components, every cobordism
can be given as the union of its connected components.

In later chapters we will
return to the notation $\sqcup$ instead of $\otimes\,$.
\medskip

Adding a sufficient number of cylinders to the right and left of
each $M_j$ in the formula in Part 3 of Lemma~\ref{lm-tens-cob},
we obtain commuting morphisms 
$\breve M_j\,=\,\id_X\otimes M_j\otimes \id_Y\,$, such that the tensor product
$M$ can be rewritten as the composite:
 \beq\lll{eq-cob-dec}
M\;=\;\pi_1^*\circ\Bigl(\breve M_1\circ\ldots\circ\breve M_b\Bigr)\circ\pi_2^*\qquad.
\eeq

With this it follows that the compositions of two morphism can be obtained
as an iteration of two types of elementary composites. The first are of
the form $M\circ\pi^*$, where $M$ is connected. The second is given by 
products as follows:

\beq\lll{eq-elprod-cob}
M\;=\;\bigl(\id\otimes M_2\bigr)\circ\bigl(M_1\otimes\id\bigr)\qquad,
\eeq
where $M_1$ and $M_2$ are both connected, and are glued over $K\geq 1$
boundary components.
After having made a gluing along one connected component, the 
identifications of the remaining $(K-1)$ boundary components are
among boundary components of the same  connected manifold. 
We will  relate them to
$(K-1)$ uncancellable one-handle attachments in either a direct 
decomposition in three-dimensions, or to a four-manifold bounding $M$.

We wish to assign the excess number of connected boundary components over 
which we glue two cobordisms as a penalty in the form of a cocycle of \Cb *. 
To this end. let us introduce the {\em interior} Betti-numbers 
\beq\lll{eq-def-effbetti}
\be_j\bigl(M\bigr)\;=\;\beta_j\bigl(\phi_0(M)\bigr)\,-\,\frac 1 2 
\beta_j\bigl(\partial\phi_0(M)\bigr)\qquad,
\eeq
where $\beta_j\,=\,dim\bigl( H_j(X)\bigr)\,$ are the usual Betti-numbers.

In this section  we are interested only in the case $j=0$,
where we can omit the filling functor $\phi_0$ from (\ref{eq-fill})
in the formula. 
Some computations for $j=1$ are given in Appendix A.2.
It is easily seen that any number $\be_0(M)\in\frac 1 2 {\bf Z}\,$ is realized.
They define a coboundary on \Cb *\ with coefficients a-priori in the half integers
by
\beq\lll{def-cocyc}
\mu_0\bigl(M_2,M_1\bigr)\;:=\;-\delta\!\be_0\bigl(M_1,M_2\bigr)\;
=\;\be_0\bigl(M_2\circ M_1\bigr)\,-
\,\be_0\bigl(M_1\bigr)\,-\be_0\bigl(M_2\bigr)\qquad.
\eeq
It is readily seen that (\ref{def-cocyc}) actually defines an integer cocycle,
which is, in fact, also an integer coboundary of, e.g., 
$\tilde\be_0(M)=\beta_0(M)-\beta_0(\Sigma_t)\,$. 
The property that makes $\mu_0$ still an interesting quantity is that it is 
non-negative on all pairs of cobordisms. 

Indeed, for an elementary composition over $K\geq 1$
connected boundary components as in (\ref{eq-elprod-cob}), we find that this 
integer is the desired excess number of boundary-component over which we glue:
\beq\lll{eq-mu0-elprod}
\mu_0\bigl(\id\!\otimes\! M_2, M_1\!\otimes\!\id\bigr)\;=\;K-1\qquad.
\eeq
Also, it is easy to see that $\mu_0(\pi^*,M)=0$ for a permutation, and
$\mu_0(\id_X\!\otimes\!M_2,M_1\!\otimes\!\id_Y)~=~0$, if $Y$ is the source of $M_2$
and $X$ the target of $M_1$, so that 
\beq\lll{eq-mu0-recurs}
\mu_0\bigl(M_1\otimes M_2,\,N\bigr)\;=\;\mu_0\bigl(\id_X\otimes M_2,\,
(M_1\!\otimes\!\id_Y)\circ N\bigr)\;+\;\mu_0(M_1\otimes\id_Y,\,N)\,
\eeq
This  formula allows us to compute $\mu_0$ recursively 
from a presentation as in Lemma~\ref{lm-tens-cob}, and also 
implies $\mu_0\geq0$ for general compositions.

For a more systematic computation of the cocycle, let us introduce the spaces
\beq\lll{eq-def-W}
W_1\,:=\,ker\bigl(H_0(\psi^t_1)\bigr)\;,\qquad\qquad
W_2\,:=\,ker\bigl(H_0(\psi^s_2)\bigr)\qquad.
\eeq

Here the maps $\psi^{s/t}_j$ are the restrictions of the charts in
(\ref{eq-psi-def}) to source and target surfaces,
\beq\lll{eq-psi-rest}
\psi^t\,:\,\Sigma_t\,\hookrightarrow\,M\;,\qquad\qquad
\psi^s\,:\,-\Sigma_s\,\hookrightarrow\,M\qquad,
\eeq
for a cobordism $M\,:\,\Sigma_s\to\Sigma_t\,$.
We now have the following:

\blm\lll{lm-mu0=dim}
For two cobordisms, $M_1\,:\,\Sigma_{s,1}\to\Sigma\,$ and  
$M_2\,:\,\Sigma\to\Sigma_{t,2}\,$, and spaces $W_j$  as above, 
the cocycle from (\ref{def-cocyc}) is given as follows:
\beq\lll{eq-cocy=dim-0}
\mu_0(M_2,M_1)\;=\;dim\Bigl(W_1\cap W_2\Bigr)\;.
\eeq
\elm

{\em Proof: } The exact sequence $0\to W_1\cap~W_2\to H_0(\Sigma)\to 
H_0(M_2)\oplus~H_0(M_1)\to H_0(M_2\circ~M_1)\to 0\,$, implies the
dimension-formula $\beta_0(M_2\circ~M_1)- \beta_0(M_2)- \beta_0(M_1)=
dim(W_1\cap W_2)- \beta_0(\Sigma)\,$. Counting also boundary components,
this yields the asserted formula for $\mu_0$.\ep

A basic property of this cocycle is that it vanishes on invertible
cobordisms, i.e., 
\beq\lll{eq-mu=0-diff}
\qquad\mu_0(M,G)\,=\,\mu_0(G,M)\,=\,0 \qquad{\rm if}\quad G\in\,
\pi_0\Bigl({\cal D}i\!f\!f(\Sigma)^+\Bigr)\qquad.
\eeq
(Here we identify $G$ with a cobordisms by picking a representing 
automorphism of the surface, $\psi_G:\Sigma\isto\Sigma\,$. The 
associated cobordisms is then as a 3-manifold given by $\Sigma\times[0,1]\,$,
and the boundary identification are $id:\Sigma\times\!\{0\}\isto\Sigma\,$
and $\psi_G:\Sigma\times\!\{1\}\isto\Sigma\,$. It is a basic fact that this 
actually establishes an isomorphism between the mapping class group of 
$\Sigma$ and the group of invertible cobordisms from $\Sigma$ to itself). 
\medskip

Another property of $\mu_0$ is found in Section 3.2 and Appendix A.2:
If  $\mu_0(M_2,M_1)>0$, then the composite $M_2\circ M_1\,$ will contain 
paths or 1-cocycles that give rise to additional, infinite generators
of the fundamental group or the first homology group, respectively,
 besides those of $M_2$ and $M_1\,$.

\ms

\prg{2.3}{Half-Projective TQFT's, and Generalizations}

As we remarked in the introduction it is not possible to construct
non-trivial TQFT's in the classical sense, which vanish on $S^1\times S^2\,$.
The purpose of this section is to define the modification that makes 
allows such a construction and discuss a few basic implications. 

\begin{defi}\lll{def-half-ext}
Suppose {\tt R} is a commutative ring with unit, and \Rm\ the symmetric 
tensor category of free {\tt R}-modules. 
Further, let  {\tt x} be an element in {\tt R}, and 
 $\mu$ a 2-cocycle  on \Cb 0, which takes values only in ${\bf Z}^{+,0}\,$.

We call $\V\,:\;\Cb 0\,\to\,\Rm\,$ a \ub{half-projective}
TQFT (with respect to {\tt x} and $\mu$), if it full fills all the
requirements for  a functor of symmetric tensor categories, except for the
preservation of compositions. Instead of this we shall assume the relaxed
condition:
\beq\lll{eq-half-ext}
\V\bigl(M_2\circ M_1\bigr)\,\;=\;\,{\tt x}^{\mu(M_2,M_1)}\,\V(M_2)\V(M_1)\qquad.
\eeq
\end{defi}

As usual a projective TQFT is one for which ${\tt x}\in\R$ is {\em invertible}
and the non-negativity of $\mu$ is dropped.
An example of the latter is the well known signature-extension of the
2+1-dimensional cobordisms, which we avoided here by passing to a central
extension of \Cb * by $\Omega_4$. In this case $\mu$ is the Wall-cocycle. The 
number ${\tt x}\in {\bf C}^{\times}\,$ 
for TQFT's that are associated to  Chern-Simons theory with Lie-algebra 
{\bf g} and level $\ell $ 
may be obtained from 
the representation theory of Kac-Moody algebras, see [KW], and is a phase
depending on $\ell$, the dimension, and the dual Coxeter number of {\bf g}. 
For quantum-group constructions starting from a double $D(\B)\,$ it is
the pairing of square roots of the modulus and comodulus of $\B$, see [Ke1].

In general, if $\mu$ is a coboundary over {\bf Z}, we may rescale a projective 
TQFT-functor, and obtain a functor in the ordinary sense. 
However, for a half-projective TQFT and a trivial cocycle we can only 
replace {\tt x} by ${\tt x}{\tt y}^{-1}\,$ for invertible ${\tt y}\in \R$. 

This is the situation, which we are interested in here, as the non-semisimple
invariants will be  associated to half-projective TQFT's with respect to
the connectivity cocycle $\mu_0$ defined in the previous section. 
It is completely separated from the signature-extension. E.g.,  it does not 
 lead to any extensions of the mapping class group. In particular, we have
 the property expressed in (\ref{eq-mu=0-diff}), and  that 
$\mu_0$ is invariant under the natural $\Omega_4$ action.

Moreover, in the usual constructions of invariants,  \R\ is assumed to be a field
so that the only non-trivial, half-projective  TQFT occurs when ${\tt x}=0$,
(with the usual convention $0^0=1$). Let us,  however, continue to assume more 
general \R\ in the following discussion, in order to give more insight into
the underlying structures and show directions of possible, further generalizations.
\ms

Let us recall a general implication about the dimension of the \R-modules associated
to a surface that has already been observed by Witten [Wi] in the context of 
ordinary TQFT's.

\blm\lll{eq-torus=dim}
Suppose \V\ is a half-projective TQFT w.r.t. ${\tt x}\in\R$ and $\mu_0$. 
If $V_{\Sigma}\,=\,\V(\Sigma)\,$ is \R-module associated to a surface $\Sigma$, then
$$
\V(S^1\times\Sigma)\quad=\quad {\tt x}\,dim(V_{\Sigma})\qquad.
$$
\elm

{\em Proof:} For every connected surface $\Sigma$ let us fix an orientation 
{\em reversing} involution
$$
\hat\chi_{\Sigma}\in\pi_0
\Bigl(Di\!f\!f(\Sigma,\partial\Sigma)\Bigr)^{\bf -}_2\;.
$$
 If we consider 
(disjoint unions of) the $\hat\chi_{\Sigma}\,$'s as boundary charts of
the cylinder $\Sigma\times[0,1]\,$, we obtain morphisms 
\beq\lll{eq-def-chi}
\chi^{\dagger}_{\Sigma}\,:\,
\emptyset\to\Sigma\otimes\Sigma\qquad\; {\rm and}\;\qquad 
\chi_{\Sigma}\,:\,\Sigma\otimes\Sigma\to\emptyset\qquad,
\eeq
that are inverses of each other, and hence define a rigidity structure on \Cb *.
With $\chi_{\Sigma}^2=1$ they are also symmetric, in the sense that 
$\chi_{\Sigma}\circ\gamma=\chi_{\Sigma}\,$. Since \V\ preserves also the 
symmetric tensor structure  $\emptyset$ has to be associated to \R\ and 
$\gamma$ is mapped to the transposition of tensor factors. Thus,
if we apply \V\ to the $\chi_{\Sigma}^{(\dagger)}\,$ we obtain maps  
\beq\lll{eq-theta}
\theta_{\Sigma}^{\dagger}\,:\;\R\,\to\,V_{\Sigma}\otimes_{\tt R}V_{\Sigma}        
\qquad\;{\rm and}\;\qquad
\theta_{\Sigma}\,:\;V_{\Sigma}\otimes_{\tt R}V_{\Sigma}\,\to\,\R\qquad
\eeq
that are  symmetric with respect to the ordinary transposition, 
and are inverses of each other. Since we assumed the $V_{\Sigma}$
to be free \R-modules, we find by straightforward algebra that 
$\theta_{\Sigma}\theta_{\Sigma}^{\dagger}$ is the dimension of $V_{\Sigma}\,$.
We also know that 
$S^1\times\Sigma\,=\,\chi_{\Sigma}\chi^{\dagger}_{\Sigma}$. The anomaly
of this product is $\mu_0=1$ so that we find the asserted formula.\ep

From the above proof we see that, in fact, we do not have to assume that the 
$V_{\Sigma}$  are free-modules.  The existence of the morphisms in (\ref{eq-theta})
implies an isomorphism $V\cong Hom_{\tt R}\bigl(V,{\tt R}\bigr)\,$, 
which applied to $\theta_{\Sigma}^{\dagger}\,$ gives an element 
$\sum_{\nu}e_{\nu}\otimes l_{\nu}\,\in
\, V\otimes_{\tt R}Hom_{\tt R}\bigl(V,{\tt R}\bigr)\,$ 
that inverts the canonical pairing.

Suppose $X$ and $Z$ are {\tt R}-modules and $f:\,V_{\Sigma}\to X\,$ and 
$p:\,Z\onto{10}\;\, X\,$ are {\tt R}-morphisms, where $p$ is onto. For 
$x_{\nu}:=f(e_{\nu})\,$ and $p(z_{\nu})=x_{\nu}\,$, we can define a map
$h:\,V_{\Sigma}\to Z\,$ by $h(v)=\sum z_{\nu}l_{\nu}(v)\,$ so that $f=p\circ h\,$.
Hence $V_{\Sigma}\,$ is projective and therefore are a direct summand of a 
free {\tt R}-module.

Note also, that if ${\tt R}={\tt R}^0\oplus {\tt R}^1\,$ then we have a direct 
sum decomposition $V_{\Sigma}=V_{\Sigma}^0\oplus V_{\Sigma}^1\,$ using the
idempotents that are given by the units in the ${\tt R}^j\,$. Moreover,
we have $V\otimes_{\R}W\,=\,V^1\otimes_{\R_1}\!W^1\,\oplus\,
V^0\otimes_{\R_0}\!W^0\,$, etc. In summary, we find the following:

\blm\lll{lm-ring-dec}
Suppose $\V$ is a TQFT into possibly non-free {\tt R}-modules, and 
${\tt R}\,=\,\bigoplus_j{\tt R}^j$, where ${\tt R}^j$ are indecomposable.

Then $\V\,=\,\bigoplus_j\V^j\,$, where each $\V^j$ is a functor into the
category of free ${\tt R}^j$-modules.
\elm

A modification  of the prerequisites that would be consistent with 
 different values of 
$\V(S^1\times \Sigma)\,$ and hence  {\tt x},  is to allow the 
symmetry-structure of \Rm\ to be
different from that induced by $\V$. Thus the $\theta_{\Sigma}\,$
are now symmetric only up to isomorphism, i.e., we have 
$\theta_{\Sigma}T\,=\,\theta_{\Sigma}(\id\otimes P_{\Sigma})\,$,
where $T$ is the ordinary transposition and 
$P_{\Sigma}\in Aut_R(V_{\Sigma})\,$. Instead of the dimension of 
$V_{\Sigma}\,$ we then obtain the trace over $P_{\Sigma}$, which may 
even be zero.

The induced symmetry structure yields in place of the transposition
the morphisms,  
$\hat\gamma_{\Sigma}\in End(V_{\Sigma}\otimes V_{\Sigma})\,$, 
which are the images of the $\gamma$ as in 
Lemma~\ref{lm-tens-cob} with $\Sigma=\Sigma_j\,$. It is not hard to see
that for $\R={\bf C}\,$ 
the structure is equivalent to the canonical one if and only if
$tr(\hat\gamma_{\Sigma})=dim(V_{\Sigma})\,$.
\ms

It is also apparent from Lemma~\ref{eq-torus=dim} that we should not write
the anomaly-term to the other side of the equation in Definition~\ref{def-half-ext}.
For ${\tt h}$``$={\tt x}^{-1}$'' this would imply that {\tt h} divides 
the dimensions of the vector spaces for every genus. Under the usual
assumption (see also V3 of Section 4.2)
that the vector space of the sphere is \R\ this would imply that {\tt h} is 
invertible.

The next lemma only uses the composition rule and applies also to
the more general settings alluded to above:

\blm\lll{lm-torus-ord}
Suppose \V\ is a half-projective TQFT w.r.t. ${\tt x}\in\R$ and $\mu_0$.
For a connected surface $\Sigma$ of genus $g$ we then have:
$$
\V(S^1\times \Sigma)\quad\in\quad {\tt x}^{max(g,1)}\,\R\qquad.
$$
\elm

{\em Proof: } Consider the two-dimensional 
four-holed sphere as a 1+1-dimensional cobordism 
$H_4:\,S^1\sqcup S^1\to S^1\sqcup S^1\,$, and, further, let
$\phi^{\dagger}:\, \emptyset\to S^1\sqcup S^1\,$ and 
$\phi:\, S^1\sqcup S^1\to\emptyset\,$ be given by two-holed 
spheres. For $g\geq 1$ we have 
$\Sigma=\, \phi\circ\bigl(H_4)^{g-1}\circ\phi^{\dagger}\,$, and thus
$$
S^1\times\Sigma\;=\;\, \bigl(S^1\times\phi\bigr)\circ
\Bigl(S^1\times H_4\Bigr)^{g-1}\circ\bigl(S^1\times\phi^{\dagger}
\bigr)\qquad.
$$
Here, every one of the $g$ compositions is over two tori in the boundaries
of two connected cobordisms. Their anomalies are thus always $\mu_0=1$
and the assertion follows from the definition of a half-projective TQFT.
\ep

If we set $d_g=dim(V_{\Sigma_g})\,$, where $\Sigma$ has genus $g\geq 1$,
the combination of Lemmas~\ref{eq-torus=dim} and  \ref{lm-torus-ord}
yields that $d_g{\tt x}\in {\tt x}^g\R\,$. Suppose for some $g\geq 2\,$,
we have already  ${\tt x}\in {\tt x}^g{\tt R}\,$. Then there is ${\tt y}\in\R$
with ${\tt x}(1-{\tt x}{\tt y})=0$. In particular, ${\tt e}={\tt x}{\tt y}$
is an idempotent in \R, which can be used to write \R\ as a direct sum
${\tt eR}\oplus(1-{\tt e})\R\,$. Now {\tt x} lies in the first summand,
and {\tt y} is an inverse of this sub-ring with identity {\tt e}.
In summary, we have the following strong restriction on the
element {\tt x} and the dimensions of the vector-spaces.

\begin{cor}\lll{cor-altern}
Suppose \V\ , {\tt x}, and $d_g$ are as above. 
Then at least one of the following two has to be true:
\ben
\item The dimensions $d_g$ are zero-divisors in \R/${\tt x}^g\R$ for every
$g\geq 2$, \ or
\item The ring is a direct sum $\R\,=\,\R_1\,\oplus\,\R_0\,$, where
the component of {\tt x} in $\R_0$ is zero, and the component in 
$\R_1$ is invertible (in $\R_1$).
\een
\end{cor}

As in Lemma~\ref{lm-ring-dec} the second possibility implies for the
TQFT-functor that 
$\V=\V_1\oplus\V_0\,$ , and  we have that $\V_1$ can be rescaled to
an ordinary TQFT. The only non-trivial half-projective TQFT we
can therefore get if the dimension condition fails to hold
(and if we stay strictly within the framework of 
Definition~\ref{def-half-ext}) is  one with ${\tt x}=0\,$.
\ms

\prg{2.4}{TQFT's for Cobordisms with Corners}
\nopagebreak

There are several ways of defining {\em extended} TQFT's, which  represent
categories of cobordisms with corners, like $\Cb *\,$. Most of them are
consistent, although not always precisely equivalent. In this section
we shall give a brief survey over  the structures that are of interest to us.
\medskip

To begin with the {\em Kazhdan-Reshetikhin ladder}
 is defined on a series of categories,
$\Cb n\,\subset\, \Cb {n,n} \,$, for which all of the cylinders in (\ref{eq-sig-cl})
start at a hole in the source surface, and end at a hole in the target surface.
The extended TQFT is then defined, for a given abelian category $\cal C$, as a
series of functors 
\beq\lll{eq-ext-TQFT}
\V_n\,:\;\Cb n\,\longrightarrow\,
\underbrace{{\cal C}\odot\ldots\odot{\cal C}}_{n\;\,{\rm times}}\qquad,
\eeq
where $\odot$ is Deligne's tensor product of categories, see [D]. In particular,
to a surface, $\Sigma\,$, with $\partial\Sigma =\sqcup^nS^1\,$, this associates an
object $X_{\Sigma}\,\in\,{\cal C}^{\odot n}\,$. We also
require compatibility with the topological tensor product, i.e., 
$\,X_{\Sigma_1\sqcup\Sigma_2}= X_{\Sigma_1}\odot X_{\Sigma_2}\,$. Moreover, 
cobordisms are mapped to morphisms in the respective category, 
${\cal C}^{\odot n}\,$.
\medskip

Quite often it is more convenient to consider fiber functors that depend on a 
{\em coloration}. By this we mean an assignment of objects, 
$\{X_j:j=1,\ldots, n\}\,$, to the special, cylindrical pieces in the boundary of
a cobordism, $M:\,\Sigma_s\to\Sigma_t\,$, in $\Cb n\,$.
Instead of $f:=\V_n(M):\,X_{\Sigma_s}\to X_{\Sigma_t}\,$, we then consider the
following linear spaces and maps:
\beq\lll{eq-???}
\V_n^{\scriptscriptstyle(X_1,\ldots,X_n)}(M):\,
Hom_{{\cal C}^{\odot n}}\Bigl(X_1\odot\ldots\odot X_n,\,X_{\Sigma_s}\Bigr)
\TO{$\scriptstyle  Hom\bigl(\bigodot_j\!X_j,\,f\bigr)\,$}{30}
Hom_{{\cal C}^{\odot n}}\Bigl(X_{\pi(1)}\odot\ldots\odot X_{\pi(n)},
\,X_{\Sigma_t}\Bigr),
\eeq
where $\pi\in S_n\,$ is the permutation of the holes, defined by the cylinders in
$\partial M$, with respect to standard orderings of the holes 
in $\Sigma_s$ and $\Sigma_t$.

Notice, that the maps in (\ref{eq-???}) full fill an obvious naturality condition,
w.r.t. any given object, $X_j\,$, which appears both in the source and the target
linear space. Conversely, suppose any  functor, ${\cal C}\to {\tt R}-mod\,$, of 
abelian categories (i.e., not necessarily tensor) is representable. Then any such
system of maps with the naturality property stems from a functor like the 
one in (\ref{eq-ext-TQFT}). Examples of categories with representable fiber functors
are all those, for which the coend $\F\,:=\,\int X^{\vee}\odot X\,\in\,
{\cal C}^{\odot 2}\,$ exists. See [M], [L1], and [Ke3] for  definitions.
\medskip

If we distinguish between in- and out-holes among the boundary components of a 
surface, we can view the objects  of $\Cb n$ as 1+1-dimensional cobordisms
themselves. Thus it is quite natural to define $Cob_{1+1+1}\,:=\,
\bigcup_{n_{in}, n_{out}}\Cb {n_{in}+n_{out}}\,$ as a 2-category, where the
objects are one-manifolds, the 1-morphisms are 1+1-cobordisms, and the
2-morphisms are 2+1-cobordisms between them. An extended TQFT is now 
a 2-functor of 2-categories:
\beq\lll{eq-2cat-TQFT}
\V\,:\;Cob_{1+1+1}\,\longrightarrow\,{\bf AbCat}\qquad.
\eeq
Here, {\bf AbCat} is the 2-tensor-category of abelian categories. I.e.,
the object associated to a one-fold, $S$, is as usual an abelian category,
${\cal C}^{(S)}={\cal C}^{\odot \beta_0(S)}\,$, but to a surface we associate a 
functor between the category of the in-holes and the category of the out-holes.
To a cobordism between two surfaces $\V$ then assigns a natural transformation 
between the respective functors. This picture may be extracted from the previous 
one,
if we construct from an object $X_{\Sigma}\in {\cal C}^{(in)}\odot
{\cal C}^{(out)}\,$ the functor
$$
{\cal F}_{\Sigma}\,:\;{\cal C}^{(in)}\,\longrightarrow\,{\cal C}^{(out)}\;:\;
X\,\mapsto\,Hom_{{\cal C}^{(in)}}\Bigl(X_{\Sigma}, X\Bigr)\;,
$$
and from a morphism $f:\,X_{\Sigma_1}\to X_{\Sigma_2}\,$ a transformation
$\,{\cal F}_{\Sigma_2}\TO{$\scriptstyle\bullet$}{20}{\cal F}_{\Sigma_1}\,$ 
in the obvious way. 

Furthermore, the 2-categorial description imposes more
constraints on $\V$ than the Kazhdan-Reshetikhin picture, 
as \V\  has to be compatible not only with
compositions of 3-dimensional cobordisms but also with those of the 2-cobordisms. 
\medskip

The observant reader might have noticed that we have suppressed here the 
permutation that appears in (\ref{eq-???}). Indeed, for a precise description
we need to consider a slightly more complicated structure for $Cob_{1+1+1}$
than that of a 2-category. Specifically, in the more general formalism 
the objects of two cobordant  1-morphisms are
not simply the same but are connected by an arrow from a special category.
Hence the 2-arrow-diagram of a 2-morphism is not simply given by a 2-gon, but
by a square as below:
\beq\lll{eq-2-arr}
\bar{cccc}
\vspace*{.2cm}&\,S^{(in)}_s&\TO{$\Sigma_s$}{45}&\,S^{(out)}_s\\
\vspace*{.2cm}&\alpha\Biggl\downarrow\;\;&M\;\Downarrow&\Biggl\downarrow\beta\\
&\,S^{(in)}_t&\TO{$\Sigma_s$}{45}&\,S^{(out)}_t
\ear
\eeq
The new vertical arrows are associated to the special cylindrical pieces of 
$\partial M$, and are defined by the permutation they induce on the numbering
of the holes, i.e., we have $\alpha\in S_{n_{in}}\,$ and $\beta\in S_{n_{out}}\,$.

Horizontal compositions are only allowed if the adjacent vertical morphisms
are identical. For vertical (ordinary) compositions of 2-morphisms 
the special, vertical arrows are also multiplied.

The TQFT functor \V\ shall now assign to a vertical permutation in $S_n$  
the obvious
functor on ${\cal C}^{\odot n}\,$ that implements the respective permutation of
tensor factors. A cobordism $M$ is then mapped by \V\ to  a natural transformation 
between
the two composites of functors that start at the category of the upper, left
corner of the square in (\ref{eq-2-arr}) and end in the lower, right corner.
\medskip

Recall, that in Section 2.1 we actually defined a more general class of cobordisms,
for which the cylindrical pieces are allowed to run from a hole in a 
component of, e.g., 
the source surface, $\Sigma_s$, to a hole in another component of $\Sigma_s\,$.

A natural way to incorporate this possibility in our description is to enlarge
the category, from which we may take the vertical arrows, from the 
symmetric groupoid to the category of singular tangles (i.e., strands for 
which an overcrossing can be changed to an undercrossing).
A closed component of such a tangle, which in this category can be isolated as
a circle, corresponds to an interior torus that can be added in a unique way
as a closed component to either $\Sigma_s$ or $\Sigma_t$ as explained in 
Section 2.1.  

In order to define \V\ on a singular tangle it suffices to give the action on
a maximum or minimum:
\beq\lll{eq-V-cup-cap}
\bar{crclrcl}
\vspace*{.2cm}&\V\Bigl(\bigcup\Bigr)\,&:&\;{\cal C}\odot{\cal C}\,\longrightarrow\,{\tt R}-mod
&\qquad\quad X\odot Y&\,\mapsto\,& Hom_{\cal C}(1,\,X\btimes Y)\,\\
&\V\Bigl(\bigcap\Bigr)\,&:&\;{\tt R}-mod\,\longrightarrow\,{\cal C}\odot{\cal C}
&{\tt R}&\,\mapsto\,&\F=\int X^{\vee}\odot X\;,
\ear
\eeq
where $\F$ is the coend as above. Here $\btimes$ is a (braided) tensor product in
$\cal C$. Note that (\ref{eq-V-cup-cap}) also implies 
$\V(\bigcirc):{\tt R}\mapsto Inv(F)\,$, where $F=\int X^{\vee}\btimes X\,$.
\medskip

As an alternative to this direct functorial description,  
we may consider also here the 
assignments of linear maps depending on  a given coloration, as in (\ref{eq-???}).
The difference is now that a cylindrical piece, which starts and end
in the source surface results in a dependence of the source vector space
on $X_j^{\vee}\odot X_j$ (instead of only $X_j$) and no dependence of the 
target linear space. The first then full fills a di-naturality condition
(instead of a naturality condition), which also plays an important role in
liftings to the coend $\F\,$.
\medskip

The notion of a half-projective, extended TQFT is most conveniently defined
for the version of functors into ${\tt R}-mod$ that depend on colorations.
The generalization from Definition~\ref{def-half-ext} is then literarily
the same. It is also not hard to generalize the construction and discussion
of half-projective TQFT for closed surfaces to the case of punctured surfaces,
since all that needs to be checked in this picture is the preservation of 
naturality. We shall thus content ourselves in this article with a construction
of half-projective functors $\V_0:\Cb 0\to{\tt R}-mod\,$.
Only in Section 4.3 will we return to the direct description of an extended TQFT
in terms
of functors as in (\ref{eq-ext-TQFT}).

\cht{3} 
{Non-Separating Surfaces, Interior Fundamental Groups,\\
 and Coordinate-Graphs}

In Lemma~\ref{lm-torus-ord} we used a decomposition of a manifold to show that
it is mapped by \V\ to a multiple of ${\tt x}^{\varrho}\,$ for some 
$\varrho\in{\bf Z}^{+,0}\,$. In this chapter we wish to generalize this 
result, and identify for any type of cobordism, 
$M\,$, orders in {\tt x} for $\V(M)\,$, that are maximal for this decomposition
argument.(It will turn out that the one in Lemma~\ref{lm-torus-ord} is indeed
maximal).  

From the discussion in Section 2.3  it seems that {\tt x}=0
is the only case we should be worried about, i.e.,  the only relevant 
question would be whether $\varrho=0$ or not, but not the precise order
for $\varrho >0\,$. Nevertheless, we shall 
stay within the more general framework, not only because of the 
possible modifications we outlined in Section 2.3, but also because
of anticipated applications to ``classical limits'', which we will
sketch in Section 4.5.

In Section 3.1 we relate the  orders, $\varrho$, 
to the  maximal number of non-separating surfaces in a cobordism. 
The subsequent sections are devoted to computing these numbers from
the fundamental groups of $M$ and $\partial M\,$. 
Specifically, we will find in Section 3.5 that $\varrho$ is the maximal rank of 
a free group, $F$, for which there is an exact sequence of the form,
$\pi_1(\partial M)\longrightarrow\pi_1(M)\longrightarrow F\longrightarrow 1\,$.
For a special case see Lemma 6.6 in [He]. As an application we find that a 
half-projective TQFT, \V, with ${\tt x}=0\,$, vanishes on cobordisms, $M$,
with nontrivial ``interior Betti-number'', i.e., $\V(M)=0$ if $\be_1(M)\neq 0\,$.

An important tool in this discussion
are the coordinate graphs, which reduce the relevant connectivity information
of cobordisms in \Cb * to that of morphisms in a graph-category \G.
A useful result, proven in Appendix A.1,
 is that decomposition along such graphs can also be realized
as decompositions of the corresponding cobordisms.

\prg{3.1} {r-Diagrams of Non-Separating Surfaces}

Let us begin with a definition the systems of non-separating surfaces we are 
interested in:

\begin{defi}
For a manifold $M$ with boundary an \ub{$r$-diagram} is an embedding of 
$r$  Riemann surfaces $\Sigma_j\hookrightarrow M$ with $j=1,\ldots,r$,
 such that
\ben
\item the surfaces are disjoint from each other, i.e., 
$\Sigma_i\cap\Sigma_j=\emptyset$ for $i\neq j$, 
\item they lie in the interior of $M$, i.e., 
$\partial M\cap\Sigma_j=\emptyset$,
\item the surfaces are closed so that the embeddings are proper,
\item every $\Sigma_j$ is  two-sided, and 
\item their union is non-separating, i.e., 
$M-\sqcup_{j=1}^r\Sigma_j$ is connected.
\een
\end{defi}

Let us reformulate the existence of an $r$-diagram for a cobordism, 
$M:\,\Sigma_s\to\Sigma_t$,
in \Cb 0 \ in a more categorial language. If we remove 
thin, two-sided collars, $\chi_j\,:=\,
\Sigma_j\times[-\varepsilon,\varepsilon]\,$,
from $M\,$, we obtain a manifold $M^*$, which  has
$2r$ additional boundary components, $\Sigma_j^{\pm}\,=\,
\Sigma\times\{\pm\varepsilon\}\,$. Hence,  as a morphism in \Cb 0,
it can be written as follows: 
$$
M^*\,:\;\Sigma_2\,\longrightarrow\,
\Sigma_t\otimes\Sigma_1^+\!\otimes\!\Sigma_1^-\!\otimes\ldots\otimes
\!\Sigma_r^+\!\otimes\!\Sigma_r^-\qquad.
$$
In terms of the rigidity
cobordisms from (\ref{eq-def-chi}) the original manifold is given by  
\beq\lll{eq-dec-sepsur}
M\;=\;\Bigl(\id_{\Sigma_t}\otimes\chi_{\Sigma_1}\!\otimes\ldots\otimes\!
\chi_{\Sigma_r}\Bigr)\circ M^*\quad.
\eeq
It is easy to see that the total anomaly of this product is $\mu_0=r$.

Let us also introduce the  quantity 
\beq\lll{eq-s}
\varrho(M)\,:=\,max\{r:\,M\, {\rm \ admits\ an\ }r{\rm-diagram}\}\qquad.
\eeq
Note at this point that  an $r$-diagram of some $M$ with $r<\varrho(M)$
can usually not be completed to a maximal diagram. An example is 
$S^1\times\Sigma_g\,$. Here we can find a $g$-diagram from the lower
curves of a  Heegaard diagram on $\Sigma_g\,$, which contains one
torus from every composition in the proof of Lemma~\ref{lm-torus-ord}.
We will, however, see
in the next section 
that after removing $\{1\}\times\Sigma_g\,$ we have for the
complement $\varrho([0,1]\times\Sigma_g)=0$.

From  (\ref{eq-dec-sepsur}) we find immediately for the following:

\blm\lll{lm-ord-rho}
Suppose \V\ is a half-projective TQFT w.r.t. ${\tt x}\in\R\,$ and $\mu_0$. 

Then we have for a cobordism $M:\Sigma_s\to\Sigma_t$ from \Cb 0, 
\beq\lll{eq-ord-rho}
\V(M)\;\in\;{\tt x}^{\varrho(M)}\,Hom_{\R}\bigl(V_{\Sigma_s},V_{\Sigma_t}\bigr)
\qquad.
\eeq
Here, $V_\Sigma=\V(\Sigma)\,$, and we mean $\bigcap_r{\tt x}^rHom(\ldots)$ if 
$\varrho(M)=\infty\,$.
\elm

The generalization of Lemma~\ref{eq-ord-rho} to extended
TQFT's as in Section 2.4 also holds true, if we  have $M\in \Cb *$ 
and replace $Hom_{\tt R}\,$
by $Hom_{{\cal C}^{\otimes N}}$ on the objects assigned by \V\ to the 
punctured source and target surfaces. For a generalization of the
previous arguments to surfaces with punctures let us make the following
observations: 

The definition of an $r$-diagram of a cobordism in \Cb *\ is simply 
an $r$-diagram of $\phi_0(M)\,$, where $\phi_0$ is as in (\ref{eq-fill}).
In order to obtain the  generalization of the composition in (\ref{eq-dec-sepsur})
we first have to make all the surfaces transversal to 
the 
external strands, $\tau=D^2\times I\,$, that connects holes of the surfaces 
$\Sigma_s\sqcup\Sigma_t\,$ to each other.
A given strand, $\tau$, is then divided by
its intersections with the $\Sigma_j$ into several components.
If one such piece, $\delta$, connects a surface $\Sigma$ to itself from
the same side, we can surger $\Sigma$ along a slightly thickened $\delta$,
and obtain a two-sided surface 
$\Sigma'=(\Sigma-D^2\times S^0)\cup (S^1\times I)\,$,
which is also disjoint from the other surfaces and, together with them,  is
non-separating. We can thus assume that strands never connect a surfaces to
itself from the same side.

By transversality we may also assume that an external strand, $\tau\,$, meets a 
collar, $\chi_j=
\Sigma_j\times[-\varepsilon,\varepsilon]\,$,
in  a vertical, cylindrical piece, $D^2_{\tau}\times[-\varepsilon,\varepsilon]\,$,
where $D^2_{\tau}\subset\Sigma_j$. Hence, if we remove all external strands, we
have presented $M$ as the composite of the unti-morphisms 
$\chi_j'=\Sigma_j'\times[-\varepsilon,\varepsilon]$, where $\Sigma_j'$ is obtained
from $\Sigma_j$ by removing the discs $D^2_{\tau}\,$, and an admissible,
connected cobordisms $M^*$ in \Cb *, as in (\ref{eq-dec-sepsur}).

\prg{3.2}{Interior Fundamental Groups, and an A-Priori Estimate on $\varrho(M)$}

The maximal number $\varrho(M)$ of non-separating surfaces in a cobordism
can be obtained from the fundamental groups of the  cobordism $M$ and its
boundary $\partial M$. It will be given by the maximal rank of a free group
onto which  the {\em interior} fundamental group factorizes.
In this section we give the definition of $\pint(M)$ and a first implication
for $\varrho(M)$.
\ms

If $M:\,\Sigma_s\to\Sigma_t$ is a cobordism in \Cb * and $\Sigma_{\nu}\,$,
with $\nu=1,\ldots, K\,$
are the connected components of $-\Sigma_s\sqcup\Sigma_t$ we 
define the {\em interior} fundamental group, $\pint(M)$ as the fundamental group of
$M$ with the
cones of the boundary components glued to it:
\beq\lll{eq-pint-cone}
\pint(M)\;=\;\pi_1\Bigl(M\cup \,{\rm C}\Sigma_1\cup\ldots \cup\,{\rm C}\Sigma_K\Bigl)
\eeq
The following shows that it is enough to consider cobordisms in \Cb 0.
\blm\lll{lm-pi0}
The inclusion $M\hookrightarrow\phi_0(M)\,$ induces an isomorphism
$$
\pint(M)\,\isto\,\pint\bigl(\phi_0(M)\bigr)\;.
$$
\elm
{\em Proof :} Suppose a cobordism $M'$ is obtained by filling a tube in $M$
that connects a holes the component $\Sigma_i$ to one in the component $\Sigma_j$.
Since for $M\in\Cb *$ we have $i\neq j$, the cones ${\rm C}\Sigma_i$ and
${\rm C}\Sigma_j$ are disjoint. Thus $M'$ with cone-attachments is obtained from
$M$ with boundary cones by gluing in a ball along a sphere, which does not
affect fundamental groups.\ep  

For the practical computation of  $\pint(M)$,
assume that  we have marked points $p_0\in M$ and $p_{\nu}\in \Sigma_{\nu}\,$.
Let us call a {\em spider}, $\lz\gamma\rz$,
of $M$ a collection of paths $\gamma_{\nu}$ 
inside $M$, with $\nu=1,\ldots, K$, that start at $p_0$ and end at $p_{\nu}$.
Thus  $X=\lz\gamma\rz\cup{\rm C}\Sigma_1\cup\ldots \cup{\rm C}\Sigma_K$ is 
a contractible space such that $M\cup X$ is the union of $M$ with its 
boundary cones as in 
(\ref{eq-pint-cone}), and $X\cap M\simeq\Sigma_1\vee\ldots\vee\Sigma_K\,$.
The interior group $\pint(M)$ is then given by 
Seifert-van Kampen (see, e.g., Theorem 7.40 in [Ro])
as the pushout of the respective fundamental groups, i.e., it is  universal 
among the solutions,  $\zeta$ and $G$,
of the following diagram:
\beq\lll{eq-pi-bdmap}
\bar{cccc}
\vspace*{.2cm}&{ ^{free}\prod^{K}_{\nu=1}}\pi_1(\Sigma_{\nu}, p_{\nu})&\TO{}{35}&0\\
\vspace*{.1cm}&I_*\Biggl\downarrow\;&&\Biggl\downarrow\\
&\pi_1(M,p_0)&\ONTO{$\zeta$}{35}&G
\ear
\;.
\eeq
The image of $I_*$ in $\pi_1(M)$
generally depends on the choice of the spider, but their generators
lie in the same conjugacy classes. Hence the smallest normal subgroup
${\cal N}[im(I_*)]\subset\pi_1(M)$ that contains the image of $I_*$ does not
depend on the choice of a spider. This yields  the formula:
\beq\lll{eq-def-pint}
\pint(M)\;:=\; {\pi_1(M)}/{\cal N}[im(I_*)]\qquad.
\eeq
We also introduce the notion of  a {\em free interior group}, $F$,  of $M$.
By this we mean a solution to the diagram (\ref{eq-pi-bdmap}), where $F=G$ is
a free group (non-abelian for rank $>1$ ), and $\zeta$ an epimorphism. 
Clearly, by universality it may also be defined by the existence of an epimorphism
\beq\lll{eq-free-int}
\tilde\zeta\,:\;\pint(M)\,\ONTO{}{22}\,F\qquad .
\eeq

Let us denote by $F(k)$ the free group in $k$ generators. The 
role of the anomaly of Section 2.2 for internal fundamental groups can be 
described as follows:

\blm\lll{lm-comp-freefac} Suppose $M$ and $N$ are cobordisms with anomaly
$\mu_0=\mu_0(M,N)\,$. Then there is an epimorphism $\xi\,$, such that the following
diagram commutes:
\beq\lll{eq-comp-freefac}
\bar{cccc}
\vspace*{.2cm}&\pint\bigl(M\circ N\bigr)&\ONTO{$\xi$}{25}&\pint(M)*\pint(N)*F(\mu_0)\\
\vspace*{.2cm}&\Biggl\uparrow&&\Biggl\uparrow\\
&\pi_1(M\circ N)&\longleftarrow&\pi_1(M)*\pi_1(N)*F(\mu_0)
\ear
\eeq
\quad .
\elm

{\em Proof:} It is enough to consider only connected cobordisms $M$ and $N$ that
are connected over $\mu_0+1$ surfaces. 
Moreover, the assertion is the same, if we glue
in the cones for the remaining boundary components of $M$ and $N$, i.e., we
may assume that $N:\emptyset\to \Sigma_0\sqcup\ldots\sqcup\Sigma_{\mu_0}\,$, and
$M$ is a cobordism in reverse direction.

Fix a point $x_0\in\Sigma_0\,$ and choose spiders $\lz \gamma^M\rz\,$ and 
$\lz \gamma^N\rz\,$
of $M$ and $N$ respectively that originate in $x_0$. Hence a leg, $\gamma^M_j\,$,
is a path in $M$ that connects $x_0\in\Sigma_0$  to a point  $x_j\in\Sigma_j\,$.

For $\tilde N\,=\,N\cup\lz \gamma^M\rz\,$ we have a natural isomorphism
$$
\pi_1(N,x_0)*F(\mu_0)\,\cong\,\pi_1\bigl(\tilde N,x_0\bigr)\qquad,
$$
in which the $j$-th free generator $a_j$ of $F(\mu_0)$ is mapped to 
$\bigl(\gamma_j^N\bigr)^{-1}\gamma_j^M\,$. 

Since $M\circ N\,=\,M\cup\tilde N\,$ and 
$M\cap\tilde N=\lz\gamma^M\rz\cup\Sigma_0\cup\ldots\cup\Sigma_{\mu_0}\,$, the group
$\pi_1(M\circ N)=\pint(M\circ N)\,$ can be computed as the push-out of the
following diagram:
\beq\lll{eq-comp-pusho}
\bar{cccc}
\vspace*{.2cm}&{ ^{free}\prod_1^{mu}\pi_1(\Sigma_j,x_j)\quad}&\TO{$\hat I^N_*$}{30}&
\pi_1(N)*F(\mu_0)\\
\vspace*{.2cm}&\Biggl\downarrow I^M_*&&\\
&\pi_1(M)&&,
\ear
\eeq
where $\hat I^N_*(g)\,=\,a_j^{-1}I^N_*(g)a_j\,$ if 
$g\in \pi_1\bigl(\Sigma_j,x_j\bigr)\,$, and the $I_*^{M/N}\,$ are defined as
in (\ref{eq-pi-bdmap}).

Clearly, $\pint(M)*\pint(N)*F(\mu_0)\,$ is a solution of (\ref{eq-comp-pusho}),
since $im\bigl[\hat I^N_*\bigr]\,$ is also in the kernel for the map
onto $\pint(N)*F(\mu_0)\,$. Hence a surjection $\xi$ exists.

The remainder of the diagram, expressing that $\xi$ is defined naturally,
can be completes in the obvious way. In the lower horizontal morphism
the free generator $a_j\in F(\mu_0)\,$ is mapped to the closed path  
$\bigl(\gamma_j^N\bigr)^{-1}\gamma_j^M\,$ in the composite $M\circ N$.
\ep

In analogy to $\varrho(M)$ from (\ref{eq-s}) let us define the maximal rank 
of a free interior group:
\beq\lll{eq-phi}
\varphi(M)\,:=\,max\bigl\{\mu:\,M\,{\rm\ has\ } F(\mu)
{\rm \ as\ free\ interior\ group}\bigr\}\;.
\eeq
The following are  easily found from Lemma~\ref{lm-comp-freefac} and 
(\ref{eq-dec-sepsur}):
\begin{cor}\lll{cor-2}\ 

\ben
\item $\varphi(M\circ N)\,\geq\, \varphi(M)\,+\,\varphi(N)\,+\,\mu_0(M,N)\;$.
\item $\varphi(M)\,\geq\,\varrho(M)\;$.
\een
\end{cor}

For $\varphi(M)$ and $\varphi(N)$ the first part implies the existence of
an other trivial cocycle with values in the non-negative integers, given by:
\beq\lll{eq-mu-pi}
\mu_{\pi}\,:=\,-\delta\varphi-\mu_0
\eeq
Its computation can be  quite intricate, and shall not be attempted here.
Instead we shall give the analogous computation for homology in Appendix A.2.

The second part of the corollary implies for one that 
$$
\varrho(M)\,\leq\,\varphi(M)\,<\,\infty\;,
$$
since the fundamental group of a compact manifold is finitely generated.
This also renders the convention made at the end of Lemma~\ref{lm-ord-rho}
superfluous. 
Another consequence is that with $\pi_1(\Sigma)\onto{20}\pi_1(M)$ being onto for 
invertible cobordisms, we have 
\beq\lll{eq-anom-inv}
\varphi(G)\,=\,\varrho(G)\,=\,0\qquad{\rm if}\quad G\in\pi_0
\Bigl({\cal D}i\!f\!f(\Sigma)^+\Bigr)\quad,
\eeq
which was used in the counter-example in Section 3.1. In the remaining  
sections of this chapter  we shall see that in fact 
$\varphi\equiv \varrho\,$.

\prg{3.3}{ The Graph-Category \G\  }

In this section we shall define a category of graphs. It will be  
used to encode the basic  connectivity properties of \Cb 0\ . 

To begin with let us fix a label set, $\S\,$,
 that is in one-to-one correspondence with
the Riemann surfaces, used as objects for \Cb 0.
The objects of the category \G\ are then given by strings of
(possibly repeated) labels, $[a_1,\ldots,a_K]\,$, $a_j\in\S\,$.

The morphisms, $\gamma\,:\,[a_1,\ldots,a_K]\,\to\,[b_1,\ldots,b_L]\,$,
are given by one-dimensional cell-complexes, $\gamma\,$, taken up to homotopy type,
 for which $\partial \gamma$ contains $K+L$ special points, that are labeled by
$a_1,\ldots,b_L\,$. The compositions is, as for cobordisms, given by 
gluings along the respective boundary components, i.e., end-points with
the same labels.

A representing cell-complex can be visualized by a graph, with $K+L$
distinguished vertices of edge-degree one. A generic example of a 
representing graph is depicted below. In this form 
the composition of two graphs is defined by placing them on top of each other. 
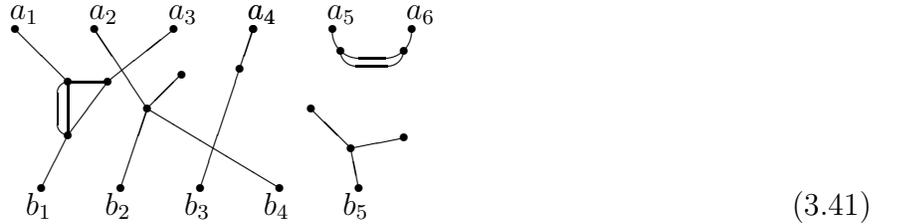
\beq\lll{fig-gr-ex}
\begin{picture}(150,81)

\put(20,10){\line(1,2){10}}
\put(30,30){\line(0,1){20}}
\put(30,40){\oval(8,20)[l]}
\put(30,50){\line(-1,1){20}}
\put(30,50){\line(1,0){15}}
\put(30,30){\line(3,4){15}}
\put(45,50){\line(5,4){25}}

\put(20,10){\circle*{3}}
\put(14,0){$b_1$}
\put(30,30){\circle*{3}}
\put(30,50){\circle*{3}}
\put(10,70){\circle*{3}}
\put(8,74){$a_1$}
\put(45,50){\circle*{3}}
\put(70,70){\circle*{3}}
\put(68,74){$a_3$}

\put(50,10){\line(1,3){10}}
\put(60,40){\line(-2,3){20}}
\put(110,10){\line(-5,3){50}}

\put(50,10){\circle*{3}}
\put(44,0){$b_2$}
\put(110,10){\circle*{3}}
\put(104,0){$b_4$}
\put(60,40){\circle*{3}}
\put(60,40){\line(1,1){13}}
\put(73,53){\circle*{3}}
\put(40,70){\circle*{3}}
\put(38,74){$a_2$}

\put(80,10){\line(1,3){20}}
\put(95,55){\circle*{3}}
\put(100,70){\circle*{3}}
\put(98,74){$a_4$}
\put(80,10){\circle*{3}}
\put(74,0){$b_3$}

\put(100,70){\circle*{3}}
\put(98,74){$a_4$}

\put(145,70){\oval(30,22)[b]}

\put(145,62){\oval(24,12)[b]}
\put(157,62){\circle*{3}}
\put(133,62){\circle*{3}}

\put(130,70){\circle*{3}}
\put(128,74){$a_5$}
\put(160,70){\circle*{3}}
\put(158,74){$a_6$}

\put(137,25){\line(-1,1){15}}
\put(122,40){\circle*{3}}
\put(157,29){\circle*{3}}
\put(137,25){\line(5,1){20}}
\put(140,10){\line(-1,5){3}}
\put(137,25){\circle*{3}}
\put(140,10){\circle*{3}}
\put(134,0){$b_5$}

\end{picture}
\eeq

\medskip 

In analogy to Lemma~\ref{lm-tens-cob} we also have a natural 
symmetric tensor structure on \G:
\blm\lll{lm-gr-ts}\  
\ben
\item \G\ is a  strict, symmetric tensor category. 
\item 
For any $\pi\in S_K\,$, there is a morphism, $\pi^*\,:\,[a_1,\ldots,a_K]\to[a_{\pi^{-1}(1)},\ldots,
a_{\pi^{-1}(K)}]$. They are natural in \G\ and give
 rise to a representation of $S_K$. Any invertible morphisms of \G\ 
is of this form.

\item For every morphism there is a unique maximal number, $b$, 
such that there is a  decomposition in the form 
$\gamma\,=\,\pi_1^*\circ\Bigl(\gamma_1\otimes\ldots\otimes\gamma_b\Bigr)
\circ\pi_2^*\,$.
A representing graph of a component, $\gamma_j$,  is connected.
\een
 \elm

The  tensor product of \G\ is given by the juxtapositions of both labels
and graphs. The permutations are given by joining a labels $a_j$ and 
$a_{\pi^{-1}(j)}$ on top and bottom by straight lines. 
Using  the triangle 
identity, this also
defines the commutativity constraint $\gamma$  on strings of arbitrary length, 
as the obvious crossing of two sets of parallel strands.

A graph, $\gamma$, that is invertible cannot connect two different 
source (or target) labels to each other, and cannot contain internal loops.
It follows for \G, that up to homotopy the permutations are the only 
possibilities. A similar statement hold true also for \Cb *, 
if we include the action
of the mapping class groups on the surfaces.  

A decomposition as in the last
part of Lemma~\ref{lm-gr-ts} is given for the $k=5$ component graph 
in (\ref{fig-gr-ex}), by the permutations 
 $\pi_1=(3,4)\,$ and $\pi_2=(2,3)\,$, and the connected 
morphisms $\gamma_1\,:\,[a_1,a_3]\to[b_1]$, $\gamma_2\,:\,[a_2]\to[b_2,b_4]$,
$\gamma_3\,:\,[a_4]\to[b_3]$, $\gamma_4\,:\,[a_5,a_6]\to[\emptyset]\,$, and 
$\gamma_5\,:\,[\emptyset]\to[b_5]\,$. 
\ms

As opposed to \Cb * \ it is easy to  give a list of the homotopy-inequivalent,
connected graphs, thus giving a complete description of the category.

The class of a  connected morphism, $\gamma$, is clearly  determined 
by the number of source labels, $K$, the number of target labels, $L$, and
the first Betti number $\beta_1(\gamma)\,=\,dim\bigl(H_1(\gamma)\bigr)$.  Canonical
representatives can be found by shrinking all of the internal edges, until we have
at most one internal vertex. The results are given in the next diagram.

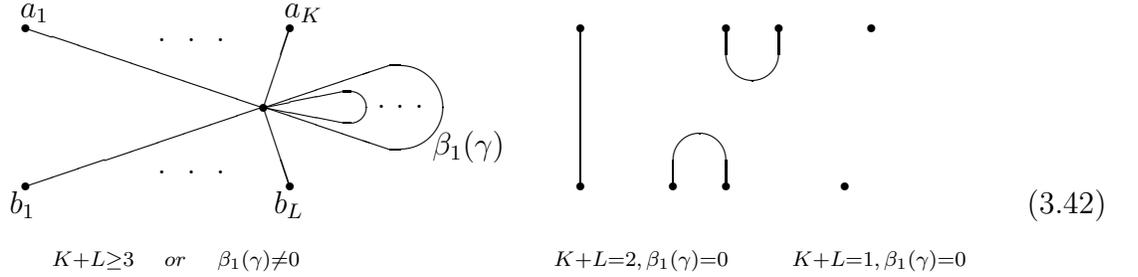
\beq\lll{fig-gr-reps}
\begin{picture}(320,80)

\put(100,40){\circle*{3}}

\put(10,70){\circle*{3}}
\put(8,74){$a_1$}
\put(100,40){\line(-3,1){90}}

\put(60,65){$.\ \ .\ \ .$}
\put(60,15){$.\ \ .\ \ .$}

\put(10,10){\circle*{3}}
\put(4,0){$b_1$}
\put(100,40){\line(-3,-1){90}}

\put(110,70){\circle*{3}}
\put(108,74){$a_K$}
\put(100,40){\line(1,3){10}}

\put(110,10){\circle*{3}}
\put(104,0){$b_L$}
\put(100,40){\line(1,-3){10}}

\put(100,40){\line(3,1){48}}
\put(100,40){\line(3,-1){48}}

\put(100,40){\line(5,1){30}}
\put(100,40){\line(5,-1){30}}

\put(130,40){\oval(18,12) [r]}
\put(148,40){\oval(40,32) [r]}

\put(143,40){$.\ .\ .$}

\put(164,22){$\beta_1(\gamma)$}

\put(20,-20){$\scriptstyle K+L\geq3 \quad or \quad  \beta_1(\gamma)\neq 0$}

\put(220,10){\line(0,1){60}}
\put(220,10){\circle*{3}}
\put(220,70){\circle*{3}}

\put(265,10){\oval(20,40)[t]}
\put(255,10){\circle*{3}}
\put(275,10){\circle*{3}}

\put(285,70){\oval(20,40)[b]}
\put(275,70){\circle*{3}}
\put(295,70){\circle*{3}}

\put(210,-20){$\scriptstyle K+L=2,\, \beta_1(\gamma)= 0$}

\put(330,70){\circle*{3}}
\put(320,10){\circle*{3}}

\put(300,-20){$\scriptstyle K+L=1,\,  \beta_1(\gamma)= 0$}

\end{picture}
\eeq
\ms

Examples for the first type of graphs, which contain an internal vertex,
are $\gamma_1$, $\gamma_2$, and $\gamma_4$ with
$\beta_1=2$, $\beta_1=0$, and $\beta_1=1$, respectively.
The graph $\gamma_3$ is of the second type with only one edge,
and $\gamma_5$ is represented by simply one external vertex without edges.
Although we will not always need the graphs to be one of the representatives,
we shall always assume below that we have no internal vertices of valency one
(as, e.g., $\gamma_2$). Hence we have $\partial \gamma=
\partial \gamma_s\sqcup\partial \gamma_t\,$, with $|\partial \gamma_s|=K$ and 
$|\partial \gamma_t|=L\,$.
\ms

As for the cobordisms in (\ref{def-cocyc}) we also have an anomaly
for the  Betti-numbers of graphs:
\beq\lll{eq-loop-add}
\beta_1(\gamma_2\circ\gamma_1)\,=\,\beta_1(\gamma_1)\,+\,
\beta_1(\gamma_2)\,+\,\,\mu_0(\gamma_2,\gamma_1)\qquad.
\eeq
Here, $\mu_0$ is defined exactly as in Lemma~\ref{lm-cocy=dim}, where the
$W_j$ are now given with respect to the inclusions 
$\partial\gamma_{s,t}\hookrightarrow \gamma\,$. There is no 
$\mu_{\partial}$-contribution, and $\be_1=\beta_1$ for graphs, 
since $H_1(\partial\gamma)=0$.
For the composite of two connected graphs over $k\geq 1$ end-points, we
obtain as in (\ref{eq-mu0-elprod}) that $\mu_0=k-1\,$.

It will be convenient to introduce a natural partial order on the  
morphisms of \G.
For two graphs, $\gamma_1$ and $\gamma_2$, we say that
\beq\lll{eq-paror-gr}
\gamma_1 \;\prec\;\gamma_2\qquad,
\eeq
iff the $\gamma_j$ belong to the same morphism set, and 
 there is an embedding $\gamma_1\hookrightarrow\gamma_2$ of some representatives, such that the corresponding
maps $H_0(\gamma_1)\isto H_0(\gamma_2)\,$ and 
$H_1(\gamma_1)\hookrightarrow H_1(\gamma_2)\,$ are an isomorphisms and a monomorphism, respectively. 
It is clear that $\gamma_2$ is obtained by adding internal edges to a given
component of $\gamma_1$. I.e., for representatives as in (\ref{fig-gr-reps})
$\gamma_2$ differs from $\gamma_1$ only by adding loops to the internal 
vertices.
An inequality as in (\ref{eq-paror-gr}) 
 obviously also implies $\beta_1(\gamma_1)<\beta_1(\gamma_2)\,$, $\gamma_1\circ\gamma\prec\gamma_2\circ\gamma\,$,
as well as $\gamma_1\otimes\gamma\prec\gamma_2\otimes\gamma\,$.

\prg{3.4}{ Coordinate Graphs of Cobordisms }

The relation between the categories $\Cb *$ and $\G$ cannot be given precisely
by a functor because the composition anomaly, $\mu_0+\mu_{\pi}$, of
$\Cb *$ is greater than the anomaly $\mu_0$ of $\G$. We shall, however,
attempt to relate the morphisms of the two categories in a way that we
can conclude from the decomposition of a graph also the decomposition 
of a cobordism. 
We begin with the definition of the notion of a 
(faithful) coordinate graph, which  will be our principal tool in the
description of the connectivity of cobordisms. 

\begin{defi}\lll{def-graph-coor}
A \ub{coordinate-graph} of a cobordism $M:\,\Sigma_s\to\Sigma_t\,$, 
is a graph $\gamma\in\G\,$, together
with  a continuous function
\beq\lll{eq-coord-gr}
f\,:\;(M,\,\Sigma_s,\,\Sigma_t)\,\longrightarrow\,
(\gamma,\,\partial \gamma_s,\,\partial \gamma_s)\;,
\eeq
such that we have a one-to-one correspondence between boundary components,
(i.e., the induced maps $H_0(\Sigma_{s/t})\isto H_0(\partial \gamma_{s/t})\,$
are isomorphisms), and the interior of $M$ is also mapped to the interior
of $\gamma$. 
\smallskip

We say that $f:\, M\to\gamma\,$ is a \ub{faithful} coordinate graph, if there
is  an embedding:
$$
{\it J}\,:\;(\gamma,\partial \gamma)\,\hookrightarrow\,(M,\partial M)\quad,
$$
which maps again connected components to each other, and for which
the composite $f\circ {\it J}$  is homotopic to the identity
(with fixed endpoints).
\end{defi}

An immediate consequence of the correspondence of boundaries is that we
can always write the morphisms of a coordinate graph as unions 
$M\,=\,M'\sqcup N\,$ and $\,\gamma\,=\,\gamma'\sqcup \kappa\,$, where
$N$ and $\kappa$ have no boundaries, and the coordinate map $f$ induces an 
isomorphism
$$
H_0(M')\,\isto\,H_0(\gamma')\;.
$$
For this reason we shall often consider only the case of connected cobordisms 
with connected coordinate graphs.

Next, let us state some obvious facts about the composition and collapse of
coordinate graphs:

\blm\lll{lm-coor-elem}\

\ben
\item If $M_1$ and $M_2$ have (faithful) coordinate graphs $\gamma_1$ and
$\gamma_2$, respectively, then $\gamma_2\circ\gamma_1$ is a (faithful)
coordinate graph of $M_2\circ M_1$.
\item
Suppose $\gamma=\gamma_2\circ\gamma_1$, where the $\gamma_j$'s are 
coordinate graphs
of the $M_j$'s as above, is a maximal (faithful) coordinate
graph of $M=M_2\circ M_1$. Then  $\gamma_j$ has to be a maximal (faithful)
coordinate graph of $M_j$, for both $j=1,2$.

\item If $\gamma$ is a coordinate graph of $M$ and $\xi:\gamma\to\gamma^:$ a
continuous map that preserves enpoints, then $\gamma^:$ is also a coordinate
graph of $M$.

If $\gamma$ is in addition faithful and there is an inclusion 
$\gamma^:\hookrightarrow\gamma\,$, whose composition with $\xi$ is homotopic
to the identity on $\gamma^:$, then $\gamma^:$ is also faithful. 
\een
\elm

Note that the converse of Part 2 is not true. A typical
application of the observation in part {\em 3.)} is given when 
$\gamma^:\prec\gamma\,$ is a subgraph, missing one edge of $\gamma$, 
and the map $\gamma\to\gamma^:$ given by collapsing the additional edge
into another path in $\gamma^:\,$, as for example in the following picture:

\beq\lll{fig-collapse}
\begin{picture}(220,60)

\put(50,5){\line(0,1){50}}

\put(65,27){$\into{30}$}

\put(110,5){\line(0,1){50}}
\put(110,15){\circle*{3}}
\put(110,45){\circle*{3}}

\put(110,35){\oval(12,20) [rt]}
\put(110,25){\oval(12,20) [rb]}
\put(116,25){\line(0,1){10}}

\put(130,27){$\ONTO{$\xi$}{30}$}

\put(180,5){\line(0,1){50}}
\put(180.5,15){\line(0,1){30}}
\put(180,15){\circle*{3}}
\put(180,45){\circle*{3}}

\end{picture}
\eeq

If $\gamma$ is as in (\ref{fig-gr-reps}) and $\gamma^:\prec \gamma$ is of the
same form with  $\beta_1(\gamma)-k$ inner loops, then $\xi$ may be defined 
by collapsing the $k$ outer loops of $\gamma$ to the interior vertex. 

Next, we assure the existence of (maximal and minimal) coordinate graphs.

\blm\lll{lm-coord-ex}\

\ben
\item Every cobordism, $M\,$, admits a (faithful) coordinate graph, $\gamma^M\,$,
which is minimal among all (faithful) coordinate graphs.\nopagebreak
\item Every cobordism has a maximal, faithful coordinate graph.
\een
\elm

{\em Proof:\ }  For a connected cobordism, $M$, we can choose $\gamma^M$,
to be the spider with $\beta_1(\gamma^M)=0$ as in Section 3.2. To define
the coordinate map choose a map 
$g:\Sigma=\Sigma_s\sqcup\Sigma_t\to\partial\gamma\,$, which maps different
components to different points. Let $p:\,M\to v\,$ be the constant map to a point $v\,$. We set
$$
f\,:\;M\,\cong\,M\sqcup_{\Sigma\times\{0\}}\bigl(\Sigma\times [0,1]\bigr)\,
\,\TO{$p\sqcup (g\times id_{[0,1]})$}{70}\,\,v\sqcup_{\sim}
\bigl(\partial\gamma^M\times [0,1]\bigr)\,\cong\,\gamma^M\quad,
$$
where we identified $\partial\gamma^M\times\{0\}\sim v\,$.

In the proof of  the second part it is clear that a faithful coordinate graph with
$\beta_1(\gamma)=k\,$, implies a surjection $H_1(M)\onto{20}{\bf Z}^{k}\,$. For
a compact $M$ we know, e.g., from a Heegaard decomposition, that $H_1(M)$ is
finitely generated so that $N$ must be bounded.
\ep

In Definition~\ref{def-graph-coor} we assumed that for a faithful coordinate 
graph, $h=f\circ J$ is only homotopic to the identity. The next lemma
asserts that we may  assume that in this case $h$ is also equal to the
identity. 

\blm\lll{lm-gen-id}
Suppose $f:\,M\to\gamma\,$ is a generic faithful coordinate graph with 
embedding $J:\,\gamma\,\into{20}M\,$.

Then there exists $f^{\$}:\,M\to\gamma\,$, such that
$$
f^{\$}\circ J\;=\;id\qquad,
$$
and $f^{\$}\,$ coincides with $f$ outside a neighborhood of $J(\gamma)\,$.
\elm

The proof, although fairly standard, is rather technical and is thus deferred
to  Appendix A.1.

An application lies in the proof of the following lemma, asserting
 that if a coordinate graph is a composite, so is the 
associated cobordism. 

\blm\lll{lm-decomp-ogr}
Let $f:\,M\to\gamma\,$ be a generic coordinate graph of a connected cobordism, $M$, and $\gamma=\gamma_2\circ\gamma_1\,$ a decomposition in $\G$. 
\ben
\item There is a graph, $\hat\gamma\succ\gamma$, with a collapse map,
$c:\,\hat\gamma\to\gamma\,$, such that 
$\hat\gamma\,=\,\hat\gamma_2\circ\hat\gamma_1\,$, and 
$c(\hat\gamma_j)\,=\,\gamma_j\,$. 

Moreover, $\hat\gamma$ has the property that there exists a coordinate map 
$\hat f:\,M\to\hat\gamma\,$ -  
with $f\,=\,c\circ\hat f\,$ - such that the 
$M_j\,:=\,\hat f^{-1}(\hat\gamma_j)\,$ are cobordisms with graphs $\hat\gamma_j$
and $M=M_2\circ M_1\,$.

\item If a composed coordinate graph, $\gamma=\gamma_2\circ\gamma_1$, 
is faithful and $\gamma_2$ and 
$\gamma_1$ are connected, then there exists a coordinate graph,
$\hat\gamma\succ\gamma$, and a map,
$\hat f$, as above (except that $f$ may be different from $c\circ\hat f\,$), 
such that the $\hat\gamma_j$ are also faithful, the $M_j$ are connected,
and the embedding of $\hat\gamma_j$ extends that of $\gamma_j\,$.
\een
\elm

{\em Proof: } 
We shall prove here only the first part of the lemma. Due to its technical 
nature the proof of the second statement is again deferred to Appendix A.1. 
It relies heavily on Lemma~\ref{lm-gen-id}.

By genericity we may assume that $f$ is differentiable, and that
$P=\gamma_2\cap\gamma_1$ consists only of regular values of $f$.
For a point, $p\in P\,$, we can find a neighborhood 
$\bigl(I_p,p)\cong\bigl([-\varepsilon,\varepsilon],0\bigr)\,$, such that
$$
f^{-1}\bigl(I_p\bigr)\;=\;
\Sigma_1^p\times[-\varepsilon,\varepsilon]\sqcup
\ldots\sqcup\Sigma_{n_p}^p\times[-\varepsilon,\varepsilon]\quad,
$$
where each $\Sigma_j^p\,$ is a connected Riemann surface, and $f$ acts on this
space as the projection on the interval $[-\varepsilon,\varepsilon]\,$.
 
In order to obtain a coordinate-graph, for each $p\in P$ we insert into $\gamma$
$n_p-1$ additional edges, $e^p_j\cong[-\varepsilon,\varepsilon]\,$, identifying
their boundary points, $\{\pm\varepsilon\}\,$, with those in $\partial I_p\,$.
For the resulting graph $\hat\gamma$ we also have a coordinate map 
$\hat f:\,M\to\hat\gamma\,$, which maps 
$\Sigma_j^p\times[-\varepsilon,\varepsilon]\to e_j^p\,$ for $j=1,\ldots,n_p-1\,$,
through a  projection onto $[-\varepsilon,\varepsilon]$, and 
which coincides with $f$ outside these regions. It is clear that $\hat f$ is 
continuous, and if we define a collapse $c:\,\hat\gamma\to\gamma\,$ by mapping 
$e^p_j\,$ to $I^p$ with fixed endpoints, then $c\circ\hat f\,=\,f\,$.

Now, $M_j\,=\,f^{-1}(\gamma_j)\,$, for $j=1,2$, are cobordisms with $M=M_2
\circ M_1\,$, where we glue over the union of all $\Sigma_j^p\,$. For the 
subgraphs $\hat\gamma_j\,=\,\hat f(M_j)\subset \hat\gamma\,$ we have of 
course $\hat\gamma=\hat\gamma_2\circ\hat\gamma_1\,$, and a one-to-one
correspondence between the set of boundary points of, e.g., $\hat \gamma_1\,$
and the set of $\Sigma_j^p$'s, since $(\partial\hat\gamma_1)_t\,$ 
 contains in addition to the points  of $(\partial\gamma_1)_t\,$ 
the interior points  $0\in e^p_j\,$. The endpoints of $\hat\gamma_1$ 
hence correspond to the boundary components of $M_1$ so that 
$\hat f:\,M_1\to\hat\gamma_1\,$ is a coordinate graph.\ep

The following are useful applications of the decomposition along a graph:

\begin{cor}\lll{cor-3}\

\ben
\item Every cobordism is given by the composite of cobordisms, whose maximal
faithful coordinate graphs are spiders (with at most three end-points).
\item If $M$ has a faithful coordinate graph  with $\beta_1(\gamma)=r$, then
$M$ admits an $r$-diagram (see Section 3.1).
\een
\end{cor}

\noindent
\parbox{4in}{{\em Proof: } Pick a maximal (see Lemma~\ref{lm-coord-ex}) faithful
coordinate-graph $\gamma$ of $M$, and write it as a composite of two trees
over $\beta_1(\gamma)+1$ points as in the diagram to the right.

Since $\gamma$ is already maximal we only have to vary the 
coordinate map in order find to the corresponding decomposition for the cobordisms according to Part {\em 2)} of Lemma~\ref{lm-decomp-ogr}.
By Lemma~\ref{lm-coor-elem} the sub-graphs  are also maximal.}
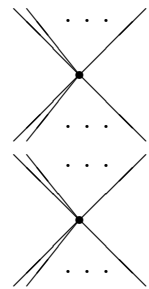
\qquad\quad\qquad\quad\parbox{2in}
{\begin{picture}(60,120)

\put(30,90){\line(-1,1){25}}
\put(30,90){\line(-4,5){20}}
\put(30,90){\line(1,1){25}}
\put(24,110){.\ .\ .\ }
\put(30,90){\circle*{3}}
\put(24,70){.\ .\ .\ }
\put(30,90){\line(-1,-1){25}}
\put(30,90){\line(-4,-5){20}}
\put(30,90){\line(1,-1){25}}

\put(30,35){\line(-1,1){25}}
\put(30,35){\line(-4,5){20}}
\put(30,35){\line(1,1){25}}
\put(24,55){.\ .\ .\ }
\put(30,35){\circle*{3}}
\put(24,15){.\ .\ .\ }
\put(30,35){\line(-1,-1){25}}
\put(30,35){\line(-4,-5){20}}
\put(30,35){\line(1,-1){25}}

\end{picture}
}\vspace*{.2cm}

The inner vertex of each tree can of course be resolved, such that
we obtain a homotopic tree with vertices that have valencies of at most three.
Repeated application of Part {\em 2)} of Lemma~\ref{lm-decomp-ogr} yields
then the decomposition into elementary pieces.

If we use the decomposition over $r+1$ points as above and reconnect 
the two cobordisms over only one surface we obtain a connected manifold. Hence
the remaining surfaces form an $r$-diagram in $M$.
\ep

\medskip

\prg{3.5}{Existence of Coordinate Graphs from Interior Groups}
  
The existence if faithful coordinate graphs and - hence $r$-diagrams - can 
be derived from projections onto the fundamental groups of graphs.
The first observation  regarding this connection follows immediately by 
picking transversal, closed paths in $M$, that represent preimages of 
the generators of $\pi_1(\gamma)\,$:

\blm\lll{lm-onto-faith}
Suppose $\,f:\,M\to\gamma\,$ is a coordinate graph of connected $M$, and 
$$
\pi_1(f)\,:\;\pi_1(M)\,\onto{30}\,\pi_1(\gamma)\;
$$
is onto. 
\nopagebreak
Then $\gamma$ is also a faithful coordinate graph.
\elm

Note, that since $f$ is constant on the boundaries, $\pi_1(f)\,$, also
factors through $\pi_1(M)\to\pint(M)\,$. Hence, with 
$\pi_1(\gamma)\,=\,F\bigl(\beta_1(\gamma)\bigr)\,$ we have that 
$\pi_1(\gamma)\,$ is in fact a free interior group in the sense of 
(\ref{eq-free-int}) of Section 3.2. Next, we show that, conversely,
a free interior group also implies the existence of a coordinate graph.

\blm\lll{lm-ex-coord} Suppose that for connected $M\,$
$$
\zeta\,:\;\pint(M)\,\onto{30}\,F(k)\cong\pi_1(\gamma)\;
$$
is a free interior group, with $k=\beta_1(\gamma)\,$.

Then there exists a faithful coordinate graph $\,f:\,M\to\gamma\,$, such that 
$\zeta=\pi_1(f)\,$.
\elm

{\em Proof:} Assume that $\gamma$ is as in diagram (\ref{fig-gr-reps}), and let
$\gamma^{\&}\subset\gamma\,$ be the bouquet of circles without the exterior edges
so that $\pi_1(\gamma)\,=\,\pi_1(\gamma^{\&})\,$.  Since 
$\pi_j(\gamma^{\&})\,=\,0$ for $j\geq 2\,$ it follows, e.g., from 
Theorem 6.39.ii) in [Sz], that there is a continuous map 
$f^{\&}:\,M\to\gamma^{\&}\,$, which induces the map $\pi_1(M)\to\pint(M)
\ONTO{$\zeta$}{25}\gamma^{\&}\,$. By construction we have 
$\pi(f^{\&}_{\partial})=0\,$ for the restriction 
$f^{\&}_{\partial}:=\,f^{\&}\Bigl |{\partial M}\,$, which in turn implies that
$f^{\&}_{\partial}\,$ is homotopic to the constant map $\partial M\to\{v\}\,$,
where $v\in\gamma^{\&}\,$ is the interior vertex of the bouquet.
Let $F^{\&}_{\partial}:\,\partial M\times [0,1]\,\to\,\gamma^{\&}\,$ be a 
corresponding homotopy, with $F^{\&}_{\partial}(x,0)= f^{\&}_{\partial}(x)\,$
and $F^{\&}_{\partial}(x,1)=v\,$. Given $\psi:\,\partial M\times [0,1]
\sqcup_{\partial M\times\{0\}}M\,\cong \,M\,$, whose restriction to $M$ is 
isotopic to the identity, we can define a function:
$$
f^c\,=\;\bigl(F^{\&}_{\partial}\sqcup_{f^{\&}_{\partial}}f^{\&}\bigr)\,\circ\,
\psi^{-1}\,\,:\;\,\bigl(M,\partial M\bigr)\,
\longrightarrow\, (\gamma^{\&},v)\;.
$$
As in the proof of Part {\em 1)} of Lemma~\ref{lm-coord-ex} we can define 
from this (using again $\psi$) a coordinate map $f:\,M\to\gamma\,$,
with $\pi_1(f)=\pi_1(f^c)\,$. Together with Lemma~\ref{lm-onto-faith} this 
implies the assertion. \ep

Let us summarize in the next theorem  the results of 
Lemma~\ref{lm-comp-freefac}, Corollary~\ref{cor-3}, 
and Lemma~\ref{lm-ex-coord}:

\btm\lll{thm-spc}

Suppose $M$ is a connected cobordism. Then the following are equivalent:
\ben
\item $M$ admits an $r$-diagram.
\item $M$ admits a free interior group of rank $r$.
\item $M$ has a faithful coordinate graph, $\gamma$, with $\beta_1(\gamma)=r\,$.
\een
\etm

In particular Theorem~\ref{thm-spc} implies the converse inequality of 
Part {\em 2)} of Corollary~\ref{cor-2}. We find 
\beq\lll{eq-rho-phi}
\varrho(M)\;=\;\varphi(M)\;,
\eeq
i.e., the maximal number of non-separating surfaces is given by the rank of 
the maximal free interior group of a cobordism. It is thus possible to 
compute the order in {\tt x} of a linear map associated to a cobordisms, $M$,
by a half-projective TQFT as it appears in Lemma~\ref{lm-ord-rho} using only
the fundamental group of $M$ and $\partial M$.

If ${\tt x}=0$ it suffices to consider only homology, since we only have to 
know whether there is a non-trivial interior group or not. More precisely,
with $\H$ and $\be$ as defined in Appendix A.2, we have:
 
\blm\lll{lm-rho=beta}
For a connected cobordism $M\in\Cb 0\,$,
$$
\be(M)=0\qquad{\rm if\ and\ only\ if\ }\qquad \varrho(M)=0\quad.
$$
\elm

{\em Proof:} By naturality of the Hurewicz map the first square of the 
following diagram commutes.

\beq\lll{diag-Hur}
\bar{cccccc}
\vspace*{.2cm}&\pi_1(\partial M)&\longrightarrow&\pi_1(M)&\onto{30}&\pint(M)\\
\vspace*{.2cm}&\Biggl\downarrow\mkern -10mu {\biggl\downarrow}
&&\Biggl\downarrow\mkern -10mu {\biggl\downarrow}&&\Biggl\downarrow\alpha\\
&H_1(\partial M,{\bf Z})&\longrightarrow&H_1(M,{\bf Z})&\onto{30}&\H(M,{\bf Z})
\ear
\eeq

Since the lower sequence is exact, and the Hurewicz maps are surjective, we
can infer the existence of a surjection $\alpha\,$ of the interior groups, 
such  that all of (\ref{diag-Hur}) commutes. From  $\H(M,{\bf Q})\neq 0\,$
we know that there is an epimorphism $\H(M,{\bf Z})\onto{20}{\bf Z}\,$, which,
composed with $\alpha$, gives rise to a free interior group pf rank one.

Conversely, if there is an epimorphism $\pint(M)\onto{20}\,{\bf Z}\,$,
the corresponding map $\pi_1(M)\onto{20}\,{\bf Z}\,$ has the commutator
sub-group in its kernel and thus factors into homology. Since the kernel also
contains the image of $\pi_1(\partial M)\,$, it follows from  diagram 
(\ref{diag-Hur}) that the epimorphism on homology factors through $\H$ so that
$\be(M)\geq 1\,$. \ep

We immediately find from this and Lemma~\ref{lm-ord-rho} the following.

\begin{cor}\lll{cor-O}
Suppose \V\ is a half-projective TQFT w.r.t. ${\tt x}=0$ and $\mu_0$, and let
$M\in\Cb 0\,$ is a cobordism with $\be (M)\neq 0\,$. Then 
$$
\V(M)\,=\,0\quad.
$$
\end{cor}

In the special case of the ``Hennings-invariant'' for $U_q(s\ell_2)\,$,
with $q$ at a root of unity, this vanishing property was observed by 
Ohtsuki, see [O]. The result there, however, is found  more-or-less from
a direct computation of the invariant.

From the discussion in Chapter 2 we have seen that ${\tt x}=0$ and 
invertible {\tt x} are probably the only possibilities so that
- from an algebraic point of view - we only have to worry whether 
there are non-trivial interior groups as, e.g.,  in Lemma~\ref{lm-rho=beta}, but 
not about the exact order, as suggested in Lemma~\ref{lm-ord-rho}. Still,
as we shall discuss in Section 4.5, the precise order can be of interest,
if we consider ``classical limits'' of TQFT's. In this case we  have 
${\tt x}\to\infty$, and $\varrho(M)$ may yield an estimate on the order
in {\tt x}, by which $\|\V(M)\|$ diverges.

\cht{4}{Construction of Half-Projective TQFT's}

The aim of this chapter is to show how half-projective TQFT's can be 
constructed from connected ones, as for example those found in [KL].
We shall organize our discussions in a deductive way, including 
additional assumptions only where needed. Specifically,  we shall
begin in Section 4.2 by introducing a set of Axioms, V1-V5,
 on a map $\V:\Cb 0\to {\tt R}-mod$, 
extract from this a list of properties, P1-P8,
that ensure the  existence of such a map, and conclude that \V\ 
necessarily has to be a 
half-projective TQFT. In Section 4.3 we show that the existence of
extended structures implies all but one of the properties automatically.
The missing Property P7 on the projectivity {\tt x}
is discussed in Section 4.4 as a consequence of the closely related 
concepts of cointegrals, semisimplicity, and the invariant on
$S^1\times S^2\,$. In the discussions of these sections we attempt to give
a clear picture of how the given assumptions influence the existence
and uniqueness of the properties we derive for $\V$. The last section then
summarizes the possible diversions from our axioms
that might lead to more general
definitions of $\V$, in particular the tensor product rule.
We also discuss a possible application of the formalism of half-projective
TQFT's to the study of ``classical limits'', where the exact orders
$\varrho(M)=\varphi(M)$ become relevant.

\prg{4.1}{Surface-Connecting Cobordisms}

For two closed, connected Riemann 
surfaces $\Sigma_j\,$, with $j=1,2$, we can think 
of their connected sum $\Sigma_1\#\Sigma_2\,$ as being the result of a 
1-surgery on  $\Sigma_1\sqcup\Sigma_2\,$, i.e., we cut away a disc from each
surface and reglue the cylinder $S^1\times I\,$ along the boundaries.
The corresponding morphism 
$\Pi:\,\Sigma_1\sqcup\Sigma_2\,\to\, \Sigma_1\#\Sigma_2\,$ is constructed by  
attaching a one-handle to the cylinder over $\Sigma_1\sqcup\Sigma_2\,$, i.e.,
$$\Pi\,=\;\Sigma_1\times[0,1]\,\sqcup_{D^2\times\{1\}}D^2\times I 
\sqcup_{D^2\times\{1\}}\,\Sigma_2\times [0,1]\;,
$$
or, equivalently, a boundary-connected sum of the $\Sigma_j\times[0,1]\,$.
More generally, we obtain for every surface $\Sigma$ with $K$ ordered 
connected components $\Sigma_j\,$ a morphism
$$
\Pi_{\Sigma}\,:\;\Sigma=\Sigma_1\sqcup\ldots\sqcup\Sigma_K\,\,\longrightarrow\,\,
\Sigma^{\#}:=\Sigma_1\#\ldots\#\Sigma_K\;.
$$
The obvious associativity condition for these morphisms is readily verified.
Changing orientations we obtain a cobordism in the opposite direction:
$$
\Pi_{\Sigma}^{\dagger}\,=\,-\Pi_{\Sigma}\,:\;\Sigma^{\#}=
\Sigma_1\#\ldots\#\Sigma_K
\,\,\longrightarrow\,\,\Sigma=\Sigma_1\sqcup\ldots\sqcup\Sigma_K\;.
$$
In the next lemma we evaluate the composites of $\Pi_{\Sigma}\,$ and 
$\Pi_{\Sigma}^{\dagger}\,$:
\blm\lll{lm-comp-Pi}
Suppose $\Sigma$ is a closed 
Riemann surface with $K$ components, $\Sigma_j\,$,
and the cobordisms $\Pi_{\Sigma}\,$ and $\Pi_{\Sigma}^{\dagger}\,$ are as 
above. Then we have
\ben
\item
$$
\Pi_{\Sigma}^{\dagger}\circ\Pi_{\Sigma}\;=\;\id_{\Sigma_1}\#\ldots
\#\id_{\Sigma_K}\;,
$$
i.e., the (interior) connected sum of the cylinders $\Sigma_j\times[0,1]\,$.
\item For
$$
\Lambda_{\Sigma}:=\Pi_{\Sigma}\circ\Pi_{\Sigma}^{\dagger}\,:\;
\Sigma^{\#}\,\,\longrightarrow\,\,\Sigma^{\#}\;.
$$
we have 
$$
\Lambda_{\Sigma}\circ\Pi_{\Sigma}\;=\;
\Pi_{\Sigma}\#\underbrace{(S^1\times S^2)\#
\ldots\#(S^1\times S^2)}_{K-1\,{\rm\ times}}
$$
and the analogous equation for $\Pi_{\Sigma}^{\dagger}\circ\Lambda_{\Sigma}\,$.
($\#$ is the interior connected sum of 3-folds).
\een
\elm

{\em Proof:} For the proof of the first part it is enough to consider only
the case $N=2$. If we think of 
$\Pi=\Sigma_1\times[0,1]\cup_{D^2}\Sigma_2\times[0,1]\,$ as a 
boundary-connected sum of two parts, the composite $\Pi^{\dagger}\circ\Pi\,$
is glued together from four pieces. The cylinders 
$\Sigma_1\times[0,1]\subset\Pi\,$ and 
$\Sigma_1\times[1,2]\subset\Pi^{\dagger}\,$ are glued together along 
$\Sigma_1\times\{1\}\,$ except at the disc $D^2\subset \Sigma_1\,$.
Their union is thus homeomorphic to $\bigl(\Sigma_1\times[0,1]\bigr)-D^3\,$,
where the union $D^2\cup_{S^1}D^2\cong S^2_1\,$ of the discs,
 that are not identified, bounds 
the removed ball $D^3\,$. In the same way the cylinders over $\Sigma_2\,$
are glued to give $\bigl(\Sigma_2\times[0,1]\bigr)-D^3\,$ with an additional 
boundary component, $S^2_2\,$. In order to find the total composite, we have 
to glue the discs together as in the definition of $\Pi$, which amounts to 
gluing $S^2_1$ onto $S^2_2\,$. This shows the first part of the lemma.

For the second part we find
$$
\Lambda_{\Sigma}\circ\Pi_{\Sigma}\,=\;
\Pi_{\Sigma}\circ\Bigl(\bigl(\id_{\Sigma_1}-D^3\bigr)\cup\bigl(S^2\times I\bigr)\cup
\bigl(\id_{\Sigma_2}-D^3\bigr)\cup\ldots\cup\bigl(S^2\times I\bigr)\cup
\bigl(\id_{\Sigma_K}-D^3\bigr)\Bigr)\;,
$$
where we rewrote the connected sum as a (3-dimensional)
index-1-surgery on the union of the cylinders.
Clearly, we can view these in the composite also as
index-1-surgeries on 
$\Pi_{\Sigma}\,\cong\,
\Pi_{\Sigma}\circ\bigl(\id_{\Sigma_1}\sqcup\ldots\sqcup\id_{\Sigma_K}\bigr)
\,$. Since  $\Pi_{\Sigma}\,$ is connected,
the surgery-points can be moved together without changing the homeomorphism
type of $\Lambda_{\Sigma}\circ\Pi_{\Sigma}\,$.
The assertion follows now from the fact that an index-1-surgery in  
a contractible neighborhood is the same as connected summing with 
$S^1\times S^2\,$. The proof for $\Pi_{\Sigma}^{\dagger}\circ\Lambda_{\Sigma}\,$
is analogous.
\ep

Besides the formulas in Lemma~\ref{lm-comp-Pi} and associativity we 
shall consider  another type of relations among the connecting morphisms.
The basic example is given next:

\blm\lll{lm-comm-Pi}
Suppose $\Sigma_1\,$, $\Sigma_2\,$, and $\Sigma_3\,$ are closed, 
connected surfaces.
Then
\beq\lll{eq-comm-Pi}
\Bigl(\Pi_{(\Sigma_1\sqcup\Sigma_2)}\sqcup\id_{\Sigma_3}\Bigr)\circ
\Bigl(\id_{\Sigma_1}\sqcup\Pi_{(\Sigma_2\sqcup\Sigma_3)}^{\dagger}\Bigr)\;=\;
\Pi^{\dagger}_{\bigl((\Sigma_1\#\Sigma_2)\sqcup\Sigma_3\bigr)}\,\circ\,
\Pi_{\bigl(\Sigma_1\sqcup(\Sigma_2\#\Sigma_3)\bigr)}\quad.
\eeq
\elm

{\em Proof:} Denote by $\Sigma_{1/3}^{\bullet}=\Sigma_{1/3}-D^2\,$, and
$\Sigma_2^{\bullet\bullet}=\Sigma_2-\bigl(D^2\sqcup D^2\bigr)\,$,
with $\partial\Sigma_2^{\bullet\bullet}=S^1_L\sqcup S^1_R\,$,
the corresponding surfaces with holes, such that, e.g., 
$\Sigma_1\#\Sigma_2=\Sigma_1^{\bullet}\cup_{S^1_L}\Sigma_2^{\bullet\bullet}
\cup_{S^1_R} D^2\,$. The morphism $\Pi_{\Sigma_1\sqcup\Sigma_2}\,$, for
example, can then be seen as 
$\Bigl(\Sigma_1^{\bullet}\cup_{S^1_L}\Sigma_2^{\bullet\bullet}
\cup_{S^1_R} D^2\Bigr)\times [1,2]\,$ with a 2-cell,
$C_{12}\cong D^2\times I\,$, attached along the connecting $S^1_L$ of the 
source surface. Similarly, $\Pi^{\dagger}_{\Sigma_2\sqcup\Sigma_3}\,$ is
$\Bigl(D^2\cup_{S^1_L}\Sigma_2^{\bullet\bullet}
\cup_{S^1_R} \Sigma_3^{\bullet}\Bigr)\times [0,1]\,$, with a 2-cell, $C_{23}$,
glued along the $S^1_R\,$ of the target surface. We construct the composite
by first gluing the pieces $\Bigl(\Sigma_1^{\bullet}\cup_{S^1_L}\Sigma_2^{\bullet\bullet}\Bigr)\times I$ and $\Bigl(\Sigma_2^{\bullet\bullet}
\cup_{S^1_R} \Sigma_3^{\bullet}\Bigr)\times I\,$ together along the 
respective $\Sigma_2^{\bullet\bullet}$ boundary pieces. The result is 
homeomorphic to the cylinder over
$\Sigma_1^{\bullet} \cup\Sigma_2^{\bullet\bullet}\cup\Sigma_3^{\bullet}\,
\cong\,\Sigma_1\#\Sigma_2\#\Sigma_3\,$. Attaching the remaining pieces, the
2-cell $C_{12}\,$ combines with $D^2\times[0,1]\subset
\Pi^{\dagger}_{\Sigma_2\sqcup\Sigma_3}\,$ to one thickened disc, as does 
$C_{23}$ with the $D^2\times I$-piece of $\Pi_{\Sigma_1\sqcup\Sigma_2}\,$.
The composite is thus $\Sigma_1\#\Sigma_2\#\Sigma_3\times I$ with a 2-cell
attached along the $S^1_R$ to the target surface and another 2-cell glued
along $S^1_L$ to the source surface. If we split the middle cylinder of
the cobordisms into two, the result is precisely the composite on the right of
(\ref{eq-comm-Pi}).\ep

The connecting morphisms allow us to express a general, connected 
cobordism by a cobordism  between connected surfaces. For this purpose let
us introduce the notations
$$
\Cc * \;\subset \;\Cb *\qquad\quad{\rm and}\qquad\quad\Cc 0 \;\subset \;\Cb 0
$$
for the subcategories, whose objects are {\em connected} surfaces. 
We have the following general presentation:

\blm \lll{lm-conmor-dec}

For any connected cobordism  $\,M:\,\Sigma_s\to\Sigma_t\,$ in $\Cb 0$ there 
exists a morphism, $\,\hat M:\,\Sigma_s^{\#}\to\Sigma_t^{\#}\,$, in \Cc 0 , 
such that 
$$
M\;=\;\Pi^{\dagger}_{\Sigma_t}\circ\hat M\circ \Pi_{\Sigma_s}\;.
$$
\elm

{\em Proof:} The proof is immediate from tangle presentations of 
cobordisms as in Appendix A.3 or [Ke2] . 
A  direct proof is given by choosing a Morse
function $f:M\to[0,1]$ with $\Sigma_s=f^{-1}(0)\,$ and  
$\Sigma_t=f^{-1}(1)\,$. It follows from the general theory of stratified
function spaces that  $f$ can be deformed such that it does not have any 
index-zero-singularities, and the 
index-one-singularities have values below all other critical values.
Also the order of the critical values of index-one can be freely permuted
so that we can assume that the one with values in an interval, $[0,\delta]\,$,
are all fusing singularities.
This means the one-handle attachement given by passing through such a singularities
is between different components of the surface. As $M$ is connected there will be
exactly $\beta_0(\Sigma_s)-1$ such handles. The handles of 
$X=f^{-1}\bigl([0,\delta]\bigr)\,$ can be freely slid using isotopies of maps on 
the  upper surfaces $f^{-1}(\delta)\,$ of $X$. Hence we can find a boundary chart,
for which $X$ is equivalent to $\Pi_{\Sigma_s}\,$.
The arguments for splitting off $\Pi^{\dagger}_{\Sigma_t}$ are analogous.
\ep

A useful application of Lemma~\ref{lm-conmor-dec} is the presentation of the
symmetric group action on the $\Pi_{\Sigma}\,$.
As in Lemma~\ref{lm-tens-cob} a permutation of the connected components
of $\Sigma$ can be given in terms of a cobordism $\pi^*\,$. It is not hard
to see that they can be induced from a corresponding braid group action on
the connected sum, $\Sigma^{\#}\,$. More precisely, for $K=\beta_0(\Sigma)\,$
we have a homomorphism
\beq\lll{eq-braid-Pi}
\rho:\,B_K\longrightarrow\pi_0\Bigl({\cal D}i\!f\!f(\Sigma^{\#})^+\Bigr)\;,
\quad\,
{\rm such\ that\ }\qquad \rho(b)\circ\Pi_{\Sigma}\,=\,
\Pi_{\Sigma}\bigl(\overline b)^*\;,
\eeq
where $b\mapsto \overline b\,$ is given by the natural map $B_K\onto{20}S_K\,$.
\medskip

It is quite helpful to consider choices of the connecting morphisms in the 
framework of tangle presentations, see Appendix A.3. 
For simplicity we consider only closed 
surfaces, denoting by $\Sigma_g$ a connected standard surface of genus $g$.

Both  $\Sigma=\Sigma_{g_1}\sqcup\ldots\sqcup\Sigma_{g_K}\,$ and 
$\Sigma^{\#}=\Sigma_{(g_1+\ldots+g_K)}\,$ are represented in a tangle category
by $g_1+\ldots+g_K\,$ pairs of end-points. However, in the first case they are 
organized in $K$ {\em groups} (i.e., we have a $\tau$-move for each group),
and for $\Sigma^{\#}$ all end-points belong to the same group. Indicating groups
by braces and admitting through-strands the connecting morphisms can be presented
as follows:

\beq\lll{fig-conn-mor}
\begin{picture}(420,90)

\put(5,40){$\Pi_{\Sigma}^{\dagger}=$}

\put(40,15){\rule{2.3in}{1mm}}
\put(40,70){\rule{2.3in}{1mm}}

\put(45,15){$\underbrace{\hspace*{1.3cm}}$}
\put(55,-2){$2g_1$}
\put(49,15){\line(0,1){55}}
\put(78,15){\line(0,1){55}}
\put(56,30){$.\,.\,.$}

\put(160,15){$\underbrace{\hspace*{1.3cm}}$}
\put(170,-2){$2g_K$}
\put(164,15){\line(0,1){55}}
\put(193,15){\line(0,1){55}}
\put(171,30){$.\,.\,.$}

\put(92,50){$.\qquad.\qquad.$}
\put(92,4){$.\qquad.\qquad.$}

\put(45,75){$\overbrace{\hspace*{5.4cm}}$}
\put(80,86){$2(g_1+\ldots+g_K)$}

\put(240,40){$\Pi_{\Sigma}=$}

\put(275,15){\rule{2.3in}{1mm}}
\put(275,70){\rule{2.3in}{1mm}}

\put(280,75){$\overbrace{\hspace*{1.3cm}}$}
\put(295,86){$2g_1$}
\put(284,15){\line(0,1){55}}
\put(313,15){\line(0,1){55}}
\put(291,60){$.\,.\,.$}

\put(327,40){$.\qquad.\qquad.$}
\put(327,82){$.\qquad.\qquad.$}

\put(395,75){$\overbrace{\hspace*{1.3cm}}$}
\put(405,86){$2g_K$}
\put(399,15){\line(0,1){55}}
\put(428,15){\line(0,1){55}}
\put(406,60){$.\,.\,.$}

\put(280,15){$\underbrace{\hspace*{5.4cm}}$}
\put(315,-2){$2(g_1+\ldots+g_K)$}

\end{picture}
\eeq

We can reproduce the first assertion of Lemma~\ref{lm-comp-Pi} using the
fact that tangles are composed in the na\"\i ve way, if the intermediate
object consists of only one group. The resulting tangle diagram of 
$\Pi_{\Sigma}^{\dagger}\circ\Pi_{\Sigma}$ consists  thus again
of $2(g_1+\ldots+g_K)\,$ vertical strands. However, now  both the
top and bottom end-points are divided into $K$ groups. 
From the three-dimensional
interpretation of the tangles, the $K-1$ dividing lines of the diagram can be
seen as $K-1$ dividing spheres of the corresponding cobordism. Since
the composites with connected objects are easily identified as the punctured
identity-cobordisms, we find  the connected sum $\id_{\Sigma_1}\#\ldots
\#\id_{\Sigma_K}\,$. 

If we consider the opposite composite, we glue over $K$ components, and thus,
by the rules given in [Ke2], we have to insert $K-1$ zero-framed annuli.
The first will surround the first $2g_1$ strands, the second the next
$2g_2$ strands, etc. Only the last group of $2g_K$ strands is not surrounded
by an annulus, as depicted in the following diagram:

\beq\lll{fig-lambda-mor}
\begin{picture}(220,90)

\put(-10,40){$\Lambda_{\Sigma}=$}

\put(40,15){\rule{2.83in}{1mm}}
\put(40,70){\rule{2.83in}{1mm}}

\put(49,15){\line(0,1){13}}
\put(78,15){\line(0,1){13}}

\put(47,30){\line(1,0){33}}
\put(47,33){\oval(6,6)[l]}
\put(80,33){\oval(6,6)[r]}

\put(49,32){\line(0,1){38}}
\put(78,32){\line(0,1){38}}

\put(56,57){$.\,.\,.$}

\put(164,15){\line(0,1){13}}
\put(193,15){\line(0,1){13}}
\put(164,32){\line(0,1){38}}
\put(193,32){\line(0,1){38}}

\put(162,30){\line(1,0){33}}
\put(162,33){\oval(6,6)[l]}
\put(195,33){\oval(6,6)[r]}

\put(171,57){$.\,.\,.$}

\put(204,15){\line(0,1){55}}
\put(233,15){\line(0,1){55}}
\put(211,57){$.\,.\,.$}

\put(92,50){$.\qquad.\qquad.$}

\put(45,75){$\overbrace{\hspace*{6.8cm}}$}
\put(95,86){$2(g_1+\ldots+g_K)$}

\put(45,15){$\underbrace{\hspace*{6.8cm}}$}
\put(95,-2){$2(g_1+\ldots+g_K)$}

\end{picture}
\eeq

Instead of the last group we could have also chosen any other group as 
the one without an annulus. 

The second assertion of Lemma~\ref{lm-comp-Pi} follows easily from 
(\ref{fig-lambda-mor}), since the  composition of a connecting morphism
from either side has simply the effect of creating $K$ groups, and hence
introduces the $\tau$-move at a particular group. It is clear that the 
$\tau$-move at the $j$-th group with $2g_j$ strands allows us to push the
annulus through  the strands, and thus separate it from the rest of the diagram.
The resulting formula for $\Lambda_{\Sigma}\circ\Pi_{\Sigma}\,$ 
follows now from the
fact that  an isolated, unframed unknot represents a connected sum with 
$S^1\times S^2\,$.
\medskip

Moreover, the relation (\ref{eq-comm-Pi}) in Lemma~\ref{lm-comm-Pi} is readily
verified in the framework of tangles. We find that in both cases the 
composition yields vertical strands going from two groups with $g_1$ and
$g_2+g_3$ pairs to two groups with $g_1+g_2$ and $g_3$ pairs of strands.
\medskip

A cobordism $\hat M$ 
as  in Lemma~\ref{lm-conmor-dec} is also easily found from a 
given tangle, representing $M$, by simply interpreting top and bottom strands
as one group. The ambiguity in choosing $\hat M$ is expressed by the 
additional $\tau$-moves that arise, when we compose with the connecting 
morphisms.
\medskip

Finally, let us illustrate the composites in (\ref{eq-braid-Pi}) for $K=2$  
in the following tangle diagram:

\beq\lll{fig-braid-Pi}  
\begin{picture}(220,90)
\put(0,40){$\rho(\sigma_1)\circ\Pi_{\Sigma}=$}

\put(75,15){\rule{2.3in}{1mm}}
\put(75,70){\rule{2.3in}{1mm}}

\put(80,75){$\overbrace{\hspace*{2.4cm}}$}
\put(110,86){$2g_1$}
\put(84,15){\line(3,2){40}}
\put(159,65){\line(3,2){10}}
\put(144,15){\line(3,2){11}}
\put(189,45){\line(3,2){39}}
\put(110,25){$.\;\;.\;\;.$}

\put(165,75){$\overbrace{\hspace*{2.4cm}}$}
\put(195,86){$2g_2$}
\put(169,15){\line(-3,2){83}}
\put(229,15){\line(-3,2){83}}
\put(173,60){$.\;\;.\;\;.$}

\put(135,44){$.\quad\,.\quad\,.$}

\put(80,15){$\underbrace{\hspace*{5.4cm}}$}
\put(132,-2){$2(g_1+g_2)$}

\end{picture}
\eeq

The cobordism $\rho(\sigma_1)\,$ of the generator $\sigma_1\in B_2\,$
is presented by the crossing of $2g_1$ parallel strands with the the
remaining set of $2g_2$ parallel strands. As  $\rho(\sigma_1)\,\subset\Cc 0\,$,
all $2(g_1+g_2)$ belong to the same group. However, if we multiply
a connecting morphism, they belong to different groups at one end of the
composite. The additional $\tau$-move allows us to change the collective 
overcrossing to a collective undercrossing. Thus 
$\rho(\sigma_1)^2\circ\Pi_{\Sigma}=\Pi_{\Sigma}\,$, and we recover the action
of the symmetric group.

\prg{4.2}{Basic Constraints on Generalized TQFT's}

The purpose of this section is to list a number of natural assumptions on 
a map from 
$\Cb 0$ to the category of {\tt R} modules, such that they necessarily
imply the notion of a half projective TQFT as  the only possible 
generalization over an ordinary TQFT.
We shall also give a list of conditions in terms of some elementary
morphisms and relations that allow us to construct a unique half-projective
TQFT from a given connected one.
\medskip

Let us begin with the axioms we shall require of an assignment
$$
\V\,:\;\Cb 0\,\longrightarrow\,{\tt R}-(free)mod\;,
$$
which maps objects to objects and morphism-sets to the corresponding 
morphism-sets. 
\ben
\item[V1)] $\V(M)\V(N)=\V(M\circ N)\,$, if $M$ and $N$ are connected, and
\newline if they are composed over a \ub{connected} surface.
\item[V2)] $\V$ respects the symmetric tensor structure.
\item[V3)] $dim\Bigl(\V(S^2)\Bigr)\,=\,1\;$.
\item[V4)] $\V(\id_{\Sigma})\;=\;\id_{{\cal V}(\Sigma)}\;$.
\een
For the remainder of the chapter we shall {\em define} for a given \V:
\beq\lll{eq-def-x}
{\tt x}\,\equiv\,{\tt x}_{\cal V}\,:=\;\V(S^1\times S^2)\;.
\eeq
Note that in a (generalized) TQFT with V1 and V2, ${\sf U}:=\V(S^2)\,$ is
a commutative algebra over {\tt R}, which acts on \V\ . I.e., if 
${\sf U}={\sf U}_0\oplus{\sf U}_1\,$, then \V\ also decomposes into a sum
in a similar way in Lemma~\ref{lm-ring-dec}. Hence V3 is, e.g.,
 the same as assuming
that ${\sf U}$ is semisimple and $\V$ is indecomposable.
Furthermore, together with V1 Axiom V3 implies
\beq\lll{eq-V3}
\V(M\#N)\,=\,\V(M\sqcup N)\;.
\eeq

It is, by V2, clearly enough to require V4 only for connected surfaces.
Also, it is easy to see [A] that $\V(\id_{\Sigma})\,$ is in any case a
projector on $\V(\Sigma)\,$, and that by a reduction to 
$\V^{red}(\Sigma)=im\Bigl(\V(\id_{\Sigma})\Bigr)\,$ we can define a 
consistent, effective TQFT $\V^{red}\,$ which obeys V4. Hence the 
last axiom simply assumes that we have already carried out this reduction.

Axiom V1 also implies that we have an honest functor
\beq\lll{eq-conn-TQFT}
\V^{conn}\,:\;\Cc 0\,\longrightarrow\,{\tt R}-mod\;,
\eeq
where $\Cc 0$ is the category of cobordisms between {\em connected}
 surfaces as in described in the previous section. Using the decomposition
in Lemma~\ref{lm-conmor-dec} and the fact that by V1 the composition in
(\ref{lm-conmor-dec}) of the lemma should be respected by $\V$, it follows that
$\V(M)$ for general $M\in\Cb 0\,$ is given by $\V^{conn}(\hat M)\,$ and 
elementary maps, to be associated to the surface connecting cobordisms:
\beq\lll{eq-inj-proj}
\bar{rcccc}
i_{\Sigma}:=\V(\Pi_{\Sigma})\,&:&\;V_{\Sigma_1}\otimes\ldots\otimes 
V_{\Sigma_K}\,&\longrightarrow&\,V_{\Sigma^{\#}}\\
p_{\Sigma}:=\V(\Pi_{\Sigma}^{\dagger})\,&:&\;V_{\Sigma^{\#}}\,
&\longrightarrow&\,V_{\Sigma_1}\otimes\ldots\otimes V_{\Sigma_K}\;,
\ear
\eeq
with $\Sigma=\Sigma_1\sqcup\ldots\sqcup\Sigma_K\,$ and 
 $\Sigma^{\#}=\Sigma_1\#\ldots\#\Sigma_K\,$.

Next, let us derive from V1-V3 a set of constraints on these maps, 
from which we infer the existence of a generalized TQFT for 
disconnected surfaces. To begin with, observe that the composite in Part 1 
of Lemma~\ref{lm-comp-Pi} is over a connected surface so that we have 
$p_{\Sigma}i_{\Sigma}=\V(\id_{\Sigma_1}\#\ldots\#\id_{\Sigma_K})\,$.
By (\ref{eq-V3})
 the connected sums can be replaced by disjoint unions, which yields the
identity. Hence, by V4,
 the $i_{\Sigma}$ and $p_{\Sigma}$ must be injections and 
projections that identify $V_{\Sigma}=\bigotimes_jV_{\Sigma_j}\,$ as a 
complemented subspace in $V_{\Sigma^{\#}}\,$. They also satisfy associativity
constraints as in 
\beq\lll{eq-elemor-assoc}
i_{\Sigma_A\sqcup\Sigma_B}\;=\;i_{\Sigma_A^{\#}\sqcup\Sigma_B^{\#}}
\Bigl(i_{\Sigma_A}\otimes i_{\Sigma_B}\Bigr)\;,
\eeq
where $\Sigma_A$ and $\Sigma_B$ are possibly disconnected surfaces, and 
$\Sigma_A^{\#}$ and $\Sigma_B^{\#}$ are as usual the connected sums of their
components.

Given a connected TQFT, $\V^{conn}\,$, we define
\beq\lll{eq-def-L}
{\bf L}_{\Sigma}\,:=\;\V^{conn}\bigl(\Lambda_{\Sigma}\bigr)\qquad\in\,
End_{\tt R}\bigl(V_{\Sigma^{\#}}\bigr)\quad.
\eeq

If \V\ were an ordinary TQFT, the definition of $\Lambda_{\Sigma}$ in 
Lemma~\ref{lm-comp-Pi} would imply that this map is identical with the
following projector:
\beq\lll{eq-P=ip}
\P_{\Sigma}\,:=\;i_{\Sigma}p_{\Sigma}\qquad\quad
\,\in\,End_{\tt R}\bigl(V_{\Sigma^{\#}}\bigr)\;.
\eeq
However, the composition for $\Lambda_{\Sigma}\,$ is over disconnected $\Sigma$
so that we cannot apply V1. Still, the findings in the second part of 
Lemma~\ref{lm-comp-Pi}  allow us to derive from Axioms V1 and (\ref{eq-V3})
the relations 
${\bf L}_{\Sigma}i_{\Sigma}\,=\,{\tt x}^{K-1}i_{\Sigma}\,$ and 
$p_{\Sigma}{\bf L}_{\Sigma}\,=\,{\tt x}^{K-1}p_{\Sigma}\,$, where
$K=\beta_0(\Sigma)\,$. This means that ${\bf L}_{\Sigma}={\tt x}^{K-1}
\P_{\Sigma}+L^c_{\Sigma}\bigl(\id-\P_{\Sigma}\bigr)$, where 
$[L^c_{\Sigma},\P_{\Sigma}]=0\,$, and 
${L^c_{\Sigma}}^2={\tt x}^{K-1}L^c_{\Sigma}\,$. So, in principle we could have a
composition anomaly of $i_{\Sigma}$ and $p_{\Sigma}$ in the form of an 
operator $L^c_{\Sigma}$ acting on the complement of the space 
$V_{\Sigma}=\bigotimes_jV_{\Sigma_j}\,$, through which the injection and  
projection normally map. At this point we cannot find a  reasonable way to 
incorporate such an anomaly into our formalism, and it does not appear in
the known construction. We shall thus add its absence to the list of axioms:
\ben
\item[V5)]\qquad\qquad\qquad\qquad\qquad\qquad $L^c_{\Sigma}=0$
\een
In the case, where each $\V(\Sigma)\,$ is generated by all $\V(M)\,$,
with $M\,:\;\emptyset\to\Sigma\,$, Axiom V5 can also be inferred
from  the condition that for any collection of $M_j:\,\emptyset\to\Sigma\,$,
$$
\sum_jc_j\V(M_j)=0\qquad\;{\rm implies }\;\qquad 
\sum_jc_j\V^{conn}(\Pi_{\Sigma}\circ M_j)=0
$$
Since the equations in Lemma~\ref{lm-comm-Pi} and in (\ref{eq-braid-Pi})
involve only compositions over connected surfaces, we shall impose by V1 the 
corresponding formulas for the maps $i_{\Sigma}$ and $p_{\Sigma}$ as 
additional conditions. Moreover, suppose that we have cobordisms,
$M_j:\Sigma_j\to\tilde\Sigma_j\,$, for connected surfaces $\Sigma_j$ and 
$\tilde\Sigma_j\,$, and $M:\,\Sigma^{\#}\to\tilde\Sigma^{\#}\,$ for their
connected sums, such that 
\beq\lll{eq-nat-Pi}
\Pi_{\tilde\Sigma}\circ\Bigl(M_1\sqcup\ldots\sqcup M_K\Bigr)\;=\;
M\circ\Pi_{\Sigma}\;.
\eeq
Then V1 and V2 imply that this relation is respected by \V.
\medskip

Let us summarize the conditions on \V\ that we have derived so far from
the Axioms V1-V5, with {\tt x} as in (\ref{eq-def-x}):
\ben
\item[P1)]There is a connected TQFT, $\V^{conn}$, as in (\ref{eq-conn-TQFT}).
\item[P2)] There are injections and projections, 
$V_{\Sigma}\INTO{$i_{\Sigma}$}{25}V_{\Sigma^{\#}}\,$ and 
$V_{\Sigma^{\#}}\ONTO{$p_{\Sigma}$}{25}V_{\Sigma}\,$, with 
$p_{\Sigma}i_{\Sigma}=\id\,$.
\item[P3)]$\tau$-Invariance: $\V(M)=p_{\Sigma_t}\V^{conn}(\hat M)i_{\Sigma_s}\,$
only depends on $M=\Pi^{\dagger}_{\Sigma_t}\circ{\hat M}\circ\Pi_{\Sigma_s}\,$.
\item[P4)] Associativity: See (\ref{eq-elemor-assoc}) and analogously for $p_{\Sigma}$.
\item[P5)] Symmetry: $\V^{conn}\bigl(\rho(b)\bigr)i_{\Sigma}= 
i_{\Sigma}\bigl(\overline b)^*\,$,\quad(see (\ref{eq-braid-Pi})),\newline
\hphantom{xxxxxxxxxx} where $S_K$ acts on $V_{\Sigma}=\bigotimes_jV_{\Sigma_j}\,$ by canonical permutation. (same for $p_{\Sigma}$).
\item[P6)] Naturality: $i_{\tilde\Sigma}\Bigl(\V^{conn}(M_1)\otimes
\ldots\otimes\V^{conn}(M_K)\Bigr)=\V^{conn}(M)i_{\Sigma}\,$, 
\newline \hphantom{xxxxxxxxxxxx} where $M$ and $M_j$ are as in (\ref{eq-nat-Pi}). (same for $p_{\Sigma}$).
\item[P7)] Projectivity: $\V^{conn}\bigl(\Lambda_{\Sigma}\bigr)\,=
\,{\tt x}^{\beta_0(\Sigma)-1}\,\P_{\Sigma}\,$. \qquad (See (\ref{eq-P=ip})).
\item [P8)] Commutation:  
$\Bigl(i_{(\Sigma_1\sqcup\Sigma_2)}\otimes\id\Bigr)
\Bigl(\id\otimes p_{(\Sigma_2\sqcup\Sigma_3)}\Bigr)\,=\,
p_{\bigl((\Sigma_1\#\Sigma_2)\sqcup\Sigma_3\bigr)}
i_{\bigl(\Sigma_1\sqcup(\Sigma_2\#\Sigma_3)\bigr)}$\newline
\hphantom{xxxxxxxxxxxxxx} with each $\Sigma_j$ connected, etc.
\een

This list of properties is, in fact, also sufficient for the existence 
of a generalized TQFT. Furthermore, this generalization comes out to be
precisely the one defined and discussed in Section 2.3, and, conversely, 
implies the Axioms V1-V5.

\btm\lll{thm-conn-half}\

Suppose there is a functor, $\V^{conn}$, and maps $i_{\Sigma}$ and 
$p_{\Sigma}$, such that properties P1-P8 are full filled.

Then there exists a unique, half-projective TQFT, w.r.t. 
${\tt x}=\V^{conn}(S^1\times S^2)\,$ and $\mu_0$,
$$
\V\,:\;\Cb 0\,\longrightarrow\,{\tt R}-mod\;,
$$
such that $i_{\Sigma}=\V(\Pi_{\Sigma})\,$, 
$p_{\Sigma}=\V(\Pi_{\Sigma}^{\dagger})\,$,
and \V\ specializes to $\V^{conn}\,$ on $\Cc 0\,$.
\etm

{\em Proof:} The assignment of a map $\V(M)$ to a connected cobordism, $M$,
is uniquely determined and well define by Axiom V1 and Property P2. Using V2
and P5 its extension to a disconnected cobordisms, $M$, is found from  the 
decompositions in Lemma~\ref{lm-tens-cob}. From the discussion in 
Section 2.2 it follows that compatibility of \V\ with the symmetric 
tensor structure allows us to consider only elementary compositions as in 
(\ref{eq-elprod-cob}). Assume that the cobordisms in this formula are
$$
M_1\,:\;\Sigma_A\,\longrightarrow\,\Sigma_B\sqcup\Sigma_C\qquad\;{\rm and}\;
\qquad M_2\,:\;\Sigma_C\sqcup\Sigma_D\,\longrightarrow\,\Sigma_E
$$
so that $(\id\otimes M_2)(M_1\otimes\id):\Sigma_A\sqcup\Sigma_D\to
\Sigma_B\sqcup\Sigma_E\,$. The connectivity cocycle is given by 
$\mu_0=\mu_0(\breve M_2,\breve M_1)=\beta_0(\Sigma_C)-1\,$. We find 
$$
\bar{rcl}
{\tt x}^{\mu_0}\V(\breve M_2)\V(\breve M_1)&=&{\tt x}^{\beta_0(\Sigma_C)-1}
 \Bigl(\id\otimes\bigl\{p_{\Sigma_E}\V^{conn}(\hat M_2)
i_{\Sigma_C\sqcup\Sigma_D}\bigr\}\Bigr)\Bigl(\bigl\{p_{\Sigma_B\sqcup\Sigma_C}
\V^{conn}(\hat M_1)i_{\Sigma_A}\bigr\}\otimes \id\Bigr)\\
\ ^{(by\, P4)}&=&
{\tt x}^{\beta_0(\Sigma_C)-1}
 \Bigl(\id\otimes\bigl\{p_{\Sigma_E}\V^{conn}(\hat M_2)
i_{\Sigma_C^{\#}\sqcup\Sigma_D^{\#}}(i_{\Sigma_C}\otimes i_{\Sigma_D})\bigr\}
\Bigr)\\
&&\hfill\times\Bigl(\bigl\{(p_{\Sigma_B}\otimes p_{\Sigma_C})
p_{\Sigma_B^{\#}\sqcup\Sigma_C^{\#}}
\V^{conn}(\hat M_1)i_{\Sigma_A}\bigr\}\otimes \id\Bigr)\\
&=&
{\tt x}^{\beta_0(\Sigma_C)-1}(p_{\Sigma_B}\otimes p_{\Sigma_E})
 \Bigl(\id\otimes\bigl\{\V^{conn}(\hat M_2)
i_{\Sigma_C^{\#}\sqcup\Sigma_D^{\#}}\bigr\}
\Bigr)(\id\otimes\P_{\Sigma_C}\otimes\id)\\
&&\hfill\times\Bigl(\bigl\{
p_{\Sigma_B^{\#}\sqcup\Sigma_C^{\#}}
\V^{conn}(\hat M_1)\bigr\}\otimes \id\Bigr)
(i_{\Sigma_A}\otimes i_{\Sigma_D})\\
\ ^{by\,(P7)}&=&
(p_{\Sigma_B}\otimes p_{\Sigma_E})
 \Bigl(\id\otimes\bigl\{\V^{conn}(\hat M_2)
i_{\Sigma_C^{\#}\sqcup\Sigma_D^{\#}}(\V^{conn}(\Lambda_{\Sigma_C})\otimes\id)
\bigr\}\Bigr)(\id\otimes\P_{\Sigma_C}\otimes\id)\\
&&\hfill\times\Bigl(\bigl\{
p_{\Sigma_B^{\#}\sqcup\Sigma_C^{\#}}
\V^{conn}(\hat M_1)\bigr\}\otimes \id\Bigr)
(i_{\Sigma_A}\otimes i_{\Sigma_D})\\
\ ^{by\,(P6)}&=&
(p_{\Sigma_B}\otimes p_{\Sigma_E})
 \Bigl(\id\otimes\V^{conn}(\hat M_2\circ\Lambda_{\Sigma_C}^{\%})
\Bigr)\\
&&\hfill\times\bigl(\id\otimes i_{\Sigma_C^{\#}\sqcup\Sigma_D^{\#}}\bigr)
\bigl(p_{\Sigma_B^{\#}\sqcup\Sigma_C^{\#}}\otimes\id\bigr)
\bigl(\V^{conn}(\hat M_1)\otimes \id\bigr)
(i_{\Sigma_A}\otimes i_{\Sigma_D})\\
\ ^{by\,(P8)}&=&
(p_{\Sigma_B}\otimes p_{\Sigma_E})
 \Bigl(\id\otimes\V^{conn}(\hat M_2\circ\Lambda_{\Sigma_C}^{\%})
\Bigr)p_{\Sigma^{\#}_B\sqcup(\Sigma^{\#}_C\#\Sigma^{\#}_D)}\\
&&\hfill\times\,\,
i_{(\Sigma^{\#}_B\#\Sigma^{\#}_C)\sqcup\Sigma^{\#}_D}
\bigl(\V^{conn}(\hat M_1)\otimes \id\bigr)
(i_{\Sigma_A}\otimes i_{\Sigma_D})\\
\ ^{by\,(P6)}&=&
(p_{\Sigma_B}\otimes p_{\Sigma_E})
p_{\Sigma^{\#}_B\sqcup\Sigma^{\#}_E}\V^{conn}(X_2)\V^{conn}(X_1)
i_{\Sigma^{\#}_A\sqcup\Sigma^{\#}_D}
(i_{\Sigma_A}\otimes i_{\Sigma_D})\\
\  ^{by\,(P4)}&=&p_{\Sigma_B\sqcup\Sigma_E}\V^{conn}(X_2\circ X_1)
i_{\Sigma_A\sqcup\Sigma_D}\\
\ ^{by\,(P3)}&=&\V\bigl(M\bigr)\hfill {\rm with}\qquad M=\Pi_{\Sigma_B\sqcup\Sigma_E}^{\dagger}\circ X_2
\circ X_1\circ
\Pi_{\Sigma_A\sqcup\Sigma_D}\;.
\ear
$$
In order to apply Naturality P6 in  this calculation we had to choose
 cobordisms, 
$\Lambda^{\%}_{\Sigma_C}:\Sigma^{\#}_C\#\Sigma^{\#}_D\to
\Sigma^{\#}_C\#\Sigma^{\#}_D\,$, 
$X_1:\Sigma^{\#}_A\#\Sigma^{\#}_D\to
\Sigma^{\#}_B\#\Sigma^{\#}_C\#\Sigma^{\#}_D\,$,
and
$X_2:\Sigma^{\#}_B\#\Sigma^{\#}_C\#\Sigma^{\#}_D\to
\Sigma^{\#}_B\#\Sigma^{\#}_E\,$, such that they full fill the 
following equations
$$
\bar{rclr}
\Pi_{\Sigma^{\#}_C\sqcup\Sigma^{\#}_D}\circ
\Bigl(\Lambda_{\Sigma_C}\sqcup\id_{\Sigma_D^{\#}}\Bigr)
&=&
\Lambda_{\Sigma_C}^{\%}\circ\Pi_{\Sigma^{\#}_C\sqcup\Sigma^{\#}_D}&\qquad\quad
(N1)\\
\Pi_{(\Sigma_B^{\#}\#\Sigma_C^{\#})\sqcup\Sigma^{\#}_D}\circ
\Bigl(\hat M_1\sqcup\id_{\Sigma_D^{\#}}\Bigr) 
&=&
X_1\circ\Pi_{\Sigma_A^{\#}\sqcup\Sigma^{\#}_D}&(N2)\\
\Bigl(\id_{\Sigma_D^{\#}}\sqcup(\hat M_2\circ\Lambda_{\Sigma_C}^{\%})\Bigr)
\circ
\Pi_{\Sigma_B^{\#}\sqcup(\Sigma_C^{\#}\#\Sigma^{\#}_D)}^{\dagger}
&=&
\Pi_{\Sigma^{\#}_B\sqcup\Sigma^{\#}_E}^{\dagger}\circ X_2&(N3)
\ear
$$
In order to complete the proof that \V\ is a half-projective functor, we
still need to show that the above cobordisms $M$, with $\hat M=X_2\circ X_1$,
is in fact $\breve M_2\circ \breve M_1\,$. 
This is accomplished by basically the same
calculation, only now for cobordisms and in reverse order:
$$
\bar{rcl}
M\qquad\,&=&\Pi_{\Sigma_B\sqcup\Sigma_E}^{\dagger}\circ X_2
\circ X_1\circ
\Pi_{\Sigma_A\sqcup\Sigma_D}\\
\ ^{assoc.}&=&
\Bigl(\Pi_{\Sigma_B}^{\dagger}\sqcup\Pi_{\Sigma_E}^{\dagger}\Bigr)\circ
\Bigl(\Pi_{\Sigma_B^{\#}\sqcup\Sigma_E^{\#}}^{\dagger}\circ X_2\Bigr)\circ
\Bigl(X_1\circ\Pi_{\Sigma_A^{\#}\sqcup\Sigma_D^{\#}}\Bigr)\circ
\Bigl(\Pi_{\Sigma_A}\sqcup\Pi_{\Sigma_D}\Bigr)\\
\ ^{N2\,\&\,N3}&=&
\Bigl(\Pi_{\Sigma_B}^{\dagger}\sqcup
\bigl(\Pi_{\Sigma_E}^{\dagger}\circ\hat M_2\circ\Lambda_{\Sigma_C}^{\%}\bigr)\Bigr)
\circ\Pi_{\Sigma_B^{\#}\sqcup(\Sigma_C^{\#}\#\Sigma^{\#}_D)}^{\dagger}\circ
\Pi_{(\Sigma_B^{\#}\#\Sigma_C^{\#})\sqcup\Sigma^{\#}_D}\circ
\Bigl(\bigl(\hat M_1\circ\Pi_{\Sigma_A}\bigr)\sqcup\Pi_{\Sigma_D}\Bigr)\\
\ ^{Lemma~\ref{lm-comm-Pi}}&=&
\Bigl(\Pi_{\Sigma_B}^{\dagger}\sqcup
\bigl(\Pi_{\Sigma_E}^{\dagger}\circ\!\hat M_2\circ\!\Lambda_{\Sigma_C}^{\%}\bigr)\Bigr)
\!\circ\!
\Bigl(\id_{\Sigma_B^{\#}}\sqcup\Pi_{\Sigma_C^{\#}\sqcup\Sigma^{\#}_D}\Bigr)
\!\circ\!
\Bigl(\Pi_{\Sigma_B^{\#}\sqcup\Sigma_C^{\#}}^{\dagger}
\sqcup\id_{\Sigma^{\#}_D}\Bigr)\!\circ\!
\Bigl(\bigl(\hat M_1\circ\!\Pi_{\Sigma_A}\bigr)\sqcup\Pi_{\Sigma_D}\Bigr)\\
&=&
\Bigl(\Pi_{\Sigma_B}^{\dagger}\sqcup
\bigl(\Pi_{\Sigma_E}^{\dagger}\circ\hat M_2\circ\Lambda_{\Sigma_C}^{\%}\circ
\Pi_{\Sigma_C^{\#}\sqcup\Sigma^{\#}_D}\bigr)\Bigr)
\!\circ\!
\Bigl(\bigl(\Pi_{\Sigma_B^{\#}\sqcup\Sigma_C^{\#}}^{\dagger}\circ
\hat M_1\circ\Pi_{\Sigma_A}\bigr)\sqcup\Pi_{\Sigma_D}\Bigr)\\
\ ^{N1}&=&
\Bigl(\id_{\Sigma_B}\sqcup
\bigl\{\Pi_{\Sigma_E}^{\dagger}\circ\hat M_2\circ
\Pi_{\Sigma_C^{\#}\sqcup\Sigma^{\#}_D}\circ\bigl(\Lambda_{\Sigma_C}
\sqcup\Pi_{\Sigma_D}\bigr)\bigr\}\Bigr)\\
&&\hfill\circ
\Bigl(\bigl\{\bigl(\Pi^{\dagger}_{\Sigma_B}\sqcup\id_{\Sigma^{\#}_C}\bigr)\circ
\Pi_{\Sigma_B^{\#}\sqcup\Sigma_C^{\#}}^{\dagger}\circ
\hat M_1\circ\Pi_{\Sigma_A}\bigr\}\sqcup\id_{\Sigma_D}\Bigr)\quad\\
\ ^{Lemma~\ref{lm-comp-Pi}}&=&
\Bigl(\id_{\Sigma_B}\sqcup
\bigl\{\Pi_{\Sigma_E}^{\dagger}\circ\hat M_2\circ
\Pi_{\Sigma_C^{\#}\sqcup\Sigma^{\#}_D}\circ\bigl(\Pi_{\Sigma_C}
\sqcup\Pi_{\Sigma_D}\bigr)\bigr\}\Bigr)\\
&&\hfill\circ
\Bigl(\bigl\{\bigl(\Pi^{\dagger}_{\Sigma_B}\sqcup\Pi_{\Sigma_C}^{\dagger}
\bigr)\circ
\Pi_{\Sigma_B^{\#}\sqcup\Sigma_C^{\#}}^{\dagger}\circ
\hat M_1\circ\Pi_{\Sigma_A}\bigr\}\sqcup\id_{\Sigma_D}\Bigr)\quad\\
\ ^{assoc.}&=&
\Bigl(\id_{\Sigma_B}\sqcup
\bigl\{\Pi_{\Sigma_E}^{\dagger}\circ\hat M_2\circ
\Pi_{\Sigma_C\sqcup\Sigma_D}\bigr\}\Bigr)
\circ
\Bigl(\bigl\{\Pi_{\Sigma_B\sqcup\Sigma_C}^{\dagger}\circ
\hat M_1\circ\Pi_{\Sigma_A}\bigr\}\sqcup\id_{\Sigma_D}\Bigr)\\
&=&
\bigl(\id_{\Sigma_B}\sqcup M_2\bigr)
\circ
\bigl(M_1\sqcup\id_{\Sigma_D}\bigr)\hfill=\,\breve M_2\circ\breve M_1\,.\\
\ear
$$
This completes the proof of the theorem. \ep

In the remaining sections we use Theorem~\ref{thm-conn-half} to construct
half-projective TQFT's from known ones. We will verify the necessary and
sufficient properties P1-P8 for a very general class of examples.

\prg{4.3}{The Example of Extended TQFT's}

In this section we show that if the connected TQFT from (\ref{eq-conn-TQFT})
originates from an extended structure, most of the properties entering
Theorem~\ref{thm-conn-half} are already full filled. Thus for the 
remainder of this and the next section we shall require that there be a 
series of functors $\V^{conn}_n\,$
as in (\ref{eq-ext-TQFT}), with $\Cb n$ replaced by $\Cc n\,$.

Connected, extended TQFT's exist for a quite general class of abelian,
categories:

\btm [{[KL]}]\lll{thm-KL} 

Suppose $\cal C$ is an abelian, rigid, balanced, modular, braided
tensor category over a field ${\tt R}=k\,$, for which the coend 
$F=\int X^{\vee}\btimes X\,$ exists. 

Then there is  a series of functors, $\V^{conn}_n\,$, on the categories
\Cc n , as in 
(\ref{eq-ext-TQFT}), which respects both types of tensor products and
2-categorial compositions, and for which 
\beq\lll{eq-TQFT-vect}
\V^{conn}_0(\Sigma_g)\;=\;
Inv\bigl(\,\underbrace{F\btimes\ldots\btimes F }_{g\;\,{\rm times}}\,\bigr)
\qquad,
\eeq
where $\btimes$ is the (braided) tensor product of $\cal C$, and 
$\Sigma_g$ is the closed
surface of genus $g\,$.\newline ( We denote $Inv(X)\equiv Hom_{\cal C}(1,X)\,$).
\etm

As explained in [Ke3] this specializes to the Reshetikhin-Turaev [RT] invariant
(often identified with the Chern-Simons quantum field theory) if $\cal C$
is semisimple, and to the Hennings-invariant [H] if ${\cal C}={\cal A}-mod\,$.
The coend is in the first case $F=\bigoplus_jj^{\vee}\btimes j\,$, where 
$j$ runs over a representative set of simple objects, and in the second
case $F={\cal A}^*\,$, equipped with the coadjoint action $ad^*\,$.
\medskip

The part of the this generalization  that will be relevant for this section
is that the connected functor needed in Theorem~\ref{thm-conn-half}
descends from the functor $\V_1:\,\Cc 1\,\to{\cal C}\,$, and that the 
latter is a functor of braided tensor categories. Including also the 
intermediate tangle presentations this is made precise in the following
commutative diagram.
\beq\lll{diag-1to0}
\bar{cccccc}
\V_1^{conn}:\vspace*{.1cm}\qquad&\Cc 1&\TO{$\cong$}{50}&\ctinf 1)\cn&\TO{}{60}&{\cal C}\\
\vspace*{.2cm}&\quad\Biggl\downarrow\phi_0&&\qquad\quad\Biggl
\downarrow \mbox{\large/}\,{\tau}{\mbox{\footnotesize-move}}&&
\quad\;\Biggl\downarrow Inv\\
\V_0^{conn}:\qquad&\Cc 0&\TO{$\cong$}{60}&\ct 0)\cn &\TO{}{60}&{\tt R}-mod
\ear
\eeq

All of the vertical arrows are surjections. The functor $\phi_0$ is the 
filling
functor that was described in (\ref{eq-fill}) of Section 2.2. It assigns to 
a once punctured surface, $\Sigma\pc\,$, a corresponding closed surface,
$\Sigma\,$, by gluing in a disc, and to a cobordisms, $M\pc$, with corners
a morphism, $M\in \Cc 0\,$, by filling in a tube, $D^2\times I$.
The $\tau$-move is described in  Appendix A.3. 
As in [KL] we denote by $\ctinf 1)\cn$ the tangle category, which has 
the same generators as $\ct 0)\cn$, but which is not subject to the $\tau$-move.
It represents isomorphically $\Cc 1\,$, which generalizes the presentation
of the mapping class groups of punctured surfaces, $\Sigma\pc$, from [MP].
The $\tau$-move accounts for isotopies over the puncture so that $\ct 0)\cn$
actually represents $\Cc 0$. 
The functor $\V^{conn}_1$ is finally constructed by assigning to a 
tangle a system of
morphisms with naturality properties, and lifting those to a morphism between 
tensor powers of the coend, $F$, see [L1], [L2], [KL], and also [Ke3] for 
a less technical summary. Let us use the  shorthand 
$X_{\Sigma}:=\V^{conn}_1\bigl(\Sigma\pc\bigr)\in obj\bigl({\cal C}\bigr)\,$ 
for the 
object assigned by $\V_1^{conn}$ to a punctured surface $\Sigma\pc\,$. 
In order to describe the condition imposed by the $\tau$-move suppose that
the strand cossing all strands emerging from one boundary of the 
tangle-diagram carries a representation $Y\,$. The crossing itself is described
by the morphism 
$\varepsilon(X_{\Sigma},Y):\,X_{\Sigma}\btimes Y\to Y\btimes X_{\Sigma}$, where
$\varepsilon$ is the braid constraint of $\cal C$. $\tau$-invariance
implies that $\mu(X_{\Sigma},Y)=\varepsilon(Y,X_{\Sigma})
\varepsilon(X_{\Sigma},Y)$ is represented trivially. This can be
 done by passing 
to the invariance, since  for any $f\in Inv(X_{\Sigma})\,$ we have by naturality
\beq\lll{eq-inv-mu}
\mu(X_{\Sigma},Y)(f\btimes \id_{Y})= (f\btimes \id_{Y})\mu(1,Y)=
(f\btimes \id_{Y})\qquad.
\eeq
Thus we can construct  $\V_0^{conn}$ for  $\Sigma$ or $M$ by first choosing 
punctured representatives $\Sigma\pc$ or $M\pc$, apply to these $\V_1^{conn}$,
and then map the result into ${\tt R}-mod$ by the $Inv$-functor so that,
e.g., $\V_0^{conn}(\Sigma)=Inv(X_{\Sigma})\,$.  The assignment
$$
X_{\Sigma_g}\;=\;\underbrace{F\btimes\ldots\btimes F}_{g\;\,{\rm times}}
$$
then explains the formula for the vector spaces in Theorem~\ref{thm-KL}.
\medskip

The benefit of this description lies in the fact that the horizontal
arrows in the top row of (\ref{diag-1to0}) 
are functors of braided tensor categories. 
The tensor product of $\Cc 1$ is given by gluing the boundaries of 
two surfaces to a three-holed sphere, 
$\Sigma_{0,3}=S^2-\bigl(D^2\sqcup D^2\sqcup D^2\bigr)\,$, and accordingly the
two cylindrical boundary components  of two cobordisms with corners to two
respective pieces in $\Sigma_{0,3}\times I\,$. We shall make identifications 
with standard surfaces that are compatible with those for the ordered
connected sums in Section 4.1, using  isomorphisms 
\beq\lll{eq-tens-surf}
\Sigma_1\pc\btimes\Sigma_2\pc\;\cong\;\bigl(\Sigma_1\#\Sigma_2\bigr)\pc\qquad.
\eeq
The category $\ctinf 1)\cn$ also admits a natural tensor product, given by the 
(opposite) juxtaposition of two tangles. The presentation functor can be chosen, 
such that the tensor structure is strictly respected. Moreover, the construction 
of $\V_1^{conn}$ is such that it is also a strict tensor functor into $\cal C$, i.e.,
\beq\lll{eq-V1-tens}
X_{\Sigma_1\#\Sigma_2}\;=\;X_{\Sigma_1}\btimes X_{\Sigma_2}\quad.
\eeq
A connected cobordism between surfaces with several components,
$M:\Sigma_1\sqcup\ldots\sqcup\Sigma_K\to\tilde\Sigma_1\sqcup\ldots
\sqcup\tilde\Sigma_L\,$ may be similarly first  described by a morphism
$$
\hat M\pc\,:\;\bigl(\Sigma_1\#\ldots\#\Sigma_K\bigr)\pc\,\longrightarrow\,
\bigl(\tilde\Sigma_1\#\ldots\#\tilde\Sigma_L\bigr)\pc\qquad,
$$
which is presented by a tangle $\T(M)\in\ctinf 1)\cn$. 
The original cobordism $M$ is
then presented by the image of $\T(M)\,$ in the tangle category $\ct 0)\cn\,$,
 where we have 
introduced a $K+L$ additional $\tau$-moves, one for every group of strands 
representing a boundary component. 

In order to guarantee invariance under  the $\tau$-moves at the source ends
of the tangle, we may proceed analogously and restrict  
$\V_1^{conn}\bigl(\hat M\pc\bigr)\,$ to the tensor product of the invariances
of the objects associated to the individual groups. If we start by carrying
 out the
reduction to the cobordism $\hat M$ of closed, connected surfaces as in 
(\ref{diag-1to0}) we find a first candidate for the inclusion from
Theorem~\ref{thm-conn-half}. Specifically, we have that
$$
\V_0^{conn}(\hat M)i^0_{\Sigma}\qquad {\rm \ only\ depends\ on}\quad 
\hat M\circ\Pi_{\Sigma}\quad,
$$
where we use the canonical injection:
\beq\lll{eq-inj-0}
i^0_{\Sigma}\,:\;Inv(X_{\Sigma_1})\otimes\ldots\otimes Inv(X_{\Sigma_K})
\;\into{25}\;Inv\bigl(X_{\Sigma_1}\btimes\ldots\btimes X_{\Sigma_K}\bigr)\quad.
\eeq
The difficulty that remains, is to find a projection in reverse direction in the
case that the target surface is also disconnected, i.e., $L>1\,$. 
In general, if
the vector spaces are given by the invariances as above, a canonical map with
these desired properties does not exist. Still, we can define canonical matrix 
elements. More precisely, for every choice of 
invariances, $f_j\in Inv(X_{\Sigma_j})\,$,
and coinvariances, $g_j\in Cov(X_{\tilde\Sigma_j})\,$, (denoting 
$Cov(Y)\equiv Hom_{\cal C}(Y,1)\,$) we have that also
$$
(g_1\btimes\ldots\btimes g_L)\V^{conn}_1\bigl(\hat M\pc\bigr)
(f_1\btimes\ldots\btimes f_K)\qquad.
$$
only depends on 
$M=\Pi_{\tilde\Sigma}^{\dagger}\circ \phi_o(\hat M\pc)\circ\Pi_{\Sigma}\,$.
This circumstance naturally leads us to first construct $\V_0$ on the morphisms
spaces, and then reconstruct the vector spaces. Generally, let us define 
for a set of objects  $A_j,B_j\in obj({\cal C})\,$ the null space:
$$
\bar{l}
H^0\,:=\;\Bigl\{h:A_1\btimes\ldots\btimes A_K\to B_1\btimes\ldots\btimes B_L\;:
\hspace*{9cm}\\
\hfill(g_1\btimes\ldots\btimes g_L)h(f_1\btimes\ldots\btimes f_K)=0\quad {\rm for\ all\ }
f_j\in Inv(A_j),\,g_j\in Cov(B_j)\,\Bigr\}\;.
\ear
$$ 
From this we define the space of matrices:
$$
H\bigl(A_1,\dots,A_K\Bigl | B_1,\dots,B_L\bigr)\;:=\;
\raise .5 ex
\hbox{$Hom_{\cal C}\bigl(A_1\btimes\ldots\btimes A_K, 
B_1\btimes\ldots\btimes B_L\bigr)$}
\mbox{\large /}\raise -.5 ex\hbox{${H^0}$}\;.
$$
For a morphism, $I$, between the tensor products of the $A_j$'s and $B_j$'s, 
let us also denote its image in the above space (i.e., its class modulo $H_0$)
 by $[I]\,$. A natural definition of the TQFT-functor for disconnected surfaces
on only the morphism spaces is thus
\beq\lll{eq-quot-mor}
\V_0(M)\,:=\;[\V_1^{conn}\bigl(\hat M\pc\bigr)]\qquad\;\;\in\;
H( X_{\Sigma_1},\dots,X_{\Sigma_K}\Bigl | X_{\tilde\Sigma_1},\dots,
X_{\tilde\Sigma_L}\bigr)\quad.
\eeq
Even in the connected case $H(A\Bigl | B)\,$ is usually going to be smaller
than $Hom_{\cal C}(A,B)\,$, if $\cal C$ is not semisimple. This is due to the fact
that the canonical pairing
\beq\lll{eq-Inv-pair}
Cov(X)\otimes Inv(X)\;\longrightarrow\;{\tt R}
\eeq
is degenerate for most objects $X\,$. Let us denote the null spaces of this
pairing by $Cov^0(X)$ and $Inv^0(X)$, respectively. It is easily seen that,
e.g, $Inv^0(X)$ is mapped to  $Inv^0(Y)$ by $Inv(f)$ 
for a morphism $f:X\to Y$, and that $[f]=0\,$, if all of $Inv(X)$ is mapped into
$Inv^0(Y)\,$. Still, we can think of the morphisms in the 
$H(\_\Bigl |\_)$-spaces as maps between vector spaces, if we pass to the 
quotients
\beq\lll{eq-def-mod-Inv}
\overline{Inv}(X)\;:=\;\frac {\,Inv(X)\,} {Inv^0(X)}\qquad\;{\rm and}\;\qquad
\overline{Cov}(X)\;:=\;\frac {\,Cov(X)\,} {Cov^0(X)}\qquad.
\eeq
Assuming that  {\tt R} acts nicely on these spaces this assertion is made 
more precise in the following lemma:

\blm\lll{lm-H-on-vs}
Suppose that for an abelian tensor category, $\cal C\,$,
 over {\tt R} the spaces $Inv(X)\,$, 
$Inv^0(X)\,$, $\overline{Inv}(X)\,$, $Cov(X)$, etc., are free
{\tt R}-modules, and that the exact sequence 
\beq\lll{eq-inv-seq}
0\,\longrightarrow\,Inv^0(X)\,\into{25}\,Inv(X)\,\onto{25}\,
\overline{Inv}(X)\,\longrightarrow\,0\qquad,
\eeq
as well as the analogous one for $Cov(X)$, are split over {\tt R}.

\noindent
Then we have 
\ben
\item Duality: $\;\overline{Cov}(X)\;\cong\;
Hom_{\tt R}\bigl(\overline{Inv}(X),{\tt R}\bigr)\,$ and vice versa.
\item $H\bigl(A_1,\dots,A_K\Bigl | B_1,\dots,B_L\bigr)\;\cong\;
Hom_{\tt R}\Bigl(\overline{Inv}(A_1)\otimes\ldots\otimes
\overline{Inv}(A_K),\,\overline{Inv}(B_1)\otimes\ldots\otimes
\overline{Inv}(B_L)\Bigr)\;$.
\een
\elm
The proof is standard, and makes use of the fact that the split sequence
in (\ref{eq-inv-seq}) allows us to choose dual basis in 
$\overline{Inv}(X)$ and $\overline{Cov}(X)\,$. If {\tt R} is a field 
the prerequisites of Lemma~\ref{lm-H-on-vs} are of course always 
full filled (assuming finite dimensions). Yet, as usual let us consider a more
general situation, in order to 
indicate the extend to which our constructions are 
the only possible ones. 

The identity in the second part of Lemma~\ref{lm-H-on-vs} show 
that we have to modify the  vector spaces of the connected TQFT from
Theorem~\ref{thm-KL} or, more generally, the diagram in (\ref{diag-1to0}) 
by assigning
\beq\lll{eq-VS}
\V_0(\Sigma)\;:=\;\overline{Inv}(X_{\Sigma})\;,\qquad\;{\rm if\ }\Sigma
{\rm \ is\ connected}\,,
\eeq
and extending this to disconnected surfaces by tensor products.
This is then compatible with the assignment of linear maps to a morphisms,
$M$, where we first choose a $\hat M\pc\in \Cc 1\,$, compute its image in
the $H(\_\Bigl |\_)$-space as in (\ref{eq-quot-mor}), and then apply the 
above isomorphism into the corresponding space of linear maps. On $\Cc 0$
the functor $\V_o$ will thus be given by the factored version of 
$\V^{conn}_0$, where we divided out the null spaces as in 
(\ref{eq-def-mod-Inv})

Using the easily verified property that a canonical injection, as in 
(\ref{eq-inj-0}), maps, e.g., $Inv^0(X)\otimes Inv(Y)\into{17}
Inv^0(X\btimes Y)\,$, we can factorize $i^0_{\Sigma}$ from (\ref{eq-inj-0})
into a map
\beq\lll{eq-inj-fac}
i_{\Sigma}\,:\;\overline{Inv}(X_{\Sigma_1})\otimes\ldots\otimes 
\overline{Inv}(X_{\Sigma_K})
\;\into{25}\;\overline{Inv}
\bigl(X_{\Sigma_1}\btimes\ldots\btimes X_{\Sigma_K}\bigr)\quad.
\eeq
Since $\overline{Inv}$ and $\overline{Cov}$ are now dual spaces, we may
define canonical projections, $p_{\Sigma}\,$, associated to the 
connecting cobordisms $\Pi^{\dagger}_{\Sigma}\,$. 
They shall be the adjoints of the corresponding
inclusion of coinvariances, i.e.,
\beq\lll{eq-proj-fac}
{p_{\Sigma}}^*\,\,:\;\overline{Cov}(X_{\Sigma_1})\otimes\ldots\otimes 
\overline{Cov}(X_{\Sigma_K})
\;\into{25}\;\overline{Cov}
\bigl(X_{\Sigma_1}\btimes\ldots\btimes X_{\Sigma_K}\bigr)\quad.
\eeq
It is also straight forward to see that with the  definitions in 
(\ref{eq-quot-mor}),
$\,\V_0(M)=p_{\tilde\Sigma}\V_0(\hat M)i_{\Sigma}\,$, if 
$M:\Sigma\to\tilde\Sigma\,$, and $\hat M$ is a corresponding morphism 
between the connected sums, such that 
$M=\Pi_{\tilde\Sigma}^{\dagger}\circ\hat M\circ\Pi_{\Sigma}\,$.
It is also clear from the  construction that 
$p_{\Sigma}i_{\Sigma}=id\,$ so that we have now the ingredients 
entering Theorem~\ref{thm-conn-half}, which full fill properties
P1, P2, and P3.

Moreover, associativity P4 follows immediately from the associativity
of the canonical inclusions of invariances and coinvariances.
Compatibility with symmetry as in P5 can be inferred from the fact 
that $\varepsilon(X_{\Sigma_1},X_{\Sigma_2})\,$ acts by naturality as
the transposition on $Inv(X_{\Sigma_1})\otimes Inv(X_{\Sigma_2})\,$
and hence also on $\overline{Inv}(X_{\Sigma_1})\otimes 
\overline{Inv}(X_{\Sigma_2})\,$. Also, property P8 is evident,
if we consider it for matrix elements. Specifically, choose
$f\in Inv(X_1)\,$, $g\in Inv(X_1\btimes X_2)\,$, 
$\alpha\in Cov(X_1\btimes X_2)\,$, and $\beta\in Cov(X_3)\,$, and denote by
$\overline f$, $\overline g$, $\overline\alpha$, and $\overline\beta$, the 
images in the quotient spaces. We then have the obvious identities
$$
\Bigl\langle\overline\alpha\otimes \overline\beta,\,
\Bigl(i_{(1\sqcup 2)}\otimes\id\Bigr)
\Bigl(\id\otimes p_{(2\sqcup 3)}\Bigr)\overline f\otimes \overline g
\Bigl\rangle
\;=\;
(\alpha\btimes \beta)(f\btimes g)
\;=\;
\Bigl\langle\overline\alpha\otimes \overline\beta,\,
p_{\bigl((1\# 2)\sqcup 3\bigr)}
i_{\bigl(1\sqcup(2\# 3)\bigr)}
\overline f\otimes \overline g
\Bigl\rangle\quad.
$$
Finally, the construction of $\V_0$ from the tensor functor, $\V_1^{conn}\,$,
as in (\ref{diag-1to0}) allows us to infer naturality P6 from the following
relation between cobordisms:
\blm\lll{lm-nat-top}
Suppose $M_j\pc$, with $j=1,2$,  are morphisms in $\Cc 1\,$, $\Sigma$ 
($\,\tilde\Sigma$) is the 
union of the closed source (target) surfaces, and the choice of the 
$\Pi_{\Sigma}\,$'s is as in Section 4.1. Then
$$
\phi_0\Bigl(M_1\pc\otimes M_2\pc\Bigr)\circ\Pi_{\Sigma}\;=\;
\Pi_{\tilde\Sigma}\circ\Bigl(\phi_0(M_1\pc)\sqcup \phi_0(M_2\pc)\Bigr)\quad.
$$
\elm

{\em Proof:} The relation is readily verified, given the tangle presentation
in (\ref{fig-conn-mor}), the fact that the tensor product in $\ctinf 1)\cn$ 
is given
by juxtaposition, and that $\phi_0$ only introduces another relation, but may be chosen
as identity on representing tangles. 

The formula may also by understood directly, by considering both sides of the equation 
as $M_1\pc\sqcup M_2\pc\,$ to which certain elementary manifolds are attached along the
two cylindrical boundary pieces. On the left hand side we glue in 
$B=S^2-(D^2\sqcup D^2)\times [0,1]\,$ to get 
$\phi_0\bigl(M_1\pc\otimes M_2\pc\bigr)\,$, and attach a 2-cell, $C_2$, to the piece,
$\cong  S^2-(D^2\sqcup D^2)$, of $X$ in the source surface,
in order to realize the composition with 
$\Pi_{\Sigma}\,$. On the right hand side we first glue in two tubes, $D^2\times J$,
to get $\phi_0(M_1\pc)\sqcup \phi_0(M_2\pc)$ and then describe the composition with
$\Pi_{\tilde\Sigma}$ by attaching a 1-cell, $C_1$, at the discs, $D^2$, in the 
target surface that belong to the tubes. 
In both cases the combined glued in piece,
$
X\cup C_2\cong (D^2\times J)\cup C_1\cup (D^2\times J)\cong D^3\,
$,
is a ball, and  the cylindrical boundary pieces of 
$M_1\pc\sqcup M_2\pc\,$ are glued to $D^3$ along two annuli that are 
embedded in its boundary $\partial D^3=S^2\,$.\ep

Let us summarize the findings of this section in the next lemma:

\blm\lll{lm-sec-43} Suppose we have a connected TQFT functor $\V_0^{conn}\,$
that descends from a functor $\V_0^{conn}\,$ of braided tensor categories as in
(\ref{diag-1to0}) (e.g., the one proposed in Theorem~\ref{thm-KL}). Then we
can construct a map $\V_0:\Cb 0\to {\tt R}\,$, which satisfies the properties
P1-P6, and P8 from Section 4.2. The vector spaces are the quotient spaces as in
(\ref{eq-VS}). 
\elm

The only thing left to investigate, in order to complete the
construction of a half-extended TQFT, is thus the projectivity property P7.
This will be done in the next section, as it 
relates to more specific properties of the constructions starting from an abelian, braided
tensor category $\cal C$ over a field $k\,$. The triviality of {\tt x} will turn out to 
determine completely the semisimplicity of $\cal C$.
\medskip

Let us conclude this section with an example of how vector spaces change, when we divide out the
null spaces as in (\ref{eq-def-mod-Inv}) for a non-semisimple category. The result in the
case of ${\cal C}=U_q(s\ell_2)-mod\,$, with $q$ a primitive $2m+1$-st root of unity, is
described, e.g.,  in [Ke3]. The dimensions of the vector spaces for a torus are given
by 
$$
dim\Bigl(Inv(F)\Bigr)\,=\,3m+1\qquad\quad{\rm and}\qquad\quad 
dim\Bigl(\overline{Inv}(F)\Bigr)=\,2m\,\;.
$$
It is a quite remarkable fact that at least for prime $2m+1$, the space $Inv^0(F)\,$ is not only 
an invariant subspace  of the  representation of the mapping class group derived from 
$\V_0^{conn}\,$, but that it is also a direct summand. I.e.,  the sequence in (\ref{eq-inv-seq})
is also split as a sequence of $SL(2,{\bf Z})$-modules.

Notice also that $\overline{Inv}(F)$
is naturally identified with an invariance of the semisimple trace sub-quotient, 
$\overline{\cal C}^{tr}\,$, which is the starting point for the TQFT extending the
Reshetikhin-Turaev invariant. The vector space of the torus of the latter is however 
only $m$-dimensional since the coend of $\overline{\cal C}^{tr}\,$ is smaller than the
image of the coend of $\cal C$.

\prg{4.4}{Integrals, Semisimplicity, and ${\tt x}=\V(S^1\times S^2)\,$}

\noindent
\parbox{4.5in}{\quad In view of the tangle presentation in (\ref{fig-lambda-mor}) 
it is clear that the 
key to understanding the representation of $\Lambda_{\Sigma}$ in a TQFT, and hence the
projectivity P7, is to explain the effect of a trivially framed annulus 
around a vertical
strand as in diagram (\ref{fig-an0}) in algebraic terms. 

Geometrically, the surgery along the annulus {\sf A} is equivalent to connecting
an $S^1\times S^2$ to the manifold we surger on. The strand {\sf S} passing through
{\sf A} then indicates a path that generates $\pi_1(S^1\times S^2)\,$.}
\parbox{2in}
{\beq\lll{fig-an0}
\begin{picture}(40,100)

\put(18,82){{\sf S}}

\put(20,15){\line(0,1){15}}
\put(20,34){\line(0,1){40}}

\put(17,45){\oval(40,26)[l]}
\put(23,45){\oval(40,26)[r]}
\put(17,32){\line(1,0){6}}

\put(40,60){{\sf A}}
\end{picture}
\eeq}

It also has a natural interpretation in the language of cobordism categories, if we 
consider the tangle in (\ref{fig-an0}) as a morphism in $\ctinf 2)$, i.e., the 
tangle category with one external strand, {\sf S}, and no $\tau$-move. 
Here the diagram from (\ref{fig-an0})
represents a cobordism, $\lambda\ppc:\,\Sigma_{0,2}\to\Sigma_{0,2}\,$, in
$\Cc 2\,$. It is explicitly given by the cylinder $\Sigma_{0,2}\times[0,1]\,$,
inside of which we have performed a surgery along the meridian of one 
of the cylindrical  boundary pieces.  It is clear that if we apply a filling functor,
which glues a tube to one of these boundary pieces, the result in $\Cc 1$ will be
$\phi_1(\lambda\ppc)=\bigl(D^2\times [0,1]\bigr)\#(S^1\times S^2)\,$.
Moreover,
we have $\lambda\ppc\circ\lambda\ppc\,=\,\lambda\ppc\#(S^1\times S^2)\,$. 
(here $\#$ is always the sum with the interior).

At this point it turns out to be  rather  instructive to include 
into our discussion the 2-categorial
picture  of cobordisms and TQFT's, that was outlined in Section 2.4. 
 For example, we can think of 
a surface with a hole as a 1+1-cobordism 
$\Sigma\pc:\emptyset\to S^1\,$, and of a
cobordism $M\pc\in \Cb 1\,$ as a 2-morphism between two such 1-morphisms. 
Also,  $\Sigma_{0,2}:S^1\to S^1$   may be seen as  the identity 1-morphism 
$\id_{S^1}\,$ on the circle, and $\lambda\ppc:\id_{S^1}\,\Rightarrow
\,\id_{S^1}\,$ 
is a 2-endomorphism.

A 2-category also implies a composition operation, $\bullet_1$, which is  the 
usual composition on the 1+1-cobordisms, and which is naturally 
extended to the  cobording 2-morphisms. Since $\id_{S^1}$ can be composed with
any $\Sigma\pc$ we can form the $\bullet_1$-composite of any
$M\pc\in\Cb 1\,$ with $\lambda\ppc$. It is clear that the result can
be obtained from $M\pc$ also by doing a surgery along a meridian that is 
pushed off the special cylindrical piece of $\partial M\pc\,$.
From this and the rules of the tangle presentation $\Cc 1\isto \ctinf 1)\cn\,$
it follows easily that if $\T(M\pc)$ presents $M\pc$ the tangle for the 
composite is given by placing a trivially framed annulus around the 
entire tangle  $\T(M\pc)\,$. 

If we introduce also $\lambda\pc_{\Sigma}:=\,
\lambda\ppc\bullet_1\id_{\Sigma}\pc\,$, we easily verify the following
identities from the 2-categorial distributive law:
\beq\lll{eq-lm-m}
\lambda\ppc\bullet_1 M\pc\;=\;\lambda\pc_{\Sigma_t}\circ M\pc\;=\; 
M\pc\circ\lambda\pc_{\Sigma_s}\qquad.
\eeq 
It is thus both  natural and useful to think of $\lambda\ppc$ as a natural
transformation on $\Cc 1\,$. In particular,
 the  associated morphism $\lambda\pc_{\Sigma_g}\,$, for the connected 
surface $\Sigma_g$  of genus $g$, is presented in $\ctinf 1)$ by
$2g$ vertical strands with a trivially framed  annulus around them. 
Comparing this to  (\ref{fig-lambda-mor}) and using the braided 
tensor structures of $\Cc 1$ and $\ctinf 1)\,$, we find
\beq\lll{eq-Lambda-tens}
\Lambda_{\Sigma}\pc\;=\;\lambda_{\Sigma_1}\pc\btimes\ldots\btimes
\lambda_{\Sigma_{K-1}}\pc\btimes \id_{\Sigma_K}\qquad,
\eeq
where the $\Sigma_j$'s are the connected components of $\Sigma\,$.
\medskip

From the definition of am extended TQFT as in (\ref{eq-2cat-TQFT}) in 
Section 2.4 - as well as what we expect from the property expressed
 in  (\ref{eq-lm-m}) - it follows that $\lambda\ppc\,$ is represented 
by a natural transformation 
of the identity functor of $\cal C$ to itself. In particular we have
the following:

\beq\lll{eq-lambda-natrsf}
\lambda(X_{\Sigma})\;:=\;\V^{conn}_1\bigl(\lambda_{\Sigma}\pc\bigr)
\eeq
(Recall that a natural transformation 
$\lambda\in Nat\bigl(id_{\cal C},id_{\cal C}\bigr)\,$ consists of an endomorphisms
$\lambda(X)\in End_{\cal C}(X)\,$ for every object, such that $f\lambda(X)=\lambda(Y)f$,
if $f:X\to Y\,$.) 
\medskip

Let us determine a few constraints on the transformation 
$X\mapsto\lambda(X)\,$. To begin with, note that the braid morphisms,
$\varepsilon(\Sigma_1\pc,\Sigma_2\pc):
\Sigma_1\pc\btimes\Sigma_2\pc\to\Sigma_2\pc\btimes\Sigma_1\pc$ of
$\Cc 1$, 
 can be presented in $\ctinf 1)\cn\,$ 
by the diagram in  (\ref{fig-braid-Pi}), except
that we combine the top ends into one group. Moreover, 
$\lambda_{\Sigma_1}\pc\btimes\id_{\Sigma_2}\,$ is given by diagram
(\ref{fig-lambda-mor}), with $K=2$. In the composite of these two tangles
we can slide the $2g_2$ lower strands one by one over the annulus and hence 
turn the overcrossing into an undercrossing. I.e., in $\Cc 1$ 
we have the following identity:
\beq\lll{eq-braid-cc1}
\varepsilon(\Sigma_1\pc,\Sigma_2\pc)
(\lambda_{\Sigma_1}\pc\btimes\id_{\Sigma_2})\;=\;
\varepsilon(\Sigma_2\pc,\Sigma_1\pc)^{-1}
(\lambda_{\Sigma_1}\pc\btimes\id_{\Sigma_2})
\quad.
\eeq
 
The corresponding conditions on the transformation of $\cal C$ is given by 
\beq\lll{eq-lambdaC-nat}
\mu(X,Y)(\lambda(X)\btimes\id_Y)\;=\;\lambda(X)\btimes\id_Y
\qquad
{\rm and}
\qquad
\mu(X,Y)(\id_X\btimes\lambda(Y))\;=\;\id_X\btimes\lambda(Y)\qquad
\eeq
where $\mu(X,Y)$ is the square of $\varepsilon$, as in (\ref{eq-inv-mu}).
To be precise, we should   impose this relation only for the special objects,
$X_{\Sigma}$. However, they will be sufficiently big so that this implies
the general statement by naturality.
 
Other topological considerations lead us to impose conjugation invariance
$\lambda(X)^t\;=\;\lambda(X^{\vee})\,$.
These two properties also correspond to generating, elementary moves in 
a ``Bridged Link
Calculus'', which replace the 2-handle slides of the conventional 
Kirby calculus, see [Ke2].
\medskip

Equation (\ref{eq-inv-mu}) also shows
that (\ref{eq-lambdaC-nat}) can be derived from a stronger condition, namely
that each $\lambda(X)$ has a decomposition into a monomorphism and an 
epimorphism, going through a multiple of the unit object:
\beq\lll{eq-moni-epi-lambda}
\lambda(X)\,:\;\,X\,\ONTO{$p_X^{\lambda}$}{40}\,\,
\underbrace{1\oplus\ldots\oplus 1}_{n_X\;\,{\rm times}}
\,\,\INTO{$i_X^{\lambda}$}{40}\,X\qquad.
\eeq
(see, e.g.,  [M] for existence and uniqueness of monic-epic-decompositions). 
Instead of the injection and projection in (\ref{eq-moni-epi-lambda}) we may 
also consider the injections $i^{\nu}_X\in Inv(X)\,$ and projections 
$p^{\nu}_X\in Cov(X)\,$, with $\nu=1,\dots, n_X\,$, for the individual
summands $\cong 1\,$ so that 
$\lambda(X)=\sum_{\nu=1}^{n_X} i^{\nu}_Xp^{\nu}_X\,$. Its action on 
invariance is determined by naturality, i.e., we have 
$\lambda(X)f=\lambda(1)f\,$ for $f:1\to X$, (using 
$\lambda(1)\in k\cong End_{\cal C}(1)\,$). We thus find the following linear
dependence between the $i^{\nu}_X$ and $f$:
\beq\lll{eq-inv-basis}
\sum_{\nu=1}^{n_X} \langle p^{\nu}_X, f\rangle\, i^{\nu}_X\;=\;
\lambda(1)\,f\qquad\qquad{\rm for\ all\ }\, f\,\in\,Inv(X)\quad.
\eeq
 
This formula allows us to establish a  relation between the number on the 
right side and   semisimplicity of the category:
\blm\lll{lm-x-Css}
Suppose $\cal C$ is an abelian, rigid, balanced tensor category over a field 
$k={\tt R}\,$, with finite dimensional morphism sets. Assume further that there
is a transformation $\lambda\in Nat\bigl(id_{\cal C},id_{\cal C}\bigr)\,$,
with $\lambda\neq 0$,  such that (\ref{eq-moni-epi-lambda}) holds for all 
objects $X$. Then 
$$
{\cal C}\quad{\it \ is\ semisimple,\ \ if\ and\ only\ if\ }\quad \lambda(1)\neq 0\qquad.
$$
\elm

{\em Proof:} In a balanced category we can construct traces, 
$tr_X:\,End_{\cal C}(X)\to k\,$, that are generally cyclic, and respect the tensor product.
As in [Ke4], we may then  define a category to be  semisimple, iff all 
pairings of the form
$$
Hom_{\cal C}(Y,X)\otimes Hom_{\cal C}(X,Y)\,\TO{\scriptsize composition}{80}\,
End_{\cal C}(X)\,\TO{$tr_X$}{40}\,\,k\;,
\qquad
$$
are non-degenerate. Rigidity allows us to reduce this to non-degeneracy of
the pairings of  invariance and coinvariance, as in  (\ref{eq-quot-mor}),
for all objects. Now, if $\lambda(1)=0\,$, (\ref{eq-inv-basis}) implies that 
$p^{\nu}_Xf=0$ for all $\nu$, $X$, and $f\in Inv(X)\,$, and hence a degeneracy 
of the pairing, if $\lambda(X)\neq 0$. Since we  assumed $\lambda\neq 0\,$,
 this proves  one implication.

If $\lambda(1)\neq 0$ and ${\tt R}=k$ is a field,
it follows from (\ref{eq-inv-basis}) that $\Bigl\{i_X^{\nu}\Bigr\}$ and 
$\Bigl\{\frac 1{\lambda(1)}p_X^{\nu}\Bigr\}$ are dual basis so that
we have non-degeneracy.\ep

For a  semisimple category it not hard to see from the proof that such a
transformations always exists, and that it is (up to a total scaling)
uniquely given by the projection $\P_X\in End(X)$ onto the maximal, trivial
sub-object. More precisely, we have
\beq\lll{eq-L=P}
\lambda(X)\;=\;\lambda(1)\P_X\quad.
\eeq

Next, we shall evaluate $\V_0(\Lambda_{\Sigma})\,$ for a surface $\Sigma$ with $K$
connected components, $\Sigma_j$. 
(We shall often use of the abbreviation $X_j\equiv X_{\Sigma_j}$).

Since $\V_1^{conn}$ is a functor of tensor categories, we can find from this 
an expression for  $\Lambda\pc_{\Sigma}$ using (\ref{eq-Lambda-tens}):
\beq\lll{eq-VL}
\V_1^{conn}(\Lambda_{\Sigma}\pc)\,=\,\lambda(X_{\Sigma_1})\btimes\ldots
\btimes\lambda(X_{\Sigma_{K-1}})\btimes\id_{X_{\Sigma_K}}\qquad.
\eeq
\medskip
If we use the factorization (\ref{eq-moni-epi-lambda}), we can write the action of this
morphism on an element, $f\in Inv(X_1\btimes\ldots \btimes X_K)\,$, as follows:
\beq\lll{eq-Vf-exp}
\V_1^{conn}(\Lambda_{\Sigma}\pc)f\,=\;\sum_{\{\nu_j\}}\,
i^{\nu_1}_{X_1}\otimes\ldots\otimes i^{\nu_{K-1}}_{X_{K-1}}\otimes 
\xi^{(\nu_1,\ldots,\nu_{K-1})}\quad,\hspace*{2cm}\;
\eeq
$$\hspace*{4cm}{\rm \ \ where\ }\qquad\quad\,\xi^{(\nu_1,\ldots,\nu_{K-1})}=
\Bigl(p^{\nu_1}_{X_1}\otimes\ldots\otimes p^{\nu_{K-1}}_{X_{K-1}}\otimes\id_{X_K}\Bigr)f\;\qquad\in
\,Inv(X_K)\quad.
$$
Here we used $X_k\cong 1\btimes\ldots\btimes 1\btimes X_K\,$. 
In the semisimple case 
we find  from (\ref{eq-Vf-exp}) that $\V_1^{conn}(\Lambda_{\Sigma})f=
(\id\btimes\ldots\btimes\P_{X_k})\V_1^{conn}(\Lambda_{\Sigma})f\,$ so that with 
(\ref{eq-L=P}) we have
$$
\V_1^{conn}(\Lambda_{\Sigma})f\;=\;\lambda(1)^{K-1}\,
\Bigl(\P_{X_1}\btimes\ldots\btimes\P_{X_k}\Bigr)f\quad.
$$
Semisimplicity also implies  $\overline{Inv}=Inv\,$ and 
$\overline{Cov}=Cov\,$ so that  we can use the bases consisting of
the $i^{\nu}_X$ and $p^{\nu}_X\,$, in order to express the 
injections and projections in (\ref{eq-inj-fac}) and  (\ref{eq-proj-fac}).
It follows immediately that the projection $\P_{\Sigma}$ in (\ref{eq-P=ip}) 
is precisely
given by the above tensor product of projections, restricted to invariance.
Hence,
\beq\lll{eq-L=lP-ss}
\V_0(\Lambda_{\Sigma})\;=\;\lambda(1)^{K-1}\P_{\Sigma}\qquad.
\eeq
In the case $\lambda(1)=0$ it follows from (\ref{eq-inv-basis})
 that the vectors
$i^{\nu}_X$ and $p^{\nu}_X$ all lie in the null spaces
$Inv^0(X)$ and $Cov^0(X)$, 
respectively. If $K>1$, this and (\ref{eq-Vf-exp}) imply that 
$$
\V_1^{conn}(\Lambda_{\Sigma})f\in Inv^0(X_1)\otimes\ldots \otimes 
Inv^0(X_{K-1})\otimes Inv(X_K)\;\subset\;
Inv^0\bigl(X_1\btimes\ldots\btimes X_K\bigr)\quad.
$$
It follows that  $\V_0(\Lambda_{\Sigma})=0$, i.e., (\ref{eq-L=lP-ss}) also
holds in the non-semisimple case. (For the case $K=1$ we have 
$\Lambda_{\Sigma}=\id_{\Sigma}\,$, and thus 
$\V_0(\Lambda_{\Sigma})=\id_{X_{\Sigma}}$ for either case.)
\medskip

In order to prove that the property P7 follows from the assignment of a natural
transformation with a decomposition as in (\ref{eq-moni-epi-lambda}),
 we still have 
to make the identification
\beq\lll{eq-x=lambda1}
\lambda(1)\;=\;{\tt x}\;\equiv\;\V_0(S^1\times S^2)\quad.
\eeq
To this end we  shall make the  assumption that the unit object
in $\cal C$ is irreducible, and that $\V_1^{conn}$ also preserves unit objects,
i.e., $X_{S^2}=1\,$. Equation (\ref{eq-x=lambda1}) follows now from 
$\lambda_{S^2}\pc=(\id_{S^2}\pc)\#(S^1\times S^2)\,$.
\medskip

For semisimple categories existence and uniqueness of transformations as in 
(\ref{eq-x=lambda1}) is obvious. Still, we wish to understand this 
assertion also in the general case, and make sure that the construction, e.g,  
in Theorem~\ref{thm-KL} actually assigns such a transformation to $\lambda\ppc\,$.
For this purpose it will be both useful and instructive to attribute another
algebraic interpretation to $\lambda\ppc$. 

This has its origins in the case ${\cal C}={\cal A}-mod$, where $\cal A$
is a finite dimensional Hopf algebra. Here a natural transformation, 
$\lambda\,$,
of the identity functor  is uniquely identified with a central element,
$\lambda\in{\tt Z}({\cal A})\,$, of the Hopf algebra. 
As the trivial representation
of $\cal A$ is given by its counit $\epsilon$, the relation 
(\ref{eq-moni-epi-lambda}) translates to 
\beq\lll{eq-Hopf-int}
y\lambda\;=\;\lambda y\;=\;\epsilon(y)\lambda\qquad {\rm \ for\ all}\quad 
y\in {\cal A}\;,\qquad\qquad\quad{\rm and\ } \;
\lambda(1)=\epsilon(\lambda)\qquad.
\eeq
Elements $\lambda$ satisfying this relation, which we  will call 
(two-sided) {\em cointegrals},  are well known in the theory
of Hopf algebras. (In  a more common convention
$\lambda$ is actually 
the cointegral of ${\cal A}^*$, which is the algebra used in the 
categorial description). Existence and uniqueness of integrals and cointegrals 
has been proven for finite dimensional Hopf algebras in [Sw]. In [LS] it is
also shown that $\cal A\,$ is semisimple, if and only if 
$\epsilon(\lambda)\neq 0\,$. Thus Lemma~\ref{lm-x-Css} is just the categorial
generalization of this result. 

The action of the associated natural transformation $X\mapsto\lambda(X)\,$
of the identity on ${\cal A}-mod\,$ to itself, is given by the canonical  
application of $\lambda$ to the $\cal A$-module, $V_X$. We have 
$\,n_X=dim\Bigl(im(\lambda(X))\Bigr)\,$. For example, if ${\cal A}=U_q(s\ell_2)'\,$,
with $q$ a root of unity, it follows from its representation theory, see [Rp],
that $n_X$ is given by the number of times  $P_0$ occurs as a direct summand 
in $V_X$,  where $P_0$ is {\em the} indecomposable, projective representation 
that  contains an invariant vector. In particular, $n_X=0$ if $V_X$ is fully
reducible. 
\medskip

The notion of cointegrals can be generalized to the type of tensor categories, 
$\cal C\,$, considered in Theorem~\ref{thm-KL}, for which the coend 
$\,F=\int X^{\vee}\otimes X\,$ exists.

In [L1] it is shown that $F$ has the structure of a categorial, braided
Hopf algebra in $\cal C$. Moreover, there is a natural Hopf algebra
pairing
$$
\omega\,:\;F\btimes F\,\longrightarrow\,1\qquad,
$$
for which the left and right null space $ker(\omega)\into{15} F$ coincide.
Let us call $\omega$ {\em balanced}, if $ker(\omega)$ is preserved by the 
balancing $v\in Nat(\id_{\cal C})$, where the action of $v$ is given by the 
lifting of $\id\btimes v(X)\in End(X^{\vee}\btimes X)\,$. This property has 
been introduced in [L1] as axiom (M2). We shall also say that $\cal C$ is
{\em strictly modular} if $ker(\omega)=0$.

The algebra of natural transformations of the 
identity of $\cal C$ is canonically isomorphic to the elements in $Cov(F)\,$.
Hence a cointegral of $F$ is defined to be a morphism,

\beq\lll{eq-CoInt}
\lambda\,:\;F\,\to\, 1\;,\qquad{\rm such\ that}\qquad
{\bar{cccc}
\vspace*{.1cm}&F&\TO{$\Delta$}{30}&F\btimes F\;\;\\
\vspace*{.1cm}&\lambda\Biggl\downarrow\;\,&&\qquad\;\Biggl\downarrow \id_F\btimes\lambda\\
&1&\TO{$e$}{30}&F\;\;
\ear}
\qquad {\rm commutes,}
\eeq
where $\Delta$ is the coproduct and $e$ is the unit of $F$.
( The previously discussed case of an ordinary Hopf algebra
is implied by the identity $Cov(F)=Cov\Bigl({\cal A}^*,ad^*\Bigr)\cong
{\tt Z}\bigl({\cal A}\bigr)\,$).

We also have the dual notion of a categorial integral for $F$, which
is given by an invariance,

\beq\lll{eq-Int}
\mu\,:\;1\,\to\,F\;,\qquad{\rm such\ that}\qquad
{\bar{cccc}
\vspace*{.1cm}&F&\TO{$ \mu\btimes\id_F$}{40}&F\btimes F\;\;\\
\vspace*{.1cm}&1^*\Biggl\downarrow\;\;\,&&\;\;\Biggl\downarrow m\\
&1&\TO{$\mu$}{40}&F\;\;
\ear}
\qquad {\rm commutes,}
\eeq
with $m$ being the multiplication, and $1^*$ the counit of $F$.

The existence and uniqueness, proven for ordinary Hopf algebras in  [Sw],
 is generalized to the categorial versions
in the following theorem:

\btm [{[L1]}]\lll{thm-exuni-L}
Suppose $\cal C$ is a category as above, which contains the coend $F$.

Then there exist (up to scalings) a unique integral, $\mu:I_{\mu}\to F$,
and a unique cointegral, $\lambda:F\to~I_{\lambda}\,$, where $I_{\mu}$
and $I_{\lambda}$ are invertible.

If the pairing, $\omega$, is balanced, then 
$I_{\mu}\cong I_{\lambda}\cong 1\,$,
and $\lambda$ and $\mu$ are as in (\ref{eq-CoInt}) and (\ref{eq-Int}),
respectively.
\etm

In the construction in [L2] and [KL] the transformation that is 
a-priori associated to the diagram in (\ref{fig-an0}) is given as
a coinvariance of $F$ by the composite
\beq\lll{eq-lmo}
\lambda^{+}\,:\;F\,\TO{$\id_F\btimes\mu$}{50}\,F\btimes F\,
\TO{$\omega$}{30}\,1\qquad.
\eeq

The following lemma compels us to add non-degeneracy of $\omega$ to
our list of requirements on $\cal C$.

\blm[{[Ke3]}]\

 Suppose $\cal C$ is as above, and $\omega$ is balanced.

Then $\lambda=\lambda^{+}\,$, if and only if $\cal C$ is strictly modular.
\elm

This now establishes the fact that the construction used in 
Theorem~\ref{thm-KL}
does in fact associate the cointegral of $F$ to the diagram in 
(\ref{fig-an0}).

Still, we need to show that the image of $\lambda(X)$ is actually
of the form $\,1\oplus\ldots\oplus 1\,$, in order to complete our
construction. For  categories  over fields, $k$,  with $char(k)=0$,
this is implied by a result of
P. Deligne, which asserts that  an 
object, on which the coaction of $F$ is trivial,
necessarily has to be the direct sum of units. It is found [D2] as a
 corollary to the existence of tensor products of abelian categories.
The same  result has been proven independently by V. Lyubashenko [L3]
by the use of {\em squared} coalgebras:

\blm[{[L3],[D2]}]\lll{lm-int=1+1}\

The cointegral of a braided tensor category over a field $k$, with $char(k)=0$,
and with coend $F$ 
factors for each object as in (\ref{eq-moni-epi-lambda}).
\elm

It thus follows that the construction in Theorem~\ref{thm-KL} also implies 
Property P7.
\medskip

\prg{4.5}{ Main Result, and Hints to Further Generalizations and Applications}
 
In this concluding section we shall summarize the various possible deviation from
our axioms, that we pointed out throughout this chapter, in order to find 
more general types of non-semisimple TQFT's for disconnected surfaces than
 half-projective TQFT's. We also discuss possible applications of the
half-projective formulation to ``classical limits'', where the {\tt x}
has the role of a renormalization parameter with ${\tt x}\to\infty\,$.
To begin with, let us state the main result of this chapter, which follows 
from results of the preceding three sections, and which completes the 
construction in [KL]:

\btm\lll{thm-main}

Suppose that $\cal C$ is a abelian, rigid, balanced, strictly modular,
braided tensor category, which is defined  over a field $k$
with $char(k)=0\,$, and
which contains the  coend $F$.

Then there exists a half-projective TQFT-functor, 
w.r.t. $\mu_0$ and ${\tt x}\in k\,$,
$$
\hspace*{2cm}
\V_0\,:\;\Cb 0\,\longrightarrow\,Vect(k)
\qquad{\rm with}\qquad
\V_0(\Sigma_g)=
\overline{Inv}\Bigl(\underbrace{\,F\btimes\ldots\btimes F\,}_{g\;{\rm times}}
\Bigr)\;,$$
such that
$$
{\tt x}\neq 0\qquad\quad{\it if\ and\ only\ if}\qquad\quad
{\cal C}\quad{\it is\ semisimple}\qquad.
$$

Moreover, for ${\tt x}=0$ the functor $\V_0$ restricted to $\Cc 0$
is the null space 
quotient of $\V_0^{conn}\,$, and for ${\tt x}\neq 0$ it is the 
TQFT extending the Reshetikhin-Turaev invariant.
\etm

The arguments given in Section 4.2 that lead up to P2, and the ones given in 
Section 4.3, yielding Lemma~\ref{lm-H-on-vs}, essentially necessitate the 
division of the original invariances by the null spaces $Inv^0(X_{\Sigma})\,$,
in order to achieve compatibility of the TQFT with the tensor structures in
$\Cb 0$ and ${\tt R}-mod$, as required in the axioms  V1-V5. Nevertheless,
by Theorem~\ref{thm-KL} these spaces do carry quite interesting representations
of mapping class groups, or other features of cobordisms in $\Cc 0$. 

In Section 2.3 we alluded to the possibility of circumventing the alternative
of Corollary~\ref{cor-altern} by admitting non-canonical symmetry structures on
${\tt R}-mod$ so that ${\tt x}\in{\tt R}\,$ may be, e.g., nilpotent. This 
would yield a richer filtration of $Inv(X)$ with respect to the action of
$\lambda(X)$ on $X$. 
\medskip

A more promising approach is to relax the Axioms V1-V5. Recall for example 
that V5 has already been put in by hand for mostly technical reasons. 
In order to explain what we might anticipate as a generalization of the tensor 
product rule, let us show that, in some cases,
 we can think of the $\overline{Inv}$-spaces as cohomologies.
In the non-semisimple case (${\tt R}=k$ a field) we have that $\lambda$
acts like a differential, i.e., $\lambda(X)^2=0\,$. Note also that
$im\Bigl(\lambda(X)\Bigr)\,\subset\,Inv^0(X)\,$ for all $X$, and that
$Inv(X)\,\subset\,ker\Bigl(\lambda(X)\Bigr)\,$. For 
${\cal C}\,=\,U_q(s\ell_2)-mod\,$, and if $X$ is the sum of only projective and
irreducible representations, it follows from the representation theory of 
$U_q(s\ell_2)$, see [Rp], that the first inclusion for the image is in fact an 
isomorphism. Moreover, if we restrict the action of $\lambda(X)$ to
the trivial weight space of $X$, 
then the second inclusion for the kernel is an isomorphism
for $X=1$ and for  $X=P_0$ 
(i.e., the unique indecomposable, projective representation
that contains an invariant vector). 

It has been worked out in [Os] that the adjoint representation does in fact
contain only projective and irreducible summands. Thus we conclude for
the vector space of the torus
$$
\V(S^1\times S^1)\,\equiv\,\overline{Inv(F)}\;=\;H^*\Bigl(F_{00},\lambda(F)\Bigr)\;,
$$
where the coend $F$ is as usual 
given by the functions on $U_q(s\ell_2)$ with coadjoint
representation, and $F_{00}$ is the intersection of the summand with trivial
Casimir value and the trivial weight space.

This formula for the vector space of a torus suggests to consider besides 
the tensor product also a derived functor,
analogous to {\bf Tor}, in order to retain some information of 
the $Inv^0$-spaces.
\bigskip

The formalism of half-projective TQFT's may also have applications in the 
case,  where $\cal C$ is semisimple so that ${\tt x}\in k\,$ is invertible.
For a fixed, finite ${\tt x}$ we recover an ordinary TQFT by rescaling the 
canonically constructed TQFT functor by
 $\V(M)\longrightarrow{\tt x}^{-\be_0(M)}\V(M)\,$. 
Yet, for ``classical limits'',
in which the renormalization parameter ${\tt x}$ will tend to infinity,
the anomaly may give  us some estimates on the divergence of the canonically
defined  $\V(M)\,$. Unlike the non-semisimple case not only triviality but
the exact value of $\varrho(M)$ should be of interest.

In the construction from [KL] the integral and cointegrals admit a 
canonical normalization:
$$
\mu\,=\,{\cal D}^{-1}\sum_{j\in{\cal J}}\,d(j)\widetilde{coev}_j\qquad\quad
{\rm and}\qquad\quad \lambda\,=\,{\cal D}ev_1\quad,
$$
where $ev$ and $\widetilde{coev}$ are the evaluation and flipped 
coevaluation, $d(j)$ are the quantum-dimensions for simple objects 
$j\in{\cal J}\,$, and
$$
{\tt x}\,=\, {\cal D}\,:=\,\pm\sqrt{\sum_{j\in{\cal J}}\,d(j)^2\,}\quad.
$$
The notation for $\cal D$ is the same as in [T]. The normalization is 
determined by $\lambda\cdot\mu=1$, and relation (\ref{eq-lmo}), which are
imposed by invariance under the local, interior moves of cancellation 
(or $\nhopf$-move) and modification, see [Ke3].
\medskip

For the Chern-Simons quantum field theory 
with a simple, connected, and simply connected  
gauge group, $G$, the set $\cal J$ is identified with the (heighest)
weights in an elementary, truncated alcove of the weight space of $G$, 
whose size depends on the
level of the theory (see, e.g., [KW]). In the classical limit the level
goes to infinity so that eventually every dominant weight will belong to 
$\cal J$. Hence $|{\cal J}|\to\infty\,$, and we have
$$
{\tt x}=\V(S^1\times S^2)\,\longrightarrow\,\infty\;.
$$
Evidently, this limit of TQFT's is ill defined on most cobordisms.

Still, there are cases, in which the limit exists in a given sense. 
For example, if $G=SU(2)$, and ${\cal J}\,=\,\{1,\ldots,N\}\,$, where
the labels are given by the dimension of the irreducible 
$SU(2)$-representation, we can consider the (projective) representation of 
$SL(2,{\bf Z})=\pi_0\Bigl({\cal D}i\!f\!f(S^1\times S^1)^+\Bigr)\,$ on 
$\V(S^1\times S^1)={\bf C}^N\,$, induced  by the TQFT. If we identify
the labels with points in ${\bf R}^+\,$, via $x_j\,=\,\frac j{\sqrt N}\,$,
then we can define the limits of the generators of $SL(2,{\bf Z})$ as 
unitary operators. In particular, {\sf S}
 is identified with the
Fourier transform on 
$\V_{\infty}\bigl(S^1\times S^1\bigr)=L^2({\bf R}^+,dx)\,$, 
and the limit of {\sf T}
 is given by the 
multiplication operator of $e^{ix^2}\,$ on the same space, where
$$
{\sf S}\,:=\,{\left[\bar{cc}0&1\\-1&0\ear\right]}\qquad\quad{\rm and}
\qquad\quad
{\sf T}\,:=\,\scriptstyle{\left[\bar{cc}1&0\\1&1\ear\right]}\qquad.
$$
For a connected cobordism, $M$, with maximal, free, interior group
of rank $\varphi(M)$, we find analogous to Lemma~\ref{lm-ord-rho}
from Theorem~\ref{thm-spc} and Lemma~\ref{lm-decomp-ogr} that
\beq\lll{=}
\V_0(M)\;=\;{\tt x}^{\varphi(M)}\,\V_0(M_2)\V_0(M_1)\qquad.
\eeq
Here, $M_2$ and $M_1$ are connected cobordisms, for which we have
$\varphi(M_j)=0$. Note that the latter also holds true for the 
invertible cobordisms, for which the classical limit existed.

For general $M$ we may ask by which order $\|\V_0(M)\|\,$ diverges in
the  limit with ${\cal J},\, {\tt x}\,\to\,\infty\,$. Assuming that 
the composition of $\V_0(M_1)$ and $\V_0(M_2)$ does not degenerate we expect 
to find from (\ref{=}) that $\varrho(M)$ is at least a lower bound on this
order. 

It is still quite crude, as we can see in the example of $S^1\times\Sigma\,$,
where $\varrho(M)$ is roughly half of the true order. More precisely, if we
use the Verlinde formula, $dim(V_{\Sigma_g})\,=\,{\cal D}^{2g-2}\sum_j
d(j)^{2-2g}\,$, (see, e.g., [T]), and that for $SU(2)\,$ we have
${\tt x}={\cal D}\sim N^{\frac 3 2}\,$, we compute from 
Lemma~\ref{eq-torus=dim}:
$$
\V_0(S^1\times\Sigma_g)\qquad\sim{\tt x}\;{\rm\ if\ }g=0\,,\qquad\;
\sim{\tt x}^{\frac 5 3}\;{\rm\ if\ }g=1\,,\;{\rm and}\qquad\;
\sim{\tt x}^{2g-1}\,\,{\rm\ if}\,g>1\quad.
$$
In contrast to that we have that for $g>0$ the number 
$\varrho(S^1\times\Sigma_g)\,$ is given precisely by the maximal number of 
non-separating curves on $\Sigma_g$, which is exactly $g$. Also,
$\varrho(S^1\times S^2)=1$ so that 
$\varrho(S^1\times\Sigma_g)\,=\,max(g,1)\,$ and the decompositions  in 
Lemma~\ref{lm-torus-ord} are  in fact maximal.
\bigskip

\ms

\setcounter{chapter}{0}
\subsection*{Appendix}\lll{pg-A}
\prg{A.1}{Proofs of Section 3.4}
\ 
\newline

\ub{\em A.1.1 Proof of Lemma~\ref{lm-gen-id}\ :\ } Generically we may 
assume that $f:\,M\to\gamma\,$ is a differentiable function, which is transversal
to the embedded graph $J(\gamma)\,$. This implies that there are no 
critical points in a vicinity of the graph, and that the edges of $J(\gamma)\,$
 are 
tangential to the level sets of $f$ only at a finite number of interior
points. Let us add the latter as vertices of valency two to $\gamma\,$.

By genericity we may also assume that for the vertices $v_j\in\gamma\,$,
the images $h(v_j)\in\gamma\,$ (where $h=f\circ J\,$) are distinct points 
in the interiors of the edges. We can choose disjoint open intervals $I_j^*$
and $\overline{I_j}\subset I_j^*\,$, such that 
 $v_j\in I_j\,$, and $I^*_j\,$ lies inside of  
an edge of $\gamma$. Hence we can find a function $\psi:\gamma\to\gamma\,$,
which is the identity outside of all of the $I^*_j\,$, which is strictly 
monotonous outside of $I_j$, and which collapses $I_j$ to the vertex $v_j\,$,
i.e., $\psi(I_j)=\{v_j\}\,$. The function $f^{\sim}:=\psi\circ f\,$ has the
property that it collapses an entire neighborhood of a vertex of the embedded
graph $J(\gamma)\,$ to a point. Since we can choose $f^{\sim}\,$ arbitrarily 
close to $f\,$ by making the intervals $I^*_j$ smaller and smaller,
we shall assume this property, by genericity,  already for $f$.

Let us thus consider the compact regions 
$B_j^o=f^{-1}\bigl(h(v_j)\bigr)\,$ in $M$, for which the 
$J(v_j)\in B_j^o$ are interior points. At points outside of the $B_j^o$ and in
a vicinity of the edges $f$ is regular, and the level surfaces
are transversal to the edges.  Hence we may assume that for each
edge $e(j,k)\subset\gamma$, joining the vertex $v_j$ to the 
vertex $v_k$, 
there is an embedding,
$\rho_{e(j,k)}:\,D^2\times[0,1]\,\into{20}\,M\,$, such that 
$t\mapsto \rho_{e(j,k)}(0,t)\,$ parametrizes $J(e(j,k))\subset M\,$,
and $f$ is constant on the disc-fibers, i.e., 
$f\bigl(\rho_{e(j,k)}(p,t)\bigr)\,=\,f\bigl(\rho_{e(j,k)}(0,t)\bigr)\,$.
Moreover, we may choose the parametrization, such that 
$\rho^{-1}_{e(j,k)}\bigl(B^o_j\bigr)\,=\,D^2\times [0,\varepsilon]\,$,
and $D^2$ is a disc of radius $\varepsilon\,$.

Next, we choose vicinities of the vertices, given by embeddings,
$\delta_j:\,\bigl(D^3,0\bigr)\,\into{20}\,\bigl(B_j^o,v_j\bigr)\,$, such that
the discs $\rho_{e(j,k)}\bigl(D^2\times\{\varepsilon\}\bigr)\,$ are 
(disjointly) contained in $\delta_j(S^2)\,$ for every edge joining $v_j\,$,
and $\rho_{e(j,k)}\bigl(D^2\times[0,\varepsilon]\bigr)\subset im(\delta_j)\,$.

The union of the images of the $\delta_j$ and the $\rho_{e(j,k)}$ forms a 
neighborhood $U(\gamma)\,$ of $J(\gamma)\,$ in $M$.
It may be given as a disjoint union of regions 
$$
R_{e(j,k)}\;:=\;\{\rho_{e(j,k)}(p,t):\,|p|<t<1-|p|,\;|p|\leq\varepsilon\,\}\;,
$$
and their complement, $\sqcup_jY_j\,$. Here, the $Y_j$ are the images of the
$\delta_j$ with the cones 
$R^{\vee}_{e(j,k)}\,:=\,\{\rho_{e(j,k)}(p,t):\,|p|<t<\varepsilon,\,
|p|\leq\varepsilon\}\,$ removed, and are therefore homeomorphic to a 
cone over a sphere, $S_j\cong S^2-\sqcup_eD^2\,$, which has a hole for every 
edge, $e$,  at $v_j\,$.
We may choose a parametrization, 
$\hat\delta_j:\,\bigl(S_j\times[0,\varepsilon],S_j\times\{0\}\bigr)\,\to\,
\bigr(CS_j,*\bigr)\cong\bigl(Y_j,*\bigr)\,\subset M\,$, such that the
second parameter is equal to $t=|p|\,$ at the common boundary with 
$R_{e(j,k)}\,$.

We define a function $\overline{\kappa_a}:\,\gamma\to\gamma\,$, which 
maps an edge, parametrized by $[0,1]\to\gamma:\,t\mapsto s(t)\,$, to itself,
such that
$$
\overline{\kappa_a}\bigl(s(t)\bigr)\,=\, s\Bigl(\frac {t-a}{1-2a}\Bigr)
\qquad\quad {\rm if}\quad a<t<1-a\;,
$$
$\overline{\kappa_a}\bigl(s([0,a])\bigr)=s(0)\,$, and 
$\overline{\kappa_a}\bigl(s([1-a,1])\bigr)=s(1)\,$.

Also define the projection
$$
\Pi_{e(j,k)}:\;R_{e(j,k)}\,\longrightarrow\,\gamma\,\,:\qquad\,
\rho_{e(j,k)}(p,t)\,\mapsto\,J^{-1}\bigl(\rho_{e(j,k)}(0,t)\bigr)\;,
$$
and from this the composite
$$
\kappa\bigl(\rho_{e(j,k)}(p,t)\bigr)\;:=\;\overline{\kappa_{|p|}}
\circ\Pi_{e(j,k)}\bigl(\rho_{e(j,k)}(p,t)\bigr)\;=\;
J^{-1}\biggl(\rho_{e(j,k)}\Bigl(0,\frac {t-|p|}{1-2|p|}\Bigr)\biggr)\;.
$$
Moreover, if  we define $\kappa\bigl(im(\hat\delta_j)\bigr)\,=\,v_j\,$, 
we obtain a continuous map
$$
\kappa\,:\;U(\gamma)\,\longrightarrow\,M\;.
$$
Now, $h:\,\gamma\to\gamma\,$ collapses, in the same way as 
$\overline{\kappa_{\varepsilon}}\,$, the $\varepsilon$-neighborhoods of a vertex
onto a point. Hence we can write
$$
h\;=\;h^{\$}\circ \overline{\kappa_{\varepsilon}}\quad,
$$
for some $h^{\$}:\,\gamma\to\gamma\,$, which is again homotopic to one.

Let $H:\,[0,\varepsilon]\times\gamma\to\gamma\,$ be such a homotopy, i.e.,
$H_{\varepsilon}\,=\,h^{\$}\,$ and $H_0=id\,$. We then define $f^{\$}:\,M\to
\gamma\,$ by 
$$
\bar{rclr}
f^{\$}(x)\;&=&\;f(x)\qquad\qquad& {\rm for \ }\, x\in M-U(\gamma)\\
f^{\$}\bigl(\rho_{e(j,k)}(p,t)\bigr)&=&H_{|p|}\Bigl(\kappa
\bigl(\rho_{e(j,k)}(p,t)\bigr)\Bigr)&{\rm on\ }\, R_{e(j,k)}\\
f^{\$}\bigl(\hat\delta_j(s,r)\bigr)&=&H_r(v_j)&{\rm on\ }\, Y_j\;\quad.
\ear
$$

For the $|p|=\varepsilon\,$-piece of $R_{e(j,k)}$ this gives
$$
\bar{ccccc}
H_{\varepsilon}\Bigl(\kappa\bigl(\rho_{e(j,k)}(p,t)\bigr)\Bigr) &=& 
h^{\$}\circ\overline{\kappa_{\varepsilon}}\circ
\Pi_{e(j,k)}\bigl(\rho_{e(j,k)}(p,t)\bigr)&=&h\Bigl(J^{-1}
\bigl(\rho_{e(j,k)}(p,t)\bigr)\Bigr)\\
&=&f\bigl(\rho_{e(j,k)}(0,t)\bigr)&=&f\bigl(\rho_{e(j,k)}(p,t)\bigr)\;,
\ear
$$
since $f$ is constant along the disc-fibers. Also, 
$H_{|p|}\Bigl(\kappa\bigl(\rho_{e(j,k)}(p,t)\bigr)\Bigr)= H_r\bigl(v_j\bigr)\,$,
if $|p|=t=r\,$, and 
$H_{\varepsilon}(v_j)\,=\,h^{\$}(v_j)\,=\,h(v_j)\,=\,f(Y_j)\,$,
so that $f^{\$}$ is continuous. Finally, on $J(\gamma)\,$ we have
$$
f^{\$}\bigl(\rho_{e(j,k)}(0,t)\bigr)\;=\;
H_0\Bigl(\kappa\bigl(\rho_{e(j,k)}(0,t)\bigr)\Bigr)\;=\;
J^{-1}\bigl(\rho_{e(j,k)}(0,t)\bigr)\;
$$
so that $f^{\$}\circ J=id\,$.
\ep

\ms

\ub{\em A.1.2 Proof of Part 2) of Lemma~\ref{lm-decomp-ogr}\ :\ } 

By Lemma~\ref{lm-gen-id} we may assume that $f\circ J\,$ is the identity.
With $\,M^o_n\,=\,f^{-1}(\gamma_n)\,$, for $n=1,2$,
 we thus obtain, as in the proof of the 
first part, cobordisms with $M=M^o_2\circ M_1^o\,$, and $\gamma_n$ embeds into
$M^o_n\,$, such that 
$\,\gamma_n\,\into{20}\,M^o_n\,\longrightarrow \,\gamma_n\,$ is the identity.
As a first step let us alter $f$, such that the $f^{-1}(\gamma_n)$ are 
connected:

Denote by $M^!\,$ the component of $M^o_1$, which contains $\gamma_1\,$, and 
by $B^{\nu}_1\,$ the other components. Each of these is a cobordism
$\,B^{\nu}_1:\,\emptyset\longrightarrow\Sigma\,$, where 
$\,\Sigma\neq\emptyset$, and $\Sigma\subset f^{-1}(P)\,$, i.e., it is a union
of the connected components $\Sigma_j^p\,$ from the proof of the first part. 
For given $B^{\nu}_1\,$, denote by $Q_{\nu}$ the set of labels $(p,j)$ that
occur in this union. Choose
for every $\nu$ a coordinate graph 
$g^{\nu}:\,B_1^{\nu}\longrightarrow\gamma_B^{\nu}\,$ 
(e.g., as in Part 1 of Lemma~\ref{lm-coord-ex}),
and extend this to
$$
\overline{g}^{\nu}\,:\;\overline{B_1}^{\nu}\,=\,B_1^{\nu}\,\sqcup
\coprod_{(p,j)\in Q_{\nu}}\Sigma_j^p\times [0,\varepsilon]\,\longrightarrow\,
\overline{\gamma_B}^{\nu}\;.
$$
Here we have added the collars of the boundary components, $\Sigma_j^p\,$,
of $B_1^{\nu}\,$ in $M_2\,$, and elongated the ends of $\gamma^{\nu}_B\,$
by intervals $[0,\varepsilon]\,$ to obtain $\overline{\gamma_B}^{\nu}\,$.
Each $\,\Sigma_j^p\times\{\varepsilon\}\,$, with $(p,j)\in Q_{\nu}\,$, is
mapped to an interior point, $p^{\#}\,$, in $\gamma_2$.
Since $\gamma_2$ is connected, there are maps 
$h^{\nu}:\,\overline{\gamma_B}^{\nu}\longrightarrow\,int(\gamma_2)\,$, such
that $(p,j)$ is assigned to the respective $p^{\#}\,$. We can thus define
a continuous function $\,\tilde f:\,M\to\gamma\,$, which is given by 
$h^{\nu}\circ\overline{g}^{\nu}\,$ on each $\overline{B_1}^{\nu}\,$, and
by $f$ on the rest of $M\,$. We now have that 
$\tilde{M}_1:= {\tilde f}^{-1}(\gamma_1)\,$ 
is connected. In the same way we can 
define from this coordinate map $\hat f:\,M\to\gamma\,$, for 
which also  ${\hat f}^{-1}(\gamma_2)\,$ is connected, and 
${\hat f}^{-1}(\gamma_1)\,$ is the composite of 
$\tilde{M}_1\,$ and the extra components of ${\tilde f}^{-1}(\gamma_2)\,$ - 
hence also connected.

We can thus assume that the $M^*_j=f^{-1}(\gamma_j)\,$ are connected, yet 
$f^{-1}(p)\,$, with $p\in P=\gamma_1\cap\gamma_2\,$, may still contain 
several components. As before we have to add an additional edge to
$\gamma\,$ for each components disjoint from the embedded $\gamma\,$. 
For such a component $\Sigma$ choose a path $\,t\mapsto q(t)\,$ in $M^*_1$,
which has  starting point $q(0)\in\Sigma\,$, and which ends in a point in 
$J(\gamma_1)\,$. For the corresponding path $r(t):=f(q(t))\,$ in $\gamma_1$,
we define the compact subsets $M_t\,:=\,f^{-1}\bigl(r([0,t])\bigr)\,$, 
with $M_t\subset M_s\subset f^{-1}(\gamma_1)\,$ for $t<s\,$. We set
$$
T\;:=\;\inf\bigl\{\,t:\,J(\gamma_1)\,\,{\rm and }\,\, \Sigma 
{\rm \ are \ connected \ in \ } M_t\,\bigr\}\;.
$$
Since $\, M_T=\bigcap_{t>T}M_t\,$ is the intersection of compact sets, in 
which $J(\gamma_1)\,$ and $\Sigma\,$ are connected, they are also connected 
in $M_T\,$.

Let $M_t^{\Sigma}\,$ be the component of $\Sigma\,$ in $M_t\,$ for $t<T\,$,
 and $M_t^{\gamma}\,$ its complement in $M_t\,$. On 
$\check M_T\,:=\,\bigcup_{t<T}M_t\,=:\,\check M^{\Sigma}_T
\cup \check M^{\gamma}_T\,$ we define a function $f^?\,$ as follows.
It shall be constant on the union $\check M^{\Sigma}_T\,$ of the 
$M^{\Sigma}_t\,$, with value $r(T)=f\bigl(q(T)\bigr)\,$, and it is $f$
on $\check M^{\gamma}_T\,$. Continuity of $r$ and $r(T)\not\in r([0,T[)\,$
imply that $\check M_T^{\Sigma}\,$ and $\check M_T^{\gamma}\,$ are disjoint 
open sets in the induced topology of $M_T\,$ so that $f^?\,$ is continuous.
It is also clear that $f^?\Bigl |\check M^{\Sigma}_T\,$ and 
$\,f^?\Bigl |\check M^{\gamma}_T\,$ have continuous extensions to
their closures in $M_T$, assigning the value $r(T)$ to the additional boundary
points in both cases. In particular, they can be glued to a continuous
function on $M_T$, and, setting $f^?=f$ on the closure of the complement 
$\overline{M^c_T}\,$, to a continuous function on $M^*_1\,$.

Define $\gamma'_1$ by attaching an interval, $[0,\varepsilon]\,$, to 
$\gamma_1$, such that $\{ 0\}\,$ is an endpoint of $\gamma_1'\,$ 
and $\{\varepsilon\}\,$
is identified with $r(T)=\gamma_1\,$. For a collar, 
$\,\Sigma\times[0,\varepsilon]\,\into{20}\,f^{-1}(\gamma_1)\,$, we can  then 
define a modified, continuous function, $f^{??}:\,M_1^*\to\gamma'_1\,$,
by setting $f^{??}=f^?\,$ on the complement of the collar, and on 
$\Sigma\times[0,\varepsilon]\,$ we define $f^{??}\,$ as the projection onto 
the additional edge, $[0,\varepsilon]\,$, of $\gamma_1'\,$. The path
$\,q:\,[0,T[\,\to\,\check M^{\Sigma}_T\,$ is thus connected to $J(r(T))\,$ in 
$M^*_1$ by a path $q':\,[T,T']\to M_1^*-\check M^{\Sigma}_T\,$. Clearly,
$r'(t):= f^{??}\bigl(q'(t))\,$ is a closed path in 
$\gamma_1\subset\gamma_1'\,$. Suppose $r''$ is the inverse path (parametrized
by $[T',T'']\,$), in $\gamma_1\,$, and set 
$\,q''=J\circ r'':\,[T',T'']\to M^*_1\,$. The path-composite 
$\hat q:= q''* q'* q:\,[0,T'']\to M^*_1\,$ starts at $\Sigma$, and ends in $r(t)\,$, and, moreover, $f^{??}\circ\hat q\,$ is homotopic to 
$[0,\varepsilon]\into{20}\gamma'_1\,$. Thus if we rescale the parametrization of $\hat q$ from $[0,T'']$ to $[0,\varepsilon]\,$, we have an embedding
$J':\gamma\into{20}\,M^*_1\,$, extending $J\,$, such that $f^{??}\circ J'\,$
is homotopic to the identity.

We apply this process to every additional surface, until we arrive at 
a faithful coordinate graph, $\gamma\,$, of $M^*\,$.
\ep

\prg{A.2}{The Spaces $\H(M,{\bf G})$, the Numbers $\be_j(M)\,$, and 
Further Anomalies}

In (\ref{eq-def-effbetti}) of Section 2.2 we introduced the interior Betti 
number $\be_j(M)$ of a cobordism, $M$. The coboundary $\mu_0=-\delta\be_0\,$
turned out to be non-negative and entered the definition of half-projective
TQFT's. In this appendix we shall be interested in the properties of the 
coboundary $\mu_1=-\delta\be_1\,$, for which we find 
$\mu_{\partial}:=\mu_1-\mu_0\geq 0\,$.  If we consider homotopy instead of
homology the analog for $\be_1$
is $\varphi(M)$, which has been introduced in Section 3.2. There we also
defined in (\ref{eq-mu-pi}) the coboundary 
$\mu_{\pi}\,$, which is in  analogy
 to $\mu_{\partial}\,$. Since the explicit computation of 
$\mu_{\pi}$ is quite hard, we shall give here instead a derivation of 
a formula for $\mu_{\partial}\,$.

We shall first identify the  numbers $\be_1(M)\in {\bf Z}^{+,0}\,$ 
with  the dimensions of the  spaces
\beq\lll{eq-def-HM}
\H\bigl(M,{\bf G}\bigr)\;:=\;coker\Bigl(H_1\bigl(\psi,{\bf G})\Bigr)\qquad,
\eeq
where we  usually consider coefficients ${\bf G}\,=\,{\bf Z}$ or $\bf Q$, and
$\psi$ is the inclusion of the punctured boundary pieces, 
$-\Sigma_s\sqcup\Sigma_t\,$, into $M\,$.
In the case, where we have no punctures, i.e., 
$M\,=\,\phi_0(M)\,$,
we know that $H_1(\psi)$
is half-rank, and the dimension is in fact given by $\be_1\,$. Hence
\beq\lll{eq-HZ-dim}
\H\bigl(M,{\bf Z}\bigr)\;=\;{\bf Z}^{\be_1(M)}\,\oplus\,\,Torsion\qquad,
\eeq
which,  (by naturality of universal coefficients)  implies
\beq\lll{eq-HQ-dim}
\H\bigl(M,{\bf Q}\bigr)\;=\;\H\bigl(M,{\bf Z}\bigr)\otimes{\bf Q}\;
=\;{\bf Q}^{\be_1(M)}\qquad.
\eeq
For punctured surfaces these formulae are still correct, since we have in 
analogy to Lemma~\ref{lm-pi0} the following:
\blm\lll{lm-H1iso-fill}
\
Suppose $M\in\Cb *$ is a cobordism, $\phi_0(M)\in\Cb 0$ is as in 
(\ref{eq-fill}), and $i_{\phi}\,:\,M\,\hookrightarrow\,\phi_0(M)\,$
is the respective inclusion.
\ben
\item The inclusion induces an isomorphism (for both types of coefficients):
$$
\H(i_{\phi})\,:\;\H\bigl(M\bigr)\,\isto\,\H\bigl(\phi_0(M)\bigr)\qquad.
$$
\item The first interior Betti-numbers are computed directly from
$$
\be_1(M)\;=\;\beta_1(M)\,-\,\frac 1 2 \beta_1(-\Sigma_s\sqcup\Sigma_t)\,+\,
dim\bigl(ker(i_M)\bigr)\qquad,
$$
where $i_M\,:\,H_1\bigl(\sqcup^N S^1\times I\bigr)\,\to\,H_1(M)\,$ is induced
by the inclusion of the cylindrical boundary pieces into the cobordism.
\een
\elm

{\em Proof: } It is enough to prove the identities for gluing a tube 
$D^2\times I$ into
only one cylinder $S^1\times I$ of a cobordism $M$. If 
$\Sigma=-\Sigma_s\sqcup\Sigma_t\,$, the corresponding surface 
$\Sigma'$ of the filled cobordism is obtained  by 
gluing discs into the holes of $\Sigma$. Since, the tube and the discs 
are contractible, the inclusions of boundary pieces induce the 
following transformation of Mayer-Vietoris sequences:
\beq\lll{eq-MV-fill}
\bar{cccccccccc}
\vspace*{.2cm}&0&\to&H_1(S^1\sqcup S^1)&\INTO{$i_{\Sigma}$}{30}
&H_1(\Sigma)&\ONTO{$p_{\Sigma}$}{30}&
H_1(\Sigma')&\to&0\\
\vspace*{.2cm}&&&p_1\Biggl\downarrow&&\psi_1\Biggl\downarrow&&\psi_1'\Biggl\downarrow&&\\
&&&H_1(S^1\times I)&\TO{$i_M$}{30}& H_1(M)&\ONTO {$p_M$}{30}&
H_1(M')&\to&0
\ear
\eeq
In the top row $p_{\Sigma}$ is onto because the discs were glued to different
components of $\Sigma$ (see definition of \Cb *), and $i_{\Sigma}$ is 
always into. As $p_1$ is onto, the two sequences factor into the isomorphism in
Part 2 of the lemma. 

The analog of (\ref{eq-MV-fill}) also holds for the union of all 
cylindrical pieces $\sqcup^N S^1\times I$ and $M'$ replaced by $\phi_0(M)$.
The formula for $\be_1$ in Part 2 is now only a matter of counting
dimensions. \ep

For a cobordism $M:\Sigma_s\to\Sigma_t\,$ there are natural 
projections $p^t\,:\,H_1(-\Sigma_s\sqcup\Sigma_t)\to 
H_1(\Sigma_t)\,$, and $p^s$ analogously. 
For a composite of two cobordisms, $M_2\circ M_1$, we can define with
$\Sigma=\Sigma_{t,1}=\Sigma_{s,2}\,$ the following sub-spaces of $H_1(\Sigma)\,$:
\beq\lll{eq-def-V}
V_1\,:=\,\,p_1^t\Bigl(ker\bigl(H_1(\psi_1)\bigr)\Bigr)\;,
\qquad\qquad
V_2\,:=\,\,p_2^s\Bigl(ker\bigl(H_1(\psi_2)\bigr)\Bigr)\qquad.
\eeq  
From these we can find  $\mu_1$  as follows:
\blm\lll{lm-cocy=dim}
For two cobordisms, $M_1\,:\,\Sigma_{s,1}\to\Sigma\,$ and  
$M_2\,:\,\Sigma\to\Sigma_{t,2}\,$, and spaces  $V_j$ as above, 
the cocycles from (\ref{def-cocyc}) are given by:
\beq\lll{eq-cocy=dim-1}
\mu_1(M_2,M_1)\,-\,\mu_0(M_2,M_1)\;=\;codim\Bigl(V_1+ V_2\Bigr)\qquad.
\eeq
\elm

{\em Proof: } 
In order to compute $\mu_1$ let us also introduce the spaces
$\, ^t\H(M)\,:=\,coker\bigl(H_1(\psi^t)\bigr)\,$, and, correspondingly,
 $ ^s\H\,$.  We have the following exact sequence:
\beq\lll{eq-Ht-seq}
0\to V_1 \to H_1(\Sigma) \to\, ^s\H(M_1) \to \H (M_1) \to 0
\eeq
The fact that $V_1$ appears here as the kernel is seen, when we 
divide $0\to ker(\psi_1)\to H_1(\Sigma_{s,1}\sqcup\Sigma)\to H_1(M_1)$
by the sub-sequence through $H_1(\Sigma_{s,1})$. 

If we divide the Mayer-Vietoris sequence for the gluing of $M_2$ and $M_1$,
by the homologies of the outer surface $\Sigma_{s,1}\sqcup\Sigma_{t,2}\,$,
we obtain an exact sequence:
\beq\lll{eq-MV-cobcomp}
0\,\to\,V_1\cap V_2 \,\to\,H_1(\Sigma)\,\to\,\, ^s\H(M_1)\oplus\, ^t\H(M_2)\,
\to\,\H(M_2\circ M_1)\,\to\,{\bf G}^{\mu_0}\,\to\,0\qquad.
\eeq
Combining this with (\ref{eq-Ht-seq}) and its analogue for $V_2$, we find
the formula for $\mu_1\,$ by counting dimensions.\ep

Lemma~\ref{lm-cocy=dim} implies that not only  $\mu_0$, and $\mu_1$, but
 also $\mu_{\partial}:=\mu_1-\mu_0$ are
non-negative integers. This means that the interior Betti-numbers
$\be_j$, and also $\be_{\partial}:=\be_1-\be_0\,$
 can only increase under compositions. An increase in $\be_0(M)$ indicates
new 1-cycles in $\H(M)$ that are created by connecting two manifolds over several
boundary components as in (\ref{eq-mu0-elprod}). An increase in $\be_{\partial}\,$
indicates 1-cycles  from the intermediate surfaces $\Sigma$ that are not killed 
in the gluing process.  The simplest example for this are two copies of a 
full torus that are  glued together along the same meridians to give 
$S^1\times S^2$. We easily see that in this case $\mu_{\partial}=1\,$. 
\ms

\prg{A.3}{Summary of Tangle Presentations}

We shall summarize  here the tangle presentations of the cobordism category
$\Cb 0$ of cobordisms between parametrized, connected surfaces, with central 
extension by $\Omega_4$, as in  [Ke2]. There are several versions, and
we shall choose here the one that involves ordinary, framed tangles, 
as opposed to bridged tangles. The bridged versions are more practical
for the construction of invariants and connected TQFT's, as in [KL],
 but are less common, and involve another type of surgeries. 

The presentation is given by an isomorphism functor:
$$
\T\,:\;\Cb 0\,\isto\,\ct 0)\quad.
$$ 
On the level of objects we assign to a surface 
$\Sigma\,=\,\Sigma_{g_1}\sqcup\ldots\sqcup\Sigma_{g_K}\,$, with $K$ ordered,
connected components, a configuration of intervals, called {\em groups}, and
points on the real line  as follows. To every component 
$\Sigma_{g_j}\,$ corresponds to a group $G_j\subset {\bf R}\,$, and to 
$\Sigma$ we associate the disjoint, ordered union 
$G_1\sqcup\ldots\sqcup G_K\subset {\bf R}\,$.
Furthermore, in each interval $G_j$ we pick $2g_j$ points, which we denote
in increasing order $\{a_{j,1},a'_{j,1},a_{j,2},\ldots,a'_{j,g_j}\}\,
\subset\,G_j\,$. 

A connected morphism in $\Cb 0\,$,
$$
M\,:\;\Sigma_s=\Sigma_{g^s_1}\sqcup\ldots\sqcup\Sigma_{g^s_{K}}\,\longrightarrow\,
\Sigma_t=\Sigma_{g^t_1}\sqcup\ldots\sqcup\Sigma_{g^t_{L}}\quad,
$$
is then represented by a generically projected 
tangle in a strip ${\bf R}\times I\,$, where we have 
marked the groups and points of $\Sigma_s$ on the upper boundary, ${\bf R}_s\,$, 
of the 
strip, and the data of $\Sigma_t$ on the lower real line, ${\bf R}_t\,$. 
The tangle diagram
has thus a general form as below, where we indicated groups by braces.

$$
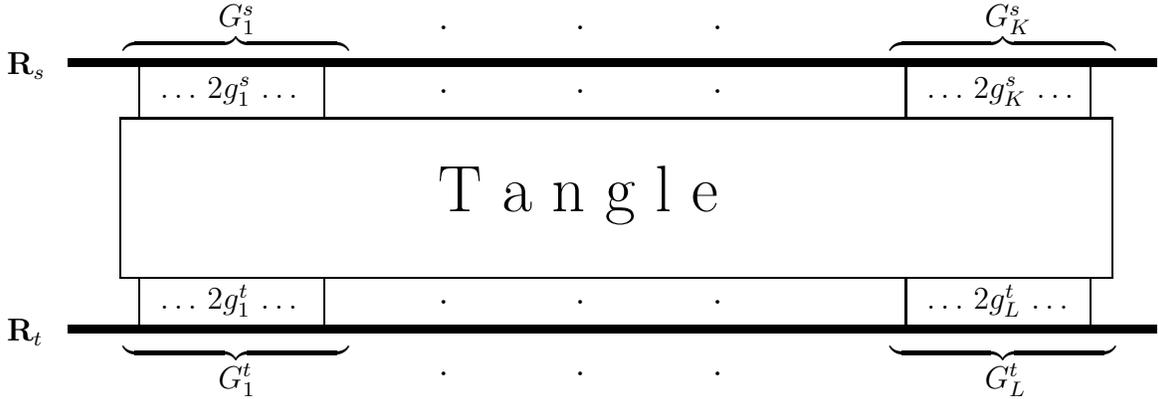
\begin{picture}(400,140)

\put(-13,16){${\bf R}_t$}

\put(-13,117){${\bf R}_s$}

\put(10,19){\rule{5.7in}{1mm}}

\put(10,120){\rule{5.7in}{1mm}}

\put(30,40){\framebox(375,60)}

\put(150,65){\Huge T\ a\ n\ g\ l\ e}

\put(31,15){$\underbrace{\hspace*{3cm}}$}
\put(67,-2){$G_1^t$}

\put(321,15){$\underbrace{\hspace*{3cm}}$}
\put(357,-2){$G_L^t$}

\put(150,3){\Large .\qquad\quad.\quad\qquad.}
\put(150,134){\Large .\qquad\quad.\quad\qquad.}
\put(150,30){\Large .\qquad\quad.\quad\qquad.}
\put(150,110){\Large .\qquad\quad.\quad\qquad.}

\put(31,125){$\overbrace{\hspace*{3cm}}$}
\put(67,135){$G_1^s$}

\put(321,125){$\overbrace{\hspace*{3cm}}$}
\put(357,135){$G_K^s$}

\put(37,100){\line(0,1){20}}
\put(107,100){\line(0,1){20}}
\put(45,108){$\dots\,2g_1^s\,\ldots$}

\put(327,100){\line(0,1){20}}
\put(397,100){\line(0,1){20}}
\put(335,108){$\dots\,2g_K^s\,\ldots$}

\put(37,20){\line(0,1){20}}
\put(107,20){\line(0,1){20}}
\put(45,28){$\dots\,2g_1^t\,\ldots$}

\put(327,20){\line(0,1){20}}
\put(397,20){\line(0,1){20}}
\put(335,28){$\dots\,2g_L^t\,\ldots$}

\end{picture}
$$

The tangle consists of possibly linked and knotted strands that end in the
$2(g_1^s+\ldots+g_K^s)$ points at ${\bf R}_s$, and the 
$2(g_1^t+\ldots+g_L^t)$ points at ${\bf R}_t\,$. To each strand we assign
in addition a framing of its normal bundle, which is the same as considering
ribbons. In diagrams we tacitly assume the framing to be in the plane of
projection (blackboard-framing).
A  component of a tangle in $\ct 0)$ has to be one of the following four
types:

\ben
\item {\em Source-ribbon}: A strand, that  connects a pair in ${\bf R}_s$, 
i.e., it starts at a point $a_{j,\nu}\in G^s_j$ and ends in $a_{j,\nu}'\in G^s_j$.
 \item {\em Target-ribbon}: A strand, that  connects a pair in ${\bf R}_s$, 
i.e., it starts at a point $a_{j,\nu}\in G^s_j$ and ends in $a_{j,\nu}'\in G^s_j$.

\item {\em Closed ribbon}: A strand 
($\cong S^1\times [-\varepsilon,\varepsilon]$)
 that is disjoint from the boundaries of
the diagram.

\item {\em Through-ribbon-pair}: A pair of strands, connecting a pair in ${\bf R}_s$
to a pair in ${\bf R}_t\,$. I.e., one strand connects $a_{j,\nu}\in G_j^s$ to
$a_{k,\mu}\in G_k^t$ and the other $a_{j,\nu}'\in G_j^s$ to
$a_{k,\mu}'\in G_k^t\,$. (Flipped version also admissible)
\een

Furthermore, we require that in all cases the ribbon, as a surface with boundary,
 has to be oriented, and the orientation shall be compatible with a fixed one
at ${\bf R}_s\,$ and ${\bf R}_t\,$, if it ends in the boundary of the strip.
\medskip

An element in $\ct 0)$ is then an equivalence class of tangles, where 
we consider two tangles equivalent, if one can be transformed into the 
other by the application of a sequence of the following five {\em moves}:

\ben
\item Isotopies with fixed endpoints. 

\item The ${\cal O}_2$-move or 2-handle slide of any type of  ribbon over 
    a closed ribbon, see [Ki].

\item The $\nhopf$-move, in which  we add or remove an isolated Hopf link,
for which one component has 0-framing, and the other either 1- or 0-framing.

\item The $\tau$-move (at a group $G^{s/t}_j\,$),
which allows us to push any type of ribbon through
the $2g_j^{s/t}$ parallel strands emerging at the group $G^{s/t}_j\,$ very close
to ${\bf R}_{s/t}\,$, as in the following diagram:
$$
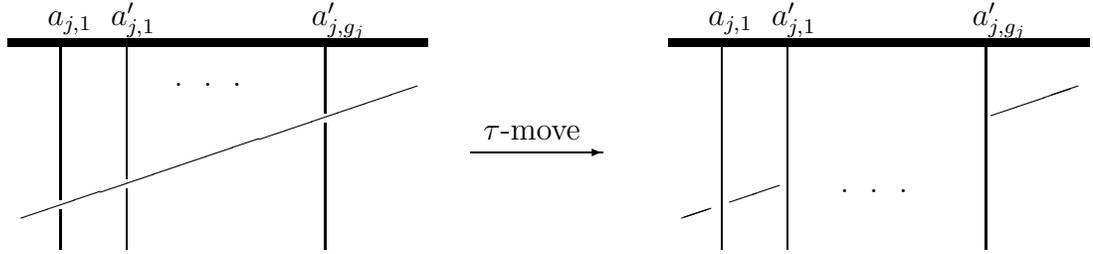
\begin{picture}(400,100)
\put(30,93){$a_{j,1}$}
\put(54,93){$a_{j,1}'$}
\put(35,85){\line(0,-1){58}}
\put(60,85){\line(0,-1){50}}

\put(35,23){\line(0,-1){15}}
\put(60,31){\line(0,-1){23}}

\put(78,70){$.\ \ .\ \ .$}

\put(130,93){$a_{j,g_j}'$}
\put(135,85){\line(0,-1){25}}
\put(135,56){\line(0,-1){48}}

\put(15,85){\rule{2.2in}{1mm}}
\put(20,20){\line(3,1){150}}

\put(190,45){\vector(1,0){50}}
\put(195,50){$\tau$-move}

\put(280,93){$a_{j,1}$}
\put(304,93){$a_{j,1}'$}
\put(285,85){\line(0,-1){77}}
\put(310,85){\line(0,-1){77}}

\put(330,30){$.\ \ .\ \ .$}

\put(380,93){$a_{j,g_j}'$}
\put(385,85){\line(0,-1){77}}

\put(265,85){\rule{2.2in}{1mm}}
\put(270,20){\line(3,1){11}}
\put(288,26){\line(3,1){19}}
\put(387,59){\line(3,1){33}}

\end{picture}
$$

\item The $\sigma$-move (at a pair $\{a_{j,\nu},a_{j,\nu}'\}\subset G_j^{s/t}\,$),
in which we replace (say for ${\bf R}_s$)
the two parallel strands at a pair by the chain of an
upward arc, a 0-framed
 annulus, and a downward arc as in the following diagram:

$$
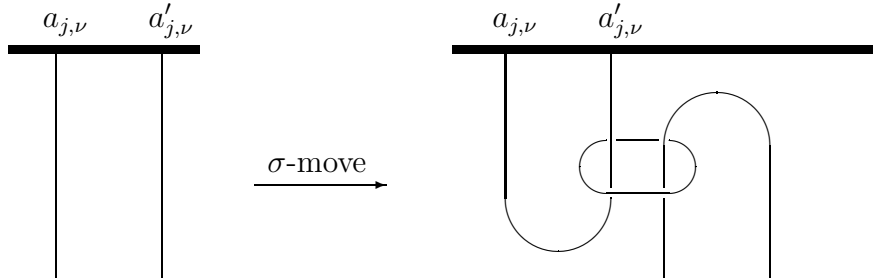
\begin{picture}(400,130)

\put(65,115){$a_{j,\nu}$}
\put(105,115){$a_{j,\nu}'$}
\put(70,105){\line(0,-1){85}}
\put(110,105){\line(0,-1){85}}

\put(52,105){\rule{1in}{1mm}}

\put(145,55){\vector(1,0){50}}
\put(150,60){$\sigma$-move}

\put(235,115){$a_{j,\nu}$}
\put(275,115){$a_{j,\nu}'$}

\put(220,105){\rule{2.2in}{1mm}}
\put(240,105){\line(0,-1){55}}
\put(280,105){\line(0,-1){51}}
\put(260,50){\oval(40,40)[b]}

\put(278,52){\line(1,0){24}}
\put(278,62){\oval(20,20)[l]}
\put(302,62){\oval(20,20)[r]}
\put(282,72){\line(1,0){16}}

\put(320,70){\oval(40,40)[t]}
\put(300,70){\line(0,-1){16}}
\put(340,70){\line(0,-1){50}}
\put(300,50){\line(0,-1){30}}

\end{picture}
$$

\een

In general the $\tau$-move has to be assumed for all groups both 
at ${\bf R}_s$
and ${\bf R}_t$. Yet, if we consider cobordisms in $\Cc 0\,$, 
where we have only one source and one target group, it is easy to see that
one $\tau$-move (say at the source) also implies the other.
We denote the respective sub-category of tangles with only one group at
top and bottom by $\ct 0)\cn\,$.
\medskip

It has been shown in [KL] that the cobordisms $\Cc n$ are represented by
the tangle category $\ct n)\cn$. The tangles in there contain in addition to
the top-, bottom-, closed, and through-Ribbons, $n$ so called {\em exterior}
 ribbons.
Each of there starts at a fixed point at ${\bf R}_s$ and ends at another 
fixed point in
${\bf R}_t$.  They can also be isotoped, and slid over closed ribbons.
However,
 there is no analog of a $\sigma$-move, and the groups, which are
 associated to one-to-one to the  
components of $\partial\phi_0(M)\,$, and  onto which we impose
the $\tau$-move,  include the external strands that represent the punctures
in the respective surface.

In the special case of $\ct 1)\cn$ the external ribbon can be isolated, and
the  class of a tangles is determined by an ordinary
 tangle without external strand, see [KL]. We  therefore introduce the
category, $\ctinf 1)\cn$, which has the same tangles as $\ct 0)\cn$ 
as generators,
but where we omitted the $\tau$-move, since we cannot isotop through the
extra strand. Although their presentations are different we have
$\ct 1)\cn\cong\ctinf 1)\cn\,$, 
which gives us the presentation of $\Cc 1$ we used in 
Section 4.3, as well as the relation in diagram (\ref{diag-1to0}).
For mapping class groups this coincides with the presentation in [MP]. 
\medskip

The composition law in $\ct 0)$, which makes $\T$ into a functor, has mostly
been described already 
 in the discussion given in Section 4.1. If we have two cobordisms,
$$
M\,:\;\Sigma_A\,\longrightarrow\,\Sigma_B\sqcup\Sigma_C\qquad\;{\rm and}\;
\qquad N\,:\;\Sigma_C\sqcup\Sigma_D\,\longrightarrow\,\Sigma_E\;,
$$
where $M$ and $N$ are connected, but the surfaces
$\Sigma_A\,$,  $\Sigma_B\,$ etc., may
be disconnected, then the  elementary composition  
$(\id\otimes N)(M\otimes\id):\Sigma_A\sqcup\Sigma_D\to
\Sigma_B\sqcup\Sigma_E\,$ is represented by the following tangle:
$$
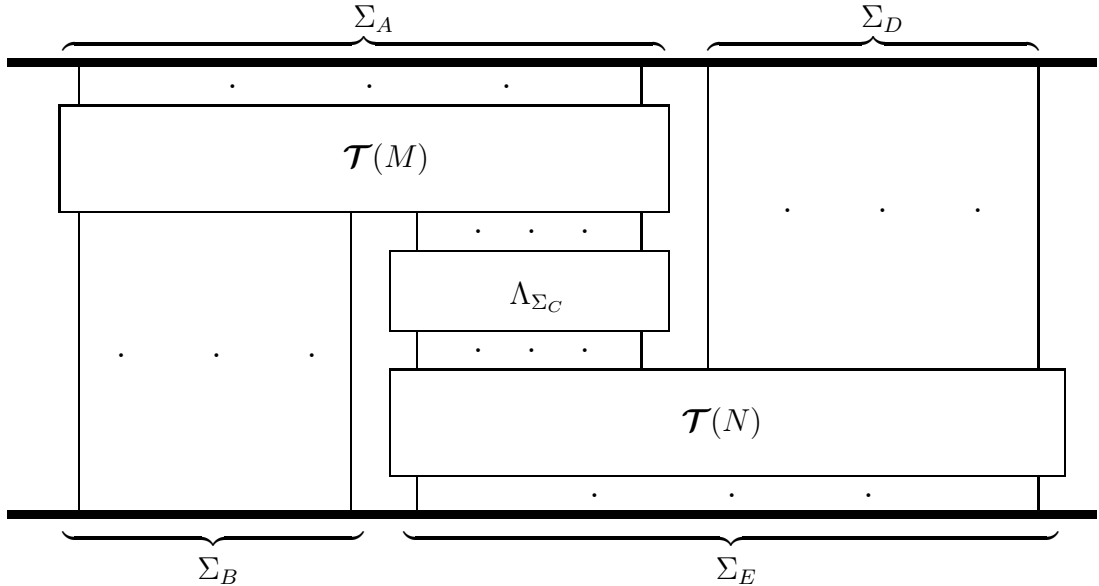
\begin{picture}(420,220)

\put(10,19){\rule{5.7in}{1mm}}

\put(31,15){$\underbrace{\hspace*{4cm}}$}
\put(82,-4){$\Sigma_B$}

\put(160,15){$\underbrace{\hspace*{8.7cm}}$}
\put(278,-4){$\Sigma_E$}

\put(230,27){\Large .\qquad\quad.\quad\qquad.}

\put(37,20){\line(0,1){115}}
\put(140,20){\line(0,1){115}}
\put(51,80){\Large .\qquad.\qquad.}

\put(165,20){\line(0,1){15}}
\put(400,20){\line(0,1){15}}

\put(155,35){\framebox(255,40)}

\put(265,52){$\T(N)$}

\put(165,75){\line(0,1){15}}
\put(250,75){\line(0,1){15}}
\put(186,82){\Large .\quad.\quad.}

\put(155,90){\framebox(105,30)}

\put(200,100){$\Lambda_{\Sigma_C}$}

\put(165,120){\line(0,1){15}}
\put(250,120){\line(0,1){15}}
\put(186,127){\Large .\quad.\quad.}

\put(30,135){\framebox(230,40)}

\put(137,152){$\T(M)$}

\put(37,175){\line(0,1){15}}
\put(250,175){\line(0,1){15}}

\put(10,190){\rule{5.7in}{1mm}}

\put(31,195){$\overbrace{\hspace*{8cm}}$}
\put(141,205){$\Sigma_A$}

\put(275,195){$\overbrace{\hspace*{4.4cm}}$}
\put(333,205){$\Sigma_D$}

\put(400,75){\line(0,1){115}}
\put(275,75){\line(0,1){115}}
\put(303,135){\Large .\qquad.\qquad.}

\put(37,175){\line(0,1){15}}
\put(250,175){\line(0,1){15}}
\put(93,182){\Large .\qquad\quad.\quad\qquad.}

\end{picture}
$$

Here $\T(M)$ and $\T(N)$ are tangles representing $M$ and $N$, and the tangle
$\Lambda_{\Sigma_C}\,$ is as in (\ref{fig-lambda-mor}) of Section 4.1.
The braces in this diagram enclose the union of the groups in the indicated
surface, i.e., we have for example $\beta_0(\Sigma_A)+\beta_0(\Sigma_D)$
groups at ${\bf R}_s\,$, with $\beta_1(\Sigma_A)+\beta_1(\Sigma_D)$ emerging
strands. 

The action of the symmetric group is as described in (\ref{fig-braid-Pi}).
The presentation of a disconnected cobordism is given simply by the
presentations of its connected components, if those are in a compatible
order.
\ms

The interior groups can also be computed from the diagrams, in an analogous
way as homotopy and homology groups are computed for link presentations of
closed manifolds. Using the explicit form of the functor $\T$, see [Ke2],
it is not  hard to see that, e.g.,  an $a$-cycle in $H_1(\Sigma_s)\,$
corresponds to a small meridian around the corresponding source-ribbon, and
the $b$-cycle of the same handle corresponds to a path that is pushed off the
same source-ribbon.

\section*{References}\lll{pg-R}
{\small

\ben

\item[[A]] Atiyah, M.: Publ. Math. Inst. Hautes Etudes Sci. {\bf 68}, 175-186 (1989).

\item[[D]] Deligne, P.: In ``The Grothendieck Festschrift''  {\bf II}, 
Progr. in Math., Birkh\"auser (1990).

\item[[D2]] Deligne, P.: Private communication. 

\item[[H]] Hennings, M.: Invariants of Links and 3-Manifolds obtained from Hopf Algebras,
 Preprint (1990).

Kauffman, L.,  Radford, D.: Invariants of 3-Manifolds Derived from Finite
Dimensional Hopf Algebras, Preprint {\tt hep-th/9406065}; 
to appear in J. Knot Theory and its Rami.

Ohtsuki, T.: see [O].

\item[[He]] Hempel, J.: ``3-Manifolds''. Annals of 
Mathematical Studies {\bf 86}. Princeton University Press, 1976.

\item[[KW]] Kac, V., Wakimoto, M.: Adv. Math. {\bf 40}, 1, 156-236 
(1988).

\item[[Ke1]]  Kerler, T.: Commun. Math. Phys. {\bf 168},  353-388 (1994).

\item[[Ke2]] Kerler, T.: Bridged Links and Tangle Presentations of
Cobordism Categories, Preprint (1994).

\item[[Ke3]] Kerler, T.: Genealogy of Nonperturbative Quantum-Invariants of 3-Manifolds: 
The Surgical Family.  To appear in Proceedings of the
Conference on Geometry and Physics, Aarhus, Denmark, 1995.
{\tt q-alg/9601021}.

\item[[Ke4]] Kerler, T.: Non-Tannakian Categories in Quantum Field Theory.
New Symmetry Principles in Quantum Field Theory, Cargese, 1991.
Physics  {\bf 295}, Plenum Press.

\item[[Ki]] Kirby, R.: Invent. Math. {\bf 45}, 35-56 (1978).

\item[[KL]] Kerler, T., Lyubashenko, V.: Non-Semisimple TQFT's for Connected
Surfaces, Preprint. 

\item[[Ku]] Kuperberg, G.: Non-Involutory Hopf algebras \& 3-Manifold Invariants, Preprint (1994).

\item[[L1]] Lyubashenko, V.: J. Pure Appl. Alg. {\bf 98}, 279-327 (1995).

\item[[L2]]  Lyubashenko, V.: Commun. Math. Phys.  {\bf 172}, 467-516 (1995).

\item[[L3]]  Lyubashenko, V.: Private communication.

\item[[LS]] Larson, R., Sweedler, M.: Amer. J. of Math., {\bf 91}, No 1., 
75-94 (1969). 

\item[[M]] MacLane, S.: ``Categories for the Working Mathematician". GTM. 
Springer-Verlag (1971).

\item[[MP]] Matveev, S.,  Polyak, M.: Commun. Math. Phys. {\bf 160}, 537 (1994).

\item[[O]] Ohtsuki, T.: Math. Proc. Cambridge Philos. Soc. {\bf 117}, 
259-273 (1995).

\item[[Os]] Ostrik, V.: Decomposition of the Adjoint Representation of the 
Small Quantum $sl_2$. {\tt q-alg/9512026}.

\item[[Ro]] Rotman, J.: 
``An Introduction to Algebraic Topology''. GTM. Springer-Verlag, 1993.
 
\item[[Rp]] Kerler, T.: Diploma thesis: Darstellungen der Quantengruppen, 1989.

\item[[RS]] Rozansky, L., Saleur, H.: Nucl. Phys. {\bf B389}, 461 (1993).

\item[[RT]] Reshetikhin, N., Turaev, V.: Invent. Math. {\bf 103}, 547-597 (1991).

\item[[Sw]] Sweedler, M.: Annals of Math. {\bf 89}, 323-335 (1969).

\item[[Sz]] Switzer, R.: ``Algebraic Topology - Homotopy and Homology''. 
Grundlehren. Springer-Verlag, 1975.

\item[[T]] Turaev, V.: ``Quantum Invariants of 3-Manifolds''. Walter de Gryter (1994). 
 
\item[[TV]] Turaev, V., Viro, O.: Topology {\bf 31}, 865-902 (1992).

\item[[Wi]] Witten, E.: Commun. Math. Phys. {\bf 121}, 351-399 (1989).

\een
}
\bigskip
{\small {\sc Institute for Advanced Study, Olden Lane,
Princeton,  NJ 08540, USA}

{\em E-mail-address:} kerler@math.ias.edu }

\end{document}